\documentclass{amsart}
\usepackage{graphicx}
\usepackage[all]{xy}
\usepackage{latexsym}
\usepackage{amssymb}
\usepackage{amsmath}
\usepackage{color}
\usepackage{xcolor}

\usepackage{fullpage}
\xyoption{arc}

 \newfam\cyrfam

  \font\tencyr=wncyr10

  \font\sevencyr=wncyr7

  \font\fivecyr=wncyr5

  \textfont\cyrfam=\tencyr \scriptfont\cyrfam=\sevencyr

    \scriptscriptfont\cyrfam=\fivecyr

  \newfam\cyifam

  \font\tencyi=wncyi10

  \font\sevencyi=wncyi7

  \font\fivecyi=wncyi5

  \textfont\cyifam=\tencyi \scriptfont\cyifam=\sevencyi

    \scriptscriptfont\cyifam=\fivecyi

  \def\id{{\mbox{1 \hskip -7pt 1}}}
\newcommand{\sgn}{{\mathit s  \mathit g\mathit  n}}
 \newcommand{\lon}{\longrightarrow}
 \newcommand{\bu}{\bullet}
 
 \newcommand{\rar}{\rightarrow}
 \newcommand{\hook}{\hookrightarrow}
 
\newcommand{\p}{{\partial}}
  \newcommand{\ot}{\otimes}
\newcommand{\Id}{{\mathrm{Id}}}

 \newcommand{\Ker}{{\mathsf K \mathsf e \mathsf r}\, }

 \newcommand{\Def}{\mathsf{Def}}
 \newcommand{\sip}{\smallskip}
 \newcommand{\bip}{\bigskip}
 \newcommand{\mip}{\vspace{2.5mm}}
 
\newcommand{\Graphs}{\mathcal{G}raphs}

\newcommand{\Holie}{\mathcal{H} \mathit{olie}}
\newcommand{\Lie}{\mathcal{L} \mathit{ie}}

\newcommand{\fGC}{\mathsf{fGC}}

\newcommand{\GC}{\mathsf{GC}}

\newcommand{\BV}{{\mathcal B}{\mathcal V}}
\newcommand{\co}{\mathit{co}}
\newcommand{\BVc}{{\mathcal B}{\mathcal V}^{\mathfrak{c}}}

\newcommand{\Grac}{\mathcal{G} ra^{\mathfrak{c}}}

\newcommand{\BVGcyc}{\mathcal{BVH}graphs}
\newcommand{\BVHgraphs}{\mathcal{BVH}graphs}
\newcommand{\wBVHgraphs}{\mathcal{\widehat{BV}H}graphs}

\newcommand{\ICH}{\mathsf{ICH}}
\newcommand{\wICH}{\widehat{\mathsf{ICH}}}
\newcommand{\HGrac}{\mathcal{H}gra}
\newcommand{\HGra}{\mathcal{H}gra}
\newcommand{\HGraphs}{\mathcal{H}graphs}

\newcommand{\ftt}{\mathfrak{t}^{\mathfrak{c}}}
\newcommand{\fcM}{f\underline{\cM}}

\newcommand{\BVGrac}{\mathcal{BVH}gra}

\newcommand{\HGCc}{\mathsf{hGC}}
\newcommand{\fHGCc}{\mathsf{fhGC}}
\newcommand{\fcHGCc}{\mathsf{fchGC}}

\newcommand{\ICG}{\mathsf{ICG}}
\newcommand{\FM}{\mathsf{FM}}
\newcommand{\FFM}{\mathsf{FFM}}

\newcommand{\wi}{\text{$\xy
(0,0)*+{_i}*\frm{o}
\endxy$}}

\newcommand{\wo}{\text{$\xy
(0,0)*+{_{0}}*\frm{o}
\endxy$}}
\newcommand{\wj}{\text{$\xy
(0,0)*+{_j}*\frm{o}
\endxy$}}
\newcommand{\wn}{\text{$\xy
(0,0)*+{_n}*\frm{o}
\endxy$}}
\newcommand{\wk}{\text{$\xy
(0,0)*+{_k}*\frm{o}
\endxy$}}

 \newcommand{\Z}{{\mathbb Z}}
 \newcommand{\bS}{{\mathbb S}}
 \renewcommand{\P}{{\mathbb P}}
 \newcommand{\C}{{\mathbb C}}
 \newcommand{\R}{{\mathbb R}}
 \newcommand{\N}{{\mathbb N}}
 \newcommand{\K}{{\mathbb K}}

\newcommand{\sF}{{\mathsf F}}

\newcommand{\sV}{{\mathsf V}}

 \newcommand{\Beq}{\begin{equation}}
 \newcommand{\Eeq}{\end{equation}}
 \newcommand{\Beqr}{\begin{eqnarray}}
 \newcommand{\Eeqr}{\end{eqnarray}}
 \newcommand{\Beqrn}{\begin{eqnarray*}}
 \newcommand{\Eeqrn}{\end{eqnarray*}}
 \newcommand{\Ba}{\begin{array}}
 \newcommand{\Ea}{\end{array}}
 \newcommand{\Bi}{\begin{itemize}}
 \newcommand{\Ei}{\end{itemize}}
 \newcommand{\Bc}{\begin{center}}
 \newcommand{\Ec}{\end{center}}
\newcommand{\fc}{{\mathfrak c}}

\newcommand{\fs}{{\mathfrak s}}

\newcommand{\ft}{{\mathfrak t}}

 \newcommand{\f}{{\mathcal O}}

 \newcommand{\cC}{{\mathcal C}}
 \newcommand{\caD}{{\mathcal D}}
 \newcommand{\cE}{{\mathcal E}}
 \newcommand{\cF}{{\mathcal F}}
 \newcommand{\cG}{{\mathcal G}}
 \newcommand{\caH}{{\mathcal H}}

 \newcommand{\cM}{{\mathcal M}}

 \newcommand{\al}{\alpha}
 
 \newcommand{\ga}{\gamma}
 
 \newcommand{\Ga}{\Gamma}
 \newcommand{\Ups}{\Upsilon}

 \newcommand{\om}{\omega}

\theoremstyle{plain}
\swapnumbers

\newtheorem{prop-def}[theorem]{Proposition-definition}

\newtheorem{main-theorem}{Main~Theorem}[section]
\newtheorem{section-theorem}{Theorem}[section]

\theoremstyle{definition}

\begin{document}

 \sloppy

 \newenvironment{proo}{\begin{trivlist} \item{\sc {Proof.}}}
  {\hfill $\square$ \end{trivlist}}

\long\def\symbolfootnote[#1]#2{\begingroup%
\def\thefootnote{\fnsymbol{footnote}}\footnote[#1]{#2}\endgroup}

%

\title{A hypergraph model for the cyclic BV operad\\ and its applications}

\author{Sergei A.\ Merkulov}
\address{Sergei~A.\ Merkulov,  Faculty of Mathematics,  The Higher School of Economics, Moscow }
\email{smerkulov25@gmail.com}

\begin{abstract}  A dg cyclic operad $\BVGcyc$ of hypergraphs is introduced which comes equipped
with an explicit quasi-isomorphism $\BVc\rar \BVGcyc$ from the { cyclic} operad $\BVc$ of Batalin-Vilkovisky algebras. A proof that the cohomology of $\BVGcyc$ equals $\BVc$ occupies most of this paper.


\sip

We use this model to construct an explicit quasi-isomorphism  $\cC hains(\FFM_2) \rar \BVGcyc$ from
the chain operad of the cyclic operad $\FFM_2$ of the compactified moduli spaces of genus zero curves
with marked framed points to  the dg cyclic operad $\BVGcyc$ which, combined with the main result mentioned above, gives a new proof of the {\it cyclic}\, formality 
of $\FFM_2$.

\sip






 \sip

%

\bip


\end{abstract}
 \maketitle
 \markboth{}{}


{\Large
\section{\bf Introduction}
}
\subsection{Operads of (framed) little disks}
The formality of the operad of little $d\geq 2$-dimensional disks $\caD_d$ 
was first proven for $d=2$ by Tamarkin in \cite{T} over a field $\K$ of characteristic zero using the theory of Kohno-Drinfeld Lie algebras. Later another proof was suggested by Kontsevich \cite{Ko} (see also Lambrechts-Voli\'{c} \cite{LV} for full details) over the field of real numbers via an explicit transcendental formula  built with the help of a beautiful graph model $\cG raphs_d$ for the homology operad $e_d:=H_\bu(\caD_d)$  of the little $d$-disks operad $\caD_d$, and the theory of semi-algebraic differential forms on semi-algebraic sets introduced by Kontsevich-Soibelman \cite{KS} and developed by  Hardt-Lambrechts-Turchin-Voli\'{c} \cite{HLTV}. The latter ``explicit formula" proof has also used heavily  the homotopy equivalence between the topological  operad $\caD_d$ and the topological operad $\FM_d$ of Fulton-MacPherson compactified configuration spaces in $\R^d$.

\sip

An elegant and independent proof of a key theorem behind Kontsevich's ``explicit formula" proof of the formality of $\caD_2$ can be found in \cite{SW};  Ševera-Willwacher's approach is based on a certain operad in the category of $\Lie_\infty$-algebras --- the so call operad ${\ICG}$ of { internally connected graphs} --- which gives us a combinatorial model \cite{Wi} for the operad of Drinfeld-Kohno Lie algebras (also known as Lie algebras of infinitesimal braids).

\sip

The topological operad $\caD_d$ has a framed version, the so called operad of  {\it framed}\, little $d$-disks $f\caD_d$ which plays an important role in homotopy algebra, algebraic topology and mathematical physics \cite{G,SW, S, CIW, Br, LW, Wi1,Wi2}. Its homology for $d=2$ equals the famous operad $\BV$ of Batalin-Vilkovisky algebras.

\sip

The operadic formality of  $f\caD_d$ for $d=2$ was established by
Giansiracusa-Salvatore [GS1] and Ševera [S]. It was proven by Khoroshkin-Willwacher \cite{KW} that for  $d\geq 3$ the operads $f\caD_d$ are formal for $d$ even and coformal for $d$ odd.

\sip

The framed little 2-discs operad $f\caD_2$ is homotopy equivalent to the Kimura-Stasheff-Voronov \cite{KSV} operad $\FFM_2$ of compactified moduli spaces of genus zero algebraic curves with at least two marked framed points. The latter topological operad has an obvious structure of a {\it cyclic}\, operad whose homology is the {\it cyclic}\ operad $\BVc$ of
Batalin-Volkovisky algebras. The results cited in the previous paragraph imply that the topological  operad $\FFM_2$ is formal as a non-cyclic operad.

\sip

The problem of formality of $\FFM_2$ as a {\it cyclic}\, operad over real numbers  was first attacked by Giansiracusa-Salvatore  [GS2] via an attempt to build a combinatorial model --- the so called complex of projective graphs --- of the operad $\BVc$, and then prove the formality of
$\FFM_2$ by an explicit transcendental formula along the lines of
 Kontsevich's approach \cite{Ko} to the formality of $\caD_2$. The work \cite{GS} is very nice and useful (e.g.\ we use below several results proven in \cite{GS} including a very important Vanishing Lemma), but it contains a mistake which was detected by Thomas Willwacher \cite{Wi4}.
  The main problem in \cite{GS} stems from their so called {\it pinwheel}\, condition on the projective graphs (see Definition 4.4 (R4) there) which 
 was imposed as an attempt to bypass the following obstacle:  given an arbitrary $BV$ algebra $(A,\Delta)$, a skew-symmetric binary operation\footnote{One can modify this binary operation into $\{a,b\}':= \Delta(\Delta(a)\cdot b) - (-1)^{|a||b|ab} - \Delta(\Delta(b)\cdot a) + \frac{(-1)^{|a|}}{2}\Delta(a)\Delta(b) $ which behaves much better under the action of the cyclic permutation $(012)$ in $\BVc$; this would change the differential  in the complex of projective graphs studied in  \cite{GS} but would not solve the pinwheel problem.} of degree $-2$
 $$
 \{a,b\}:= \Delta(\Delta(a)\cdot b) - (-1)^{|a||b|}\Delta(\Delta(b)\cdot a), \ \ \ \forall a,b\in A,
 $$
 does {\em not}\, in general satisfy the Jacobi identity.

 \sip
 
  The {\it cyclic}\, formality of 
$\FFM_2$ was first proven by Campos-Idrissi-Willwacher \cite{CIW} by extending Tamarkin's \cite{T} and  Ševera's [S] approaches to the (non-cyclic) formalities  of $\caD_2$ and $f\caD_2$ respectively.

\subsection{An outline of main results} 
In this paper we introduce  a dg cyclic operad $\BVGcyc$ of {hyper}graphs  and prove that it comes equipped with an explicit {\it quasi-isomorphism}\,
of dg cyclic operads
$$
F: \BVc\rar \BVGcyc.
$$ 
 A proof of this claim occupies most of this paper. The hypergraphs we use are of a relatively simple type: they have ordinary edges
(in cohomological degree $-1$) 
$$
\Ba{c}\resizebox{10mm}{!}{\xy
 (0,1)*{\circ}="a",
(7,1)*{\circ}="b",
\ar @{-} "a";"b"
\endxy}\Ea
$$
and trivalent hyperedges (in cohomological degree $-2$)
$$
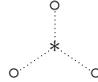
\Ba{c}\resizebox{13.5mm}{!}{
\xy
(0,2)*{\ast}="1";
(-6,-2)*{\circ}="2";
(0,8)*{\circ}="3";
(6,-2)*{\circ}="4";
\ar @{.} "1";"2" <0pt>
\ar @{.} "1";"3" <0pt>
\ar @{.} "1";"4" <0pt>
\endxy} \Ea
$$
only. Under cyclic operadic compositions degree -1 edges may combine into degree -2 hyperedges, so these two structures interact. 
\sip

Our second  result is a construction of an explicit quasi-isomorphism
$$
\Omega: \BVGcyc^* \lon \Omega^\bu_{\FFM_2}
$$
from the dual dg (Hopf) cooperad $\BVGcyc^*$ to the De Rham algebra of semi-algebraic differential forms on the geometric operad 
$\FFM_2$. This result implies in turn a quasi-isomorphism of dg cyclic operads
$$
Chains(\FFM_2) \lon \BVGcyc
$$
and hence proves the formality of $\FFM_2$ via an explicit zig-zag of quasi-isomorphisms of dg cyclic operads
$$
Chains(\FFM_2) \lon \BVGcyc \longleftarrow \BVc.
$$

Our approach to this main result goes through  a cyclic operad of internally connected hypergraphs $\ICH$ and its close relation to the cyclic operad of infinitesimal { ribbon} braids. It is much inspired by Appendix B in the Ševera-Willwacher paper \cite{SW} which gives a nice one-page proof of the quasi-isomorphism $e_2\rar \Graphs_2$.
Unfortunately, our  technical details are not that short, and occupy most of this paper.

\sip

\subsection{Why {\it hyper}graphs?} Let $\FFM_2((n+1))^0$ stand for the moduli spaces of genus zero algebraic curves with $(n+1)$ different framed points.
 For any $i \neq j\in \{0,1,\dots, n\}$ there is a map
$$
\pi_{ij}: \FFM_2((n+1))^0 \lon \FFM_2 \simeq S^1
$$
which forgets all framed points except the two labelled by $i$ and $j$. Let $\text{Vol}_{S^1}$ the standard homogeneous volume form on the circle $S^1$, and set
$$
\al_{ij}:=\pi_{ij}^*(\text{Vol}_{S^1}).
$$
These 1-forms extend nicely to the Kimura-Stasheff-Voronov \cite{KSV} compactification $\FFM_2((n+1))$ of $\FFM_2((n+1))^0$.
Consider also the following  ($S^1$-basic) 2-forms on $\FFM_2((n+1))$ with $n\geq 2$
\Beq\label{1: Omega_ijk}
\Omega_{ijk}:={\al}_{ij}{\al}_{jk} + {\al}_{ki}{\al}_{ij} + {\al}_{jk}{\al}_{ki},
\Eeq
and let $\bar{\al}_{ij}$ and
$$
\bar{\Omega}_{ijk}:=\bar{\al}_{ij}\bar{\al}_{jk} + \bar{\al}_{ki}\bar{\al}_{ij} + \bar{\al}_{jk}\bar{\al}_{ki}.
$$
stand for the associated cohomology classes in $H^\bu(\FFM_2((n+1)))$. There is an isomorphism of graded commutative algebras \cite{GS}
$$
H^\bu(\FFM_2((n+1)))=\frac{\R [\bar{\al}_{ij}]_{0\leq i\neq j\leq n}}{I},
$$
where $\R [\bar{\al}_{ij}]$ is the polynomial algebra generated by the degree 1 formal variables $\bar{\al}_{ij}$ and  $I$ is the ideal generated by the following 12-term {\em cyclic Arnold relations},
$$
\bar{\Omega}_{jkl} - \bar{\Omega}_{ikl} + \bar{\Omega}_{ijl} -\bar{\Omega}_{ijk} \equiv 0  \ \ \ \forall i,j,k,l\in \{0,1,\ldots,n\}, \
\# \{i,j,k,l\}=4.
$$

\sip

In our construction of a combinatorial model for $\BVc$ we treat the  
1-forms $\al_{ij}$ and the two forms $\Omega_{ijk}$ as independent of each other but encode their behavior under the ``cyclic co-compositions" (that is, their restrictions on  the codimension 1 boundary strata in $\FFM_2$) into the definition of the cyclic compositions in a certain cyclic operad $\BVGrac$ (see \S {\ref{4: Subsec on BVGrac}}  for its precise definition). That cyclic operad $\BVGrac$ is generated by hypergraphs which have  edges $
\Ba{c}\resizebox{10mm}{!}{\xy
 (0,1)*{\circ}="a",
(6,1)*{\circ}="b",
\ar @{-} "a";"b"
\endxy}\Ea
$ (corresponding to the 1-form propagators $\al_{ij}$) and hyperedges
$
\Ba{c}\resizebox{11.5mm}{!}{
\xy
(0,2)*{\ast}="1";
(-4,-1)*{\circ}="2";
(0,6)*{\circ}="3";
(4,-1)*{\circ}="4";
\ar @{.} "1";"2" <0pt>
\ar @{.} "1";"3" <0pt>
\ar @{.} "1";"4" <0pt>
\endxy} \Ea
$ (corresponding to the 2-form propagators $\Omega_{ijk}$);
the main point is that  cyclic operadic compositions can create in general hyperedges from edges. The toy cyclic operad $\BVGrac$ comes equipped with a morphism of cyclic operads\footnote{The superscript $3$ in the symbol $\Lie_3$ means that the Lie bracket generator of $\Lie_3$ is assigned the cohomological degree $1-3=-2$.}
$$
\Lie_3 \rar \Grac
$$
which sends the Lie generator of $\Lie_3$ precisely into the hypergraph with one hyperedge
 $$
 \Ba{c}\resizebox{15mm}{!}
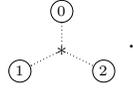
{ \xy
(-7,0)*+{_1}*\frm{o}="0";
 (7,0)*+{_2}*\frm{o}="1";
 (0,10)*+{_0}*\frm{o}="2";
  (0,3.5)*{*}="s";
 \ar @{.} "s";"2" <0pt>
 \ar @{.} "s";"0" <0pt>
 \ar @{.} "s";"1" <0pt>
\endxy} \Ea.
$$
Hence one can apply to $\BVGrac$ Thomas Willwacher's twisting endofunctor \cite{Wi} (adopted to cyclic operads)  to get --- after some extra twisting (in \S 4.6) reflecting the above mentioned interplay between edges and hyperdges --- a well-defined {\it differential}\, graded (dg for short) cyclic operad $f\BVHgraphs$ from which we extract the required combinatorial model $\BVHgraphs$  for
the cyclic operad $\BVc$ (see \S\S\,7-8).  A straightforward calculation (see \S {\ref{5: Ex 1}}) shows the relation (\ref{1: Omega_ijk}) between a priori independent propagators $\al_{ij}$ and $\Omega_{ijk}$ must hold true up to $d$-exact terms, and then another simple calculation (see \S {\ref{5: Ex 2}})
recovers the cyclic Arnold relations, again up to $d$-exact terms.

\sip

The dg operad $\BVHgraphs$ comes equipped with an obvious morphism of cyclic operads
$$
\BVc \lon \BVHgraphs
$$
which we prove in \S \S\,8-10 to be a   quasi-isomorphism. We also construct in \S\S\. 7-8 an explicit quasi-isomorphism 
of dg algebras
$$
\BVHgraphs^* \lon \Omega_{\FFM_2}
$$
from the dual (dg Hopf) cooperad  to the semialgebraic de Rham algebra on  $\Omega_{\FFM_2}$ which behaves nicely enough  to assure
that the associated morphism of cyclic operads $Chains(\FFM_2) \lon \BVGcyc$ is a well-defined quasi-isomorphism.

\subsection{Some notation}
 The set $\{1,2, \ldots, n\}$ is abbreviated to $[n]$ while the set $\{0,1,2,\ldots, n\}$ is abbreviated to $((n+1))$.
 The group of automorphisms of $[n]$ (resp., of $((n+1))$) is
denoted by $\bS_n$ (resp., $\bS_{(n+1)}$). The trivial (resp., sign) one-dimensional representation of
 $\bS_n$  is denoted by $\id_n$ (resp.,  $\sgn_n$).
 The cardinality of a finite set $A$ is denoted by $\# A$.

\sip

If $V=\oplus_{i\in \Z} V^i$ is a graded vector space, then
$V[n]$ stands for the graded vector space with $V[n]^i:=V^{i+n}$; for $v\in V^i$ we set $|v|:=i$.
 The canonical degree $-1$ isomorphism $V\rar V[1]$ is denoted by $\fs$;
one has $|\fs^k v|=|v| -k$ for any homogeneous element $v\in V$.






\mip

{\bf Acknowledgement}. It is a great pleasure to thank Vladimir Dotsenko, Anton Khoroshkin, Sergey Shadrin and especially Thomas Willwacher for many valuable discussions.

\bip


{\large
\section{\bf  Cyclic operads: definitions and examples}
}
\sip

\subsection{Cyclic (co)operads} A cyclic operad $\f^\fc$ in a symmetric monoidal category $\cC$ is a functor \cite{GK}
$$
\Ba{rccc}
\f^\fc: &  \mbox{Groupoid of finite finite sets} & \lon & \cC\\
& I &\lon & \f^\fc(I)
\Ea
$$
equipped, for any $i\in I$, $j\in J$, with a natural in $I$ and $J$ composition morphism
$$
_i\circ_j: \f^\fc(I)\ot \f^\fc(J) \lon \f^\fc(I\sqcup J\setminus \{i,j\})=: \f^\fc(I)\mbox{$_i\circ_j$} \f^\fc(J)
$$
satisfying the associativity axiom,
$$
\left(\f^\fc(I)\mbox{$_i\circ_{j_1}$} \f^\fc(J)\right) \mbox{$_{j_2}\circ_k$} \f^\fc(K)=
 \f^\fc(I)\mbox{$_i\circ_{j_1}$}\left( \f^\fc(J) \mbox{$_{j_2}\circ_k$} \f^\fc(K)\right),  \
 \forall i\in I, \ j_1\neq j_2 \in J,\  k\in K
$$
 and the commutativity axiom 
 which --- in the most important for us case when $\cC$ is a category of dg vector spaces --- reads as follows
 $$
 a\,\mbox{$_i\circ_{j}$} b = (-1)^{|a||b|} b\mbox{$_j\circ_{i}$} a, \ \ \forall\ a\in \f(I),\ b\in \f(J).
 $$
For any non-negative integer $n$ we set $\f((n+1)):=\f(\{0,1,\ldots , n\})$. Thus a cyclic operad
$\f$ can be understood as a countable collection $\{\f^\fc((n+1))\}_{n\geq 0}$ of $\bS_{(n+1)}$ modules equipped
with a countable collection of composition morphisms
$$
_i\circ_j: \f^\fc((n+1))\ot \f^\fc((n+1)) \lon \f^\fc((m+n)),
$$
 parameterized by pairs of integers $i\in \{0,1,\ldots, n\}$ and $j\in \{0,1, \ldots , m\}$,  and
 satisfying the associativity and commutativity axioms. All cyclic operads we consider in this
 paper have the ``curvature" component $\f(1)$ vanishing so we assume from now on that $n\geq 1$.

\sip

Dually, a cyclic cooperad $co \f^\fc$ is a functor
$$
\Ba{rccc}
\co\f^\fc: &  \mbox{Groupoid of finite finite sets} & \lon & \cC\\
& I &\lon & \co\f^\fc(I)
\Ea
$$
equipped, for any finite set $I$ and its decomposition into a disjoint union $I=I'\sqcup I''$, with a natural map
$$
\Delta_{I',I''}: \co\f^\fc(I) \lon \co\f^\fc(I'\sqcup x') \ot  \co\f^\fc(x''\sqcup I'')
$$
such that for any finite set $I$ and its decomposition $I=I'\sqcup I'' \sqcup I'''$ one has the equality
$
\Phi^L_{I',I'',I'''}=\Phi^R_{I',I'',I'''}
$
 of the following two  maps
\Beq\label{3: Phi^R}
\xymatrix{
 \Phi^R_{I',I'',I'''}: \co\f^\fc(I)  \ar[r]^-{^{\Delta_{I',I''\sqcup I'''}}}
 &
 \co\f^\fc(I'\sqcup x') \ot
 \co\f^\fc(x''\sqcup I'' \sqcup I''') \ar[d]^{^{1\ot \Delta_{x''\sqcup I'', I'''}}}\\
 &
  \co\f^\fc(I'\sqcup x') \ot
 \co\f^\fc(x''\sqcup I''\sqcup y')\ot \co\f^\fc(y''\sqcup  I''')
}
\Eeq
\Beq\label{3: Phi^L}
\xymatrix{
 \Phi^L_{I',I'',I'''}: \co\f^\fc(I)  \ar[r]^-{^{\Delta_{I'\sqcup I'', I'''}}}
 &
 \co\f^\fc(I'\sqcup I'' \sqcup y') \ot
 \co\f^\fc(y'' \sqcup I''') \ar[d]^{^{\Delta_{ I', I''\sqcup y'}}}\\
 &
  \co\f^\fc(I'\sqcup x') \ot
 \co\f^\fc(x''\sqcup I''\sqcup y')\ot \co\f^\fc(y''\sqcup  I''')
}
\Eeq

\subsubsection{\bf Derivations of cyclic operads} Let $\f^\fc=\{\f^\fc((n+1))\}_{n\geq 1}$ be a acyclic operad in
the category of dg vector spaces. Its derivation of degree $d\in \Z$ is, by definition, a collection of linear maps
$$
D: \f^\fc((n+1)) \rar \f^\fc((n+1)),\ \ \forall \ n\geq 1,
$$
such that
$$
D(a\mbox{$_i\circ_{j}$} b)= D(a)\mbox{$_i\circ_{j}$} b + (-1)^{d|a|} a\mbox{$_i\circ_{j}$} D(b).
$$
A differential in $\f$ is a derivation $\delta$ of degree $+1$ such that $\delta^2=0$.

\sip

Let $D$ be an element of $\f^\fc((2))$ satisfying the skew-symmetry condition under the action of $\bS_{(2)}$,
$$
(01)(D)=-D.
$$
Then the linear map given, for any $n\geq 1$, by
$$
\Ba{rccc}
D: & \f^\fc((n+1)) & \lon & \f^\fc((n+1)) \\
   &   a & \lon & D(a):= \sum_{i=0}^n D\mbox{$_1\circ_{i}$} a
\Ea
$$
is a derivation of degree $|D|$.

\subsection{\bf Example: a cyclic operad of commutative algebras} Let $\cC om^\fc=\{\cC om^\fc((n+1))\}_{n\geq 1}$
 be a collection of one-dimensional trivial representations of $\bS_{(n+1)}$,
$$
\cC om^\fc((n+1)):=\id_{n+1},
$$
with the composition maps $_i\circ_{j}$ being the standard isomorphisms $\K\ot \K\rar \K$.
This is a cyclic operad. It is suitable to represent pictorially a generator of $\bS_{n+1}$
as a corolla with $n+1$ legs (or hairs) attached which are labelled,
$$
\cC om^\fc((n+1))= \K \left\langle
 \Ba{c}\resizebox{11mm}{!}{\xy
 (0,0)*{\circ}="a",
(2.3,5)*{^{1}}="1",
(-2.3,5)*{^0}="2",
(5,1)*{^{\, 2}}="3",
(-5,1)*{^n}="4",
(-3.6,-3.6)*{}="5",
(3.6,-3.6)*{}="6",
(0,-4.5)*{}="7",
\ar @{-} "a";"1" <0pt>
\ar @{-} "a";"2" <0pt>
\ar @{-} "a";"3" <0pt>
\ar @{-} "a";"4" <0pt>
\ar @{-} "a";"5" <0pt>
\ar @{-} "a";"6" <0pt>
\ar @{-} "a";"7" <0pt>
\endxy}
\Ea
\right\rangle
$$
It is assumed that this corolla is non-planar and legs are undistinguishable  (so that $\bS_{(n+1)}$
acts trivially via relabeling). Then the operadic composition $_i\circ_{j}$ can be understood
pictorially as gluing $i$-labelled leg of one corolla to $j$-th labelled leg of another corolla
and then contracting the resulting internal edge (and a sutable relabelling of legs if necessary).

\sip

We denote by $\overline{\cC om^\fc}$  the {\em augmentation}\,  of $\cC om^\fc$ which is, by definition,  a
cyclic sub-operad of
$\cC om^\fc$ which has $\overline{\cC om^\fc}(2)=0$ and $\overline{\cC om^\fc}(n+1)=\cC om(n+1)$
for all $n\geq 2$.

\subsection{\bf Example: cyclic homology operad $H_\bu(S^1)$} Let $S^1$ be the topological circle,
and $H_\bu(S^1)=\K \oplus \K[1]$ be its homology group. The degree $-1$ generator
of $\K[1]$ is denoted by $\Delta$. This cohomology group  can be made into a cyclic operad $H_{\bu}(S^1)=\{ H_\bu(S^1)((n+1))\}$
 with all $\bS_{(n+1)}$-modules  $H_\bu(S^1)((n+1))$ equal to zero except the following one
$n=1$,
$$
H_\bu(S^1)(2)=\K\oplus \K[1],
$$
on which $\bS_2$ acts trivially. The cyclic composition rules are given by the isomorphisms
$$
\K \ot \K \rar \K, \ \  \K\ot \K[1] \rar \K[1], \ \ \ \K[1]\ot \K \rar \K[1],
$$
and the zero map $\K[1]\ot \K[1]\rar 0$.

\sip

Similarly to $\overline{\cC om^\fc}$ above we denote by $\bar{H}_\bu(S^1)$ the augmentation of
 $H_{\bu}(S^1)$; this is a cyclic sub-operad of $H_{\bu}(S^1)$ with $\bar{H}_\bu(S^1)((2))=\K[1]$ whose all operadic
 compositions are trivial. In pictures the generator of $\bar{H}_\bu(S^1)((2))$ is represented by an edge
 $\Ba{c}\resizebox{12mm}{!}{\xy
 (-4,1)*+{_0}="a",
(4,1)*+{_1}="b",
(0,0.7)*{_{^\Delta}},
\ar @{-} "a";"b" <0pt>
\endxy
}\Ea$ decorated with the symbol $\Delta$.

\subsection{\bf Example: a cyclic operad of $BV$-algebras} Let  $\BV'$ be the cyclic operad generated
freely by the cyclic operads $\overline{\cC om}$ and  $\bar{H}_{\bu}(S^1)$. Its elements can be identified with
 trees $T$ whose vertices are corollas from $\overline{\cC om}$
and whose
\Bi
\item[--]
  internal edges (if any) are all decorated with $\Delta$,
  \item[--] each labelled leg of $T$
may have none or precisely one $\Delta$ decoration, e.g.
\Ei
$$
 \Ba{c}\resizebox{26mm}{!}{\xy
 (0,0)*{\circ}="a",
(2.3,4)*{_{^0}}="1",
(-2.3,4)*{_{^1}}="2",
(5,1)*{}="3",
(-6,1)*{_{^3}}="4",
(-3.6,-3.6)*{_{^7}}="5",
(3.6,-3.6)*{_{^8}}="6",
(0,-4.5)*{_{_{10}}}="7",
 (9,2)*{\circ}="b",
(10,6)*{_{^4}}="1b",
(16,1)*{_{^9}}="2b",
(12.7,1.1)*{^{_\Delta}},
(4.7,0.6)*{^{_\Delta}},
(10.1,-2.4)*{^{_\Delta}},
(-3,0.1)*{^{_\Delta}},
 (11,-6)*{\circ}="c",
(15,-4)*{_{^6}}="1c",
(7.3,-8)*{_{^2}}="2c",
(13,-10)*{_{_5}}="3c",
\ar @{-} "a";"1" <0pt>
\ar @{-} "a";"2" <0pt>
\ar @{-} "a";"4" <0pt>
\ar @{-} "a";"5" <0pt>
\ar @{-} "a";"6" <0pt>
\ar @{-} "a";"7" <0pt>
\ar @{-} "a";"b" <0pt>
\ar @{-} "c";"b" <0pt>
\ar @{-} "1b";"b" <0pt>
\ar @{-} "2b";"b" <0pt>
\ar @{-} "1c";"c" <0pt>
\ar @{-} "2c";"c" <0pt>
\ar @{-} "3c";"c" <0pt>
\endxy}
\Ea\in \BV'((10+1))
$$
The operadic composition \mbox{$_i\circ_{j}$} of such trees is
\Bi
\item[(i)] given as in $\cC om$ if legs $i$ and $j$ are both not decorated by $\Delta$,
\item[(ii)] given by gluing legs $i$ and $j$ into a new internal edge if precisely one of
these legs is decorated by $\Delta$,
   \item[(iii)] set to be zero of both legs $i$ and $j$ are decorated by $\Delta$.
\Ei
The cyclic operad  $\BVc=\{\BVc((n+1))\}_{n\geq 1}$ of Batalin-Vilkovisky algebras  is
defined to be the quotient of $\BV'$ by the operadic ideal generated by
the following relation
\Beq\label{2: BV relation}
 \Ba{c}\resizebox{12mm}{!}{\xy
 (0,0)*{\circ}="a",
(2.3,5)*{_{^1}}="2",
(-5,2)*{_{^0}}="1",
(-2.3,-5)*{_{^4}}="4",
(5,-2)*{_{^3}}="3",
(-2.3,0.66)*{^{_\Delta}},
\ar @{-} "a";"1" <0pt>
\ar @{-} "a";"2" <0pt>
\ar @{-} "a";"3" <0pt>
\ar @{-} "a";"4" <0pt>
\endxy}
\Ea
+
 \Ba{c}\resizebox{12mm}{!}{\xy
 (0,0)*{\circ}="a",
(-5,2)*{_{^0}}="1",
(2.3,5.2)*{_{^1}}="2",
(5,-2)*{_{^3}}="3",
(-2.3,-5)*{_{^4}}="4",
(1.28,2.)*{^{_{\Delta}}},
\ar @{-} "a";"1" <0pt>
\ar @{-} "a";"2" <0pt>
\ar @{-} "a";"3" <0pt>
\ar @{-} "a";"4" <0pt>
\endxy}
\Ea
+
 \Ba{c}\resizebox{12mm}{!}{\xy
 (0,0)*{\circ}="a",
(-5,2)*{_{^0}}="1",
(2.3,5.2)*{_{^1}}="2",
(5,-2)*{_{^3}}="3",
(-2.3,-5)*{_{^4}}="4",
(2.5,-1.2)*{^{_{\Delta}}},
\ar @{-} "a";"1" <0pt>
\ar @{-} "a";"2" <0pt>
\ar @{-} "a";"3" <0pt>
\ar @{-} "a";"4" <0pt>
\endxy}
\Ea
+
 \Ba{c}\resizebox{12mm}{!}{\xy
 (0,0)*{\circ}="a",
(-5,2)*{_{^0}}="1",
(2.3,5.2)*{_{^1}}="2",
(5,-2)*{_{^3}}="3",
(-2.3,-5.2)*{_{^4}}="4",
(-0.8,-2.5)*{^{_\Delta}},
\ar @{-} "a";"1" <0pt>
\ar @{-} "a";"2" <0pt>
\ar @{-} "a";"3" <0pt>
\ar @{-} "a";"4" <0pt>
\endxy}
\Ea
=
 \Ba{c}\resizebox{14mm}{!}{\xy
 (9,1)*{\circ}="b",
(4,3)*{_{^0}}="1b",
(13,5)*{_{^1}}="2b",
(10.2,-2.4)*{^{_\Delta}},
 (11,-5)*{\circ}="c",
(16,-8.5)*{_{^2}}="1c",
(8,-10)*{_{^3}}="2c",
\ar @{-} "c";"b" <0pt>
\ar @{-} "1b";"b" <0pt>
\ar @{-} "2b";"b" <0pt>
\ar @{-} "1c";"c" <0pt>
\ar @{-} "2c";"c" <0pt>
\endxy}
\Ea
+
 \Ba{c}\resizebox{14mm}{!}{\xy
 (9,1)*{\circ}="b",
(4,3)*{_{^0}}="1b",
(13,5)*{_{^2}}="2b",
(10.2,-2.4)*{^{_\Delta}},
 (11,-5)*{\circ}="c",
(16,-8.5)*{_{^1}}="1c",
(8,-10)*{_{^4}}="2c",
\ar @{-} "c";"b" <0pt>
\ar @{-} "1b";"b" <0pt>
\ar @{-} "2b";"b" <0pt>
\ar @{-} "1c";"c" <0pt>
\ar @{-} "2c";"c" <0pt>
\endxy}
\Ea
+
 \Ba{c}\resizebox{14mm}{!}{\xy
 (9,1)*{\circ}="b",
(4,3)*{_{^0}}="1b",
(13,5)*{_{^3}}="2b",
(10.2,-2.4)*{^{_\Delta}},
 (11,-5)*{\circ}="c",
(16,-8.5)*{_{^1}}="1c",
(8,-10)*{_{^2}}="2c",
\ar @{-} "c";"b" <0pt>
\ar @{-} "1b";"b" <0pt>
\ar @{-} "2b";"b" <0pt>
\ar @{-} "1c";"c" <0pt>
\ar @{-} "2c";"c" <0pt>
\endxy}
\Ea.
\Eeq
The cyclic operad $\BV$ has essentially two generators
$$
\Ba{c}\resizebox{9mm}{!}{\xy
 (0,1)*{_0}="a",
(7,1)*{_1}="b",
(3.5,0.6)*{^{_\Delta}},
\ar @{-} "a";"b" <0pt>
\endxy
}\Ea\in \BVc(2), \ \ \
\Ba{c}\begin{xy}
 <0mm,0.66mm>*{};<0mm,3mm>*{}**@{-},
 <0.39mm,-0.39mm>*{};<2.2mm,-2.2mm>*{}**@{-},
 <-0.35mm,-0.35mm>*{};<-2.2mm,-2.2mm>*{}**@{-},
 <0mm,0mm>*{\circ};<0mm,0mm>*{}**@{},
   <0mm,0.66mm>*{};<0mm,3.4mm>*{^{^0}}**@{},
   <0.39mm,-0.39mm>*{};<2.9mm,-4mm>*{^{_2}}**@{},
   <-0.35mm,-0.35mm>*{};<-2.8mm,-4mm>*{^{_1}}**@{},
\end{xy}\Ea \in \overline{\cC om}((3))\subset \BVc((3)),
$$
subject to the above relation as well as the following one
$$
\Ba{c}\resizebox{9mm}{!}{\xy
 (0,1)*{_0}="a",
(7,1)*{_1}="b",
(3.5,0.6)*{^{_\Delta}},
\ar @{-} "a";"b" <0pt>
\endxy}
\Ea
\mbox{$_1\circ_0$}
\Ba{c}\resizebox{9mm}{!}{\xy
 (0,1)*{_0}="a",
(7,1)*{_1}="b",
(3.5,0.6)*{^{_\Delta}},
\ar @{-} "a";"b" <0pt>
\endxy}
\Ea
=0
$$
which encodes $\Delta^2=0$.


%

\subsection{\bf Example: a topological cyclic operad $\FFM_2$} A nice example of a cyclic operad in the category of the topological spaces is given the collection  $\FFM_2=\{\FFM_2(n+1)\}_{n\geq 1}$  of compactified
moduli spaces of genus zero algebraic curves with framed points introduced in \cite{KSV}; it has been shown in that paper that is homology operad $H_\bu(\FFM_2)$ is precisely the cyclic operad $\BVc$.

\subsection{\bf A cyclic operad of homotopy Lie algebras $\Holie^\fc_{d}$ for $d$ odd} The subscript $d$ means that the cyclic operad $\Holie^\fc_{d}=\{Holie^\fc_d((n+1))\}_{n\geq 2}$ is a minimal resolution of a cyclic operad of Lie algebras $\Lie_d$ whose Lie bracket generator has degree $1-d$; thus the case $d=1$  corresponds to  the ordinary
cyclic $\Lie_\infty$ operad. 

\sip

It is suitable to identify the generators of the cyclic operad
$\Holie_d^\fc((n+1))$ with the unrooted trees whose legs are labelled by integers $0,1,\ldots, n$ and
 whose internal edges (shown in pictures as dotted ones) are equipped with a direction (up to a flip and multiplication by $-1$), e.g.
$$
 \Ba{c}\resizebox{26mm}{!}{\xy
 (0,0)*{\bu}="a",
(2.3,4)*{_{^0}}="1",
(-2.3,4)*{_{^1}}="2",
(5,1)*{}="3",
(-5,1)*{_{^3}}="4",
(-3.6,-3.6)*{_{^7}}="5",
(3.6,-3.6)*{_{^8}}="6",
%
 (9,2)*{\bu}="b",
(10,6)*{_{^4}}="1b",
(16,1)*{_{^9}}="2b",
 (11,-6)*{\bu}="c",
(15,-4)*{_{^6}}="1c",
(7.3,-8)*{_{^2}}="2c",
(13,-10)*{_{_5}}="3c",
\ar @{-} "a";"1" <0pt>
\ar @{-} "a";"2" <0pt>
\ar @{-} "a";"4" <0pt>
\ar @{-} "a";"5" <0pt>
\ar @{-} "a";"6" <0pt>
%
\ar @{.>} "a";"b" <0pt>
\ar @{.>} "c";"b" <0pt>
\ar @{-} "1b";"b" <0pt>
\ar @{-} "2b";"b" <0pt>
\ar @{-} "1c";"c" <0pt>
\ar @{-} "2c";"c" <0pt>
\ar @{-} "3c";"c" <0pt>
\endxy}
\Ea
=
- \Ba{c}\resizebox{26mm}{!}{\xy
 (0,0)*{\bu}="a",
(2.3,4)*{_{^0}}="1",
(-2.3,4)*{_{^1}}="2",
(5,1)*{}="3",
(-5,1)*{_{^3}}="4",
(-3.6,-3.6)*{_{^7}}="5",
(3.6,-3.6)*{_{^8}}="6",
%
 (9,2)*{\bu}="b",
(10,6)*{_{^4}}="1b",
(16,1)*{_{^9}}="2b",
 (11,-6)*{\bu}="c",
(15,-4)*{_{^6}}="1c",
(7.3,-8)*{_{^2}}="2c",
(13,-10)*{_{_5}}="3c",
\ar @{-} "a";"1" <0pt>
\ar @{-} "a";"2" <0pt>
\ar @{-} "a";"4" <0pt>
\ar @{-} "a";"5" <0pt>
\ar @{-} "a";"6" <0pt>
%
\ar @{<.} "a";"b" <0pt>
\ar @{.>} "c";"b" <0pt>
\ar @{-} "1b";"b" <0pt>
\ar @{-} "2b";"b" <0pt>
\ar @{-} "1c";"c" <0pt>
\ar @{-} "2c";"c" <0pt>
\ar @{-} "3c";"c" <0pt>
\endxy}
\Ea
\in \Holie_d((9+1)).
$$
Each such unrooted tree comes equipped with an {\em orientation}\, which is, by definition,  an ordering of its
vertices and legs (up to a permutation and multiplication by its sign), and a choice of  a particular direction on each internal (dotted) edge. It is useful to
understand such a graph as composed from  corollas
\Beq\label{3: generators of Holie cyc}
  \Ba{c}\resizebox{13mm}{!}{\xy
 (0,0)*{\bu}="a",
(2.3,5)*{^{1}}="1",
(-2.3,5)*{^0}="2",
(5,1)*{^{\, 2}}="3",
(-5,1)*{^n}="4",
(-3.6,-3.6)*{}="5",
(3.6,-3.6)*{}="6",
(0,-4.5)*{}="7",
\ar @{-} "a";"1" <0pt>
\ar @{-} "a";"2" <0pt>
\ar @{-} "a";"3" <0pt>
\ar @{-} "a";"4" <0pt>
\ar @{-} "a";"5" <0pt>
\ar @{-} "a";"6" <0pt>
\ar @{-} "a";"7" <0pt>
\endxy}
\Ea
\Eeq
 whose black vertex has degree $1+2d$ while the attached legs (or hairs) have degree $-d$ as it explains the
 aforementioned orientation immediately; the internal edges have degree $-2d$ as they are assumed to be composed from two legs  of the associated vertices; these two legs must be ordered in one or another way which leads to a choice of a direction on each edge (up to the flip and multiplication by $-1$). We prefer showing the internal edges  in ``dotted color" to emphasize their different grading from the legs (which are shown in ``solid color").

 \sip

 The cyclic composition $\Ga'  \mbox{$_i\circ_j$} \Ga''$
of such unrooted trees is given by gluing $i$-th labelled leg (of degree $-d$) of $\Ga'$ to the $j$-th labelled leg (of degree $-d$)
of $\Ga''$ and creating thereby a degree $-2d$ internal edge, e.g.
\Beq\label{3: comp in cyc Holieb}
  \Ba{c}\resizebox{11mm}{!}{\xy
 (0,0)*{\bu}="a",
(2.3,5)*{}="1",
(-2.3,5)*{}="2",
(5,1)*{^{\, i}}="3",
(-5,1)*{}="4",
(-3.6,-3.6)*{}="5",
(3.6,-3.6)*{}="6",
(0,-4.5)*{}="7",
\ar @{-} "a";"1" <0pt>
\ar @{-} "a";"2" <0pt>
\ar @{-} "a";"3" <0pt>
\ar @{-} "a";"4" <0pt>
\ar @{-} "a";"5" <0pt>
\ar @{-} "a";"6" <0pt>
\ar @{-} "a";"7" <0pt>
\endxy}
\Ea
\mbox{$_i\circ_{j}$}
  \Ba{c}\resizebox{11mm}{!}{\xy
 (0,0)*{\bu}="a",
(2.3,5)*{}="1",
(-2.3,5)*{}="2",
(5,1)*{}="3",
(-5,1)*{^j}="4",
(-3.6,-3.6)*{}="5",
(3.6,-3.6)*{}="6",
(0,-4.5)*{}="7",
\ar @{-} "a";"1" <0pt>
\ar @{-} "a";"2" <0pt>
\ar @{-} "a";"4" <0pt>
\ar @{-} "a";"5" <0pt>
\ar @{-} "a";"6" <0pt>
\endxy}
\Ea
=
  \Ba{c}\resizebox{18mm}{!}{\xy
 (0,0)*{\bu}="a",
(2.3,5)*{}="1",
(-2.3,5)*{}="2",
(-5,1)*{}="4",
(-3.6,-3.6)*{}="5",
(3.6,-3.6)*{}="6",
(0,-4.5)*{}="7",
 (10,0)*{\bu}="a'",
(12.3,5)*{}="1'",
(7.7,5)*{}="2'",
(15,1)*{}="3'",
(7.4,-3.6)*{}="5'",
(13.6,-3.6)*{}="6'",
(10,-4.5)*{}="7'",
\ar @{-} "a";"1" <0pt>
\ar @{-} "a";"2" <0pt>
\ar @{-} "a";"4" <0pt>
\ar @{-} "a";"5" <0pt>
\ar @{-} "a";"6" <0pt>
\ar @{-} "a";"7" <0pt>
\ar @{.>} "a";"a'" <0pt>
\ar @{-} "a'";"1'" <0pt>
\ar @{-} "a'";"2'" <0pt>
\ar @{-} "a'";"5'" <0pt>
\ar @{-} "a'";"6'" <0pt>
\endxy}
\Ea
\Eeq
 The differential on unrooted trees $T\in \Holie_d$ is given 
 by the sum
 $$
 \delta T =\sum_{v\in V(T)} \delta_v T
 $$ 
  where $\delta_v$ acts on the vertex $v=\bu$ by splitting it as follows,
 $$
\delta_v: \bu \ \ \text{goes into} \  \Ba{c}\resizebox{12mm}{!}{ \xy
  (4,1)*{\bu}="r";
  (-4,1)*{\bu}="l";
 \ar @{.>} "l";"r" <0pt>
\endxy}\Ea.
 $$
 and then taking the sum over all possible reattachments of dangling half-edges/legs --- the ones which were attached originally to $v$ --- among the two newly created vertices in such a way that each new vertex  has valency at least 3.
 For example, the action of $\delta$ on a one-vertex graph looks as follows,
$$
\delta:
 \Ba{c}\resizebox{9mm}{!}{\xy
 (0,0)*{\bullet}="a",
(2.3,4)*{}="1",
(-2.3,4)*{}="2",
(5,1)*{}="3",
(-5,1)*{}="4",
(-3.6,-3.6)*{}="5",
(3.6,-3.6)*{}="6",
(0,-4.5)*{}="7",
\ar @{-} "a";"1" <0pt>
\ar @{-} "a";"2" <0pt>
\ar @{-} "a";"3" <0pt>
\ar @{-} "a";"4" <0pt>
\ar @{-} "a";"5" <0pt>
\ar @{-} "a";"6" <0pt>
\ar @{-} "a";"7" <0pt>
\endxy}
\Ea
\ \lon \ \sum
\Ba{c}\resizebox{10mm}{!}{\xy
%
%
 (0,-3.3)*{\bullet}="a", 
 (0,3.3)*{\bullet}="b",
(-7,-6)*{}="1",
(7,-7)*{}="2",
(-3,-9)*{}="3",
(3,-9)*{}="4",
(5,8)*{}="5",
(-5,8)*{}="6",
(0,9)*{}="7",
\ar @{.>} "a";"b" <0pt>
\ar @{-} "a";"1" <0pt>
\ar @{-} "a";"2" <0pt>
\ar @{-} "a";"3" <0pt>
\ar @{-} "a";"4" <0pt>
\ar @{-} "b";"5" <0pt>
\ar @{-} "b";"6" <0pt>
\ar @{-} "b";"7" <0pt>
\endxy}
\Ea.
$$
The cohomology of this cyclic dg operad $\Holie_d$, $d\in 2\Z+1$, is the operad  $\Lie_d^\fc$ of Lie algebras which is generated by trivalent graphs as in the picture
$$
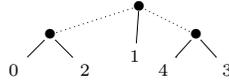
\Ba{c}\resizebox{33mm}{!}{ \xy
(-10,1)*+{_0}="0";
(0,1)*+{_2}="1";
(7,3)*+{_1}="2";
(11,1)*+{_4}="2'";
(20,1)*+{_3}="n";
  (-5,6)*{\bu}="bu";
   (7.5,10)*{\bu}="r";
     (15.5,6)*{\bu}="rr";
 \ar @{-} "0";"bu" <0pt>
  \ar @{-} "1";"bu" <0pt>
   \ar @{.} "bu";"r" <0pt>
    \ar @{-} "rr";"2'" <0pt>
     \ar @{-} "rr";"n" <0pt>
  \ar @{.} "rr";"r" <0pt>
 \ar @{-} "r";"2" <0pt>
\endxy} \Ea 
$$
defined modulo the so called IHX-relation for each dotted edge, that is, modulo the 3-term relation obtained by applying $\delta$  to a black vertex $v$ of valency 4,
$$
0= \left(\delta
 \Ba{c}\resizebox{9mm}{!}{\xy
 (0,0)*{\bullet}="a",
(2.3,4)*{}="1",
(-2.3,4)*{}="2",
(5,1)*{}="3",
(-5,1)*{}="4",
(-3.6,-3.6)*{}="5",
(3.6,-3.6)*{}="6",
(0,-4.5)*{}="7",
\ar @{~} "a";"1" <0pt>
\ar @{~} "a";"2" <0pt>
\ar @{~} "a";"5" <0pt>
\ar @{~} "a";"6" <0pt>
\endxy}
\Ea\right)
=\left( \sum
\Ba{c}\resizebox{9mm}{!}{\xy
%
%
 (0,-3.3)*{\bullet}="a", 
 (0,3.3)*{\bullet}="b",
(-5,-8)*{}="3",
(5,-8)*{}="4",
(5,8)*{}="5",
(-5,8)*{}="6",
\ar @{.>} "a";"b" <0pt>
\ar @{~} "a";"3" <0pt>
\ar @{~} "a";"4" <0pt>
\ar @{~} "b";"5" <0pt>
\ar @{~} "b";"6" <0pt>
\endxy}
\Ea\right).
$$
where any wavy hair can stand for a solid leg or a dotted half-edge
attached to $v$. For example, if all wavy hairs are solid legs, 
one obtains a Jacobi relation
$$
0=\Ba{c}\resizebox{10mm}{!}{\xy
 (0,-3.8)*{\bullet}="a", 
 (0,3.8)*{\bullet}="b",
(-5,-8)*{_1}="3",
(5,-8)*{_2}="4",
(5,8)*{_3}="5",
(-5,8)*{_0}="6",
\ar @{.>} "a";"b" <0pt>
\ar @{-} "a";"3" <0pt>
\ar @{-} "a";"4" <0pt>
\ar @{-} "b";"5" <0pt>
\ar @{-} "b";"6" <0pt>
\endxy}
\Ea
+
\Ba{c}\resizebox{10mm}{!}{\xy
 (0,-3.8)*{\bullet}="a", 
 (0,3.8)*{\bullet}="b",
(-5,-8)*{_2}="3",
(5,-8)*{_3}="4",
(5,8)*{_1}="5",
(-5,8)*{_0}="6",
\ar @{.>} "a";"b" <0pt>
\ar @{-} "a";"3" <0pt>
\ar @{-} "a";"4" <0pt>
\ar @{-} "b";"5" <0pt>
\ar @{-} "b";"6" <0pt>
\endxy}\Ea
+
\Ba{c}\resizebox{10mm}{!}{\xy
 (0,-3.8)*{\bullet}="a", 
 (0,3.8)*{\bullet}="b",
(-5,-8)*{_3}="3",
(5,-8)*{_1}="4",
(5,8)*{_2}="5",
(-5,8)*{_0}="6",
\ar @{.>} "a";"b" <0pt>
\ar @{-} "a";"3" <0pt>
\ar @{-} "a";"4" <0pt>
\ar @{-} "b";"5" <0pt>
\ar @{-} "b";"6" <0pt>
\endxy}
\Ea
$$
The case $d=3$ is special importance for this paper.

\sip



\subsection{\bf  Example: a cyclic operad of graphs $\Grac_d$}\label{3: subsec on Gra_cyc}
For any $d\in \Z$ there is a cyclic operad 
$$
\Grac_d:=\{\Grac_d((n+1))\}_{n\geq 1},
$$ 
where the
$\bS_{(n+1})$-module $\Grac_d((n+1))$
  is spanned by graphs $\Ga$ with labelled $n+1$ vertices and edges of degree $1-d$, e.g.
$$
  \Ba{c}\resizebox{11mm}{!}{ \xy
(-5,0)*+{_0}*\frm{o}="0";
 (5,0)*+{_1}*\frm{o}="1";
 \ar @{-} "1";"0" <0pt>
\endxy} \Ea \in \Grac_d((1+1)), \ \ \Ba{c}\resizebox{8mm}{!}{ \xy
(-3,0)*+{_0}*\frm{o}="0";
 (3,0)*+{_1}*\frm{o}="1";
 (0,6)*+{_2}*\frm{o}="2";
\endxy} \Ea,  \Ba{c}\resizebox{11mm}{!}{ \xy
(-5,0)*+{_0}*\frm{o}="0";
 (5,0)*+{_1}*\frm{o}="1";
 (0,8)*+{_2}*\frm{o}="2";
 \ar @{-} "1";"0" <0pt>
\endxy} \Ea, \Ba{c}\resizebox{11mm}{!}{ \xy
(-5,0)*+{_0}*\frm{o}="0";
 (5,0)*+{_1}*\frm{o}="1";
 (0,8)*+{_2}*\frm{o}="2";
 \ar @{-} "0";"2" <0pt>
 \ar @{-} "1";"0" <0pt>
 \ar @{-} "1";"2" <0pt>
\endxy} \Ea\in \Grac_d((2+1)).
$$  
 For  $d$ odd it is tacitly assumed that each edge of $\Ga$ comes equipped with a direction which can be flipped (together with the multiplication
 of  $\Ga$ by $-1$). For $d$ even the edges are undirected, but instead they are assumed to be ordered (up to a permutation and multiplication by the sign of that permutation). We do not show directions on edges for $d$ odd in our pictures assuming that some choice has been made.
 Each graph $\Ga\in \Grac_d((n+1))$ is assigned the cohomological degree 
 $$
 |\Ga|:=(1-d)\# E(\Ga),
 $$
 where $E(\Ga)$ is the set of edges of $\Ga$. 
  The permutation group $\bS_{(n+1)}$ acts on the generators of $\Grac_d((n+1))$ be relabeling the vertices.
  
  \sip

For a vertex $v\in V(\Ga)$ labelled by an integer $i\in ((n+1))$ we denote by $H(i)$  the set of all  half-edges  attached to $v$.
 
\sip

  The operadic
 composition of such graphs, for any $i\in ((n+1))$ and $j\in ((m+1))$,
$$
\Ba{rccc}
_i\circ_j: & \Grac_d((n+1))\ot \Grac_d((m+1)) & \lon & \Grac_d((m+n))\\
           &  \Ga' \ot \Ga''&\lon &  \Ga' \mbox{$_i\circ_j$} \Ga''
\Ea
$$
is defined on the generators as a linear combination of graphs
\Beq\label{3: ij-composition in Grac}
\Ga'  \mbox{$_i\circ_j$} \Ga'':= \sum_{f': H() \rar V(\Ga'')\atop
f'': H(j)\rar  V(\Ga')} \Ga' \mbox{$_{f'}\circ_{f''}$} \Ga''
\Eeq
where the summand $\Ga'_{f'}\circ_{f''} \Ga''$  corresponding to each  pair $(f',f'')$ of maps is given by the graph  obtained from $\Ga'$ and $\Ga''$ by erasing vertices labelled by
$i$ and, respectively, $j$ and reattaching the resulting dangling half-edges to the vertices of $\Ga''$
and, respectively, of $\Ga'$ in accordance with the rule specified by  $f'$ and, respectively, $f''$.
Here are two examples of the cyclic operadic composition in $\Grac_d$ for $d$ even
(when graphs can not have multiple edges)
$$
\Ba{c}\resizebox{8mm}{!}{ \xy
(-3,0)*+{_0}*\frm{o}="0";
 (3,0)*+{_1}*\frm{o}="1";
 (0,6)*+{_2}*\frm{o}="2";
\endxy} \Ea\
 \mbox{$_1\circ_0$}
 \Ba{c}\resizebox{11mm}{!}{ \xy
(-5,0)*+{_0}*\frm{o}="0";
 (5,0)*+{_1}*\frm{o}="1";
 \ar @{-} "1";"0" <0pt>
\endxy} \Ea
=\Ba{c}\resizebox{11mm}{!}{ \xy
(-5,0)*+{_0}*\frm{o}="0";
 (5,0)*+{_1}*\frm{o}="1";
 (0,8)*+{_2}*\frm{o}="2";
 \ar @{-} "1";"0" <0pt>
\endxy} \Ea
+
\Ba{c}\resizebox{11mm}{!}{ \xy
(-5,0)*+{_0}*\frm{o}="0";
 (5,0)*+{_1}*\frm{o}="1";
 (0,8)*+{_2}*\frm{o}="2";
 \ar @{-} "1";"2" <0pt>
\endxy} \Ea
\ \ , \ \ \
\Ba{c}\resizebox{11mm}{!}{ \xy
(-5,0)*+{_0}*\frm{o}="0";
 (5,0)*+{_1}*\frm{o}="1";
 (0,8)*+{_2}*\frm{o}="2";
 \ar @{-} "0";"2" <0pt>
 \ar @{-} "1";"0" <0pt>
\endxy} \Ea
 \mbox{$_1\circ_0$}
 \Ba{c}\resizebox{11mm}{!}{ \xy
(-5,0)*+{_0}*\frm{o}="0";
 (5,0)*+{_1}*\frm{o}="1";
 \ar @{-} "1";"0" <0pt>
\endxy} \Ea  =\Ba{c}\resizebox{11mm}{!}{ \xy
(-5,0)*+{_0}*\frm{o}="0";
 (5,0)*+{_1}*\frm{o}="1";
 (0,8)*+{_2}*\frm{o}="2";
 \ar @{-} "0";"2" <0pt>
 \ar @{-} "1";"0" <0pt>
 \ar @{-} "1";"2" <0pt>
\endxy} \Ea
$$
The induced orientations are obvious.

\sip

This cyclic operad of graphs is of interest to us because there exists a morphism of cyclic operads 
\Beq\label{2: BVc to Grac}
f: \BVc \lon \Grac_2
\Eeq
given on the generators by
\Beq\label{2: BVc to Grac formulae}
\Ba{ccc}
\Ba{c}\resizebox{9mm}{!}{\xy
 (0,1)*{_{_0}}="a",
(7,1)*{_{_1}}="b",
(3.5,0.6)*{^{_\Delta}},
\ar @{-} "a";"b" <0pt>
\endxy
}\Ea & \lon &  \Ba{c}\resizebox{11mm}{!}{ \xy
(-5,0)*+{_0}*\frm{o}="0";
 (5,0)*+{_1}*\frm{o}="1";
 \ar @{-} "1";"0" <0pt>
\endxy} \Ea \vspace{2mm} \\
\Ba{c}\begin{xy}
 <0mm,0.66mm>*{};<0mm,3mm>*{}**@{-},
 <0.39mm,-0.39mm>*{};<2.2mm,-2.2mm>*{}**@{-},
 <-0.35mm,-0.35mm>*{};<-2.2mm,-2.2mm>*{}**@{-},
 <0mm,0mm>*{\circ};<0mm,0mm>*{}**@{},
   <0mm,0.66mm>*{};<0mm,3.4mm>*{^{_0}}**@{},
   <0.39mm,-0.39mm>*{};<2.9mm,-4mm>*{^{_2}}**@{},
   <-0.35mm,-0.35mm>*{};<-2.8mm,-4mm>*{^{_1}}**@{},
\end{xy}\Ea
& \lon & \Ba{c}\resizebox{10mm}{!}{ \xy
(-4,0)*+{_0}*\frm{o}="0";
 (4,0)*+{_1}*\frm{o}="1";
 (0,5.8)*+{_2}*\frm{o}="2";
\endxy} \Ea
\Ea
\Eeq
This map is not a quasi-isomorphism, but in \S 7 below we extend $\Grac_2$ to a larger {\it dg}\, cyclic
operad $\BVGcyc$ of hypergraphs for which the map $i$ shown above becomes a quasi-isomorphism.

%

\subsection{\bf  Example: a cyclic operad of infinitesimal ribbon braids}\label{2: subsec on T_cyc}
 Let $\tilde{\ft}=\{\tilde{\ft}(n)\}_{n \geq 1}$ be the ({\it non-cyclic}) operad
of infinitesimal ribbon braids introduced in \cite{S}; this is an operad in the symmetric monoidal category of Lie
algebras  such that $\tilde{\ft}(n)$ is the completed free Lie algebra generated by symbols $t_{ij}=t_{ji}$, $i,j\in [n]$ modulo the Lie ideal generated by the
following relations
$$
\Ba{rcl}
[t_{ij}, t_{kl}] &=& 0\ \ \  \text{for}\ \{i,j\}\cap \{k,l\}=\emptyset  ,
 \\
\mbox{$[ t_{ij}, t_{ik} + t_{jk} ]$} &=& 0.   \\
\Ea
$$
In particular, the elements $t_{ii}$ (denoted by $2s_i$ in \cite{S}) are central in $\tilde{\ft}(n)$ for any $i\in [n]$.
The (non-cyclic) operadic compositions are given explicitly in \S 3.1 of \cite{S}.

The Lie subalgebra of $\tilde{\ft}(n)$  generated by $n$ elements $\{t_{in}\}_{i\in [n]}$ can be identified with
with ${\mathfrak{freeLie}}_{n-1}\oplus \K \langle t_{nn}\rangle$ where $\mathfrak{freeLie}_{n-1}$ is
the completed free Lie algebra on $n-1$ generators $\{t_{in}\}_{i\in [n-1]}$, and $\K \langle t_{nn}\rangle$ is the 1-dimensional center. 
It forms a Lie ideal and gives rise to a short exact sequence of Lie algebras
$$
0 \lon \widehat{\mathfrak{freeLie}}_{n-1}\oplus \K\langle t_{nn}\rangle \lon \tilde{\ft}(n) \lon \tilde{\ft}(n-1) \lon 0.
$$

\sip

 It is shown in \cite{CIW}
that this operad has a cyclic structure with the action of the permutation
$(01)\in \bS_{(n+1)}$ on the generators given, for any $i,j\geq 2$, by the formulae
$$
t_{ij} \rar t_{ij}, \ \
t_{1i} \rar - \sum_{k=1}^n t_{ik},  \ \
t_{11} \rar \sum_{k,l=1}^n t_{kl}.
$$
The action of $\bS_{(n+1)}$ on $\tilde{\ft}(n)$ can be made simpler if one introduces in $\tilde{\ft}(n)$ a
new set of generators,
\Beqrn
T_{0i} &=& -t_{ii} - 2 \sum_{k=1}t_{ik} - \sum_{k,l=1}^n t_{kl}\ \ \ \text{for any}\ i\in [n],\\
T_{ij} &=& 2t_{ij} - t_{ii} - t_{jj},\ \ \ \text{for any}\ i,j\in [n], i\neq j.
\Eeqrn
Then the action becomes just a relabeling of indices, i.e. for any $A,B\in ((n+1))$ with $A\neq B$ and any
$\sigma\in \bS_{(n+1)}$
one has
$$
T_{AB} \lon T_{\sigma(A)\sigma(B)}.
$$
Hence $\tilde{\ft}(n)$ can be equivalently defined as the completed Lie algebra $\ftt((n+1))$ generated by
the symbols
 $\{T_{AB}\}_{A,B\in ((n+1)), A\neq B}$ modulo the following relations
\Beqrn
 &&\hspace{-20mm}[T_{AB}, T_{CD}]= 0\ \ \ \text{for}\ \{A,B\}\cap \{C,D\}=\emptyset,   \\
 && \hspace{-20mm}\mbox{$\displaystyle [\sum_{B=0}^n T_{AB}, T_{CD}]$}=0 \ \ \text{for any $A,C,D\in ((n+1))$}.   
\Eeqrn
Here and elsewhere we always assume that $T_{AA}\equiv 0$ for any $A\in ((n+1))$.
These relations imply, for any $i,j,k\in [n]$,
$$
\Ba{rcl}
[T_{ij}, T_{kl}] &=& 0\ \ \text{if}\ \ \ \#\{i,j,k,l\}=4   \\
\mbox{$[ T_{ij}, T_{ik} + T_{jk} ]$} &=& 0 \ \ \text{if}\ \ \ \#\{i,j,k\}=3   \\
\mbox{$[ T_{ij}, T_{0i} + T_{0j} ]$} &=& 0    \\
\mbox{$[ T_{0k}, T_{0i} + T_{ik} ]$} &=& 0   \\
\Ea
$$
which in turn imply  that the subspace $\ft_n((n+1))$ of the Lie algebra $\ft^c((n+1))$ generated by the set
$\{T_{An}\}_{A\in ((n+1))}$
is a  Lie ideal so that there a short exact sequence  of Lie algebras
$$
0\lon \ft_n((n+1)) \lon \ftt((n+1)) \lon \ftt((n)) \lon 0.
$$

Thus the Lie ideal $\ft_n((n+1))$ can be identified the quotient of the completed free Lie algebra generated by 
symbols $T_{0n},T_{1n}, \ldots, T_{n-1\,n}$ modulo one relation saying that the sum
$$
T_{0n}+ T_{1n} + \ldots + T_{n-1\, n}
$$
is a central element in $\ft_n((n+1))$. If one considers a filtration of $\ft_n((n+1))$ by the number of letters
$T_{ij}$, $i,j\in [n]$ in the Lie words, then the associated graded Lie algebra
can be identified with ${\mathfrak{freeLie}}_{n-1}\oplus \K\langle T_{0n}\rangle$, where ${\mathfrak{freeLie}}_{n-1}$ is the completed  free Lie algebra in $T_{in}$, $i\in \{1,2,\ldots,n-1\}$.

\sip

The collection
$$
\ftt:=\{\ftt((n+1))\}_{n\geq 1}
$$
is a cyclic operad (see \S 3.1 in \cite{CIW}) in the category of Lie algebras. Hence 
$$
\odot^\bu(\ftt[1]):=\{\odot^\bu(\ftt((n+1))[1])\}_{n\geq 1}
$$
is the associated dg cyclic operad in the category of codefifferential  cocommutative coalgebras. There is a quasi-isomorphism of cyclic operads \cite{S}
\Beq\label{2: BVc to CE(T^c)}
\BVc \lon \odot^\bu(\ftt[1])
\Eeq
given on the generators as follows (cf.\ (\ref{2: BVc to Grac formulae}))
$$
\Ba{ccc}
\Ba{c}\resizebox{9mm}{!}{\xy
 (0,1)*{_{_0}}="a",
(7,1)*{_{_1}}="b",
(3.5,0.6)*{^{_\Delta}},
\ar @{-} "a";"b" <0pt>
\endxy
}\Ea & \lon & T_{01} \in \ftt((2))[1],   \\
\Ba{c}\begin{xy}
 <0mm,0.66mm>*{};<0mm,3mm>*{}**@{-},
 <0.39mm,-0.39mm>*{};<2.2mm,-2.2mm>*{}**@{-},
 <-0.35mm,-0.35mm>*{};<-2.2mm,-2.2mm>*{}**@{-},
 <0mm,0mm>*{\circ};<0mm,0mm>*{}**@{},
   <0mm,0.66mm>*{};<0mm,3.4mm>*{^{_0}}**@{},
   <0.39mm,-0.39mm>*{};<2.9mm,-4mm>*{^{_2}}**@{},
   <-0.35mm,-0.35mm>*{};<-2.8mm,-4mm>*{^{_1}}**@{},
\end{xy}\Ea
& \lon &  1\in   \odot^\bu(\ftt((3))[1]). 
\Ea
$$
This result was obtained in \cite{S} along the lines of Tamarkin's proof  of its version for the operad of Gerstenhaber algebras in \cite{T}.

\bip


{\large
\section{\bf  An auxiliary cyclic operad of hypergraphs and its twisting}
}

\sip

\subsection{Graphs and hypergraphs}\label{2: subsec on hypergraphs}  A {\em hypergraph}\, $\Ga$ can be understood as a triple of sets
($V(\Ga)$, $V_*(\Ga)$, $F(\Ga)$) together with an injective map $\iota: F(\Ga) \rar V(\Ga)\times V_*(\Ga)$. Composing that map $\iota$ with the canonical projections one gets the associated {\em structure maps}\,
$\iota_V: F(\Ga)\rar V(\Ga)$ and $\iota_*: F(\Ga)\rar V_*(\Ga)$.
The set $V(\Ga)$ is called the set of {\em vertices}, $V_*(\Ga)$ is called the set of its {\em hyperedges} (or $*$-{\em vertices})\footnote{The term {\em $*$-vertices}\, is useful in pictorial presentations of hypergraphs. In this paper we use both terms --- the $*$-vertices and hyperedges --- for one and the same structure depending on the context.}, and $F(\Ga)$ is called the set of {\em flags}\, of $\Ga$.

\sip

For any vertex $v\in V(\Ga)$ (resp., any hyperedge $v_*\in V_*(\Ga)$) the pre-image
$\iota^{-1}_V(v)$ (resp., $\iota^{-1}_*(v_*)$) is called the set of flags attached to $v$ (resp., to $v_*$), and their cardinalities are called {\em valencies}. An {\em automorphism} of a hypergraph $\Ups$ is a triple of bijections
$$
f_V: V(\Ga) \rar V(\Ga), \ \ \ f_{*}: V_*(\Ga) \rar V_*(\Ga), \ \ \ f_F: F(\Ga) \rar F(\Ga)
$$
which respect the structure maps,
$$
\iota_V \circ f_F= f_V \circ \iota_V, \ \iota_* \circ f_F= f_* \circ \iota_*.
$$

\sip

There is an obvious geometric representation of any hypergraph as a 1-dimensional CW-complex,
whose 0-cells are vertices and $*$-vertices, and whose 1-cells are flags which always
connect vertices of different types, e.g.

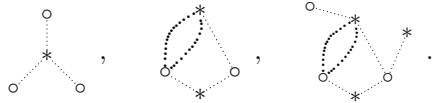
$$
\Ba{c}\resizebox{11.5mm}{!}{
\xy
(0,3)*{\ast}="1";
(-5,-2)*{\circ}="2";
(0,9)*{\circ}="3";
(5,-2)*{\circ}="4";
\ar @{.} "1";"2" <0pt>
\ar @{.} "1";"3" <0pt>
\ar @{.} "1";"4" <0pt>
\endxy} \Ea, \ \ \
\Ba{c}\resizebox{12mm}{!}{
\xy
(0,5)*{\ast}="1";
(0,-7)*{\ast}="2";
(-5,-4)*{\circ}="3";
(5,-4)*{\circ}="4";
"1";"3" **\crv{~*=<2pt>{.} (-3,5) & (5,4)};
"1";"3" **\crv{~*=<2pt>{.} (-5,2) & (-5,2)};
\ar @{.} "1";"4" <0pt>
\ar @{.} "2";"3" <0pt>
\ar @{.} "2";"4" <0pt>
\endxy}
\Ea,
 \ \
\Ba{c}\resizebox{15.4mm}{!}{
\xy
(0,5)*{\ast}="1";
(0,-7)*{\ast}="2";
(8,3)*{\ast}="5";
(-7,7)*{\circ}="6";
(-5,-4)*{\circ}="3";
(5,-4)*{\circ}="4";
"1";"3" **\crv{~*=<2pt>{.} (-3,5) & (5,4)};
"1";"3" **\crv{~*=<2pt>{.} (-5,2) & (-5,2)};
\ar @{.} "1";"4" <0pt>
\ar @{.} "2";"3" <0pt>
\ar @{.} "2";"4" <0pt>
\ar @{.} "5";"4" <0pt>
\ar @{.} "6";"1" <0pt>
\endxy}
\Ea.
$$

Hypergraphs whose $*$-vertices are all bivalent are called {\em graphs}, and their flags are called {\em half-edges}; in this case $*$-vertices are erased in pictures and the corresponding hyperdges  are called {\em edges}, e.g.
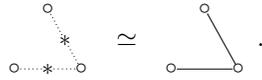
$$
\Ba{c}\resizebox{11.5mm}{!}{
\xy
(2.59,3.25)*{\ast}="1r";
(0,-1)*{\ast}="1d";
(0,8)*{\circ}="3";
(-5,-1)*{\circ}="2";
(5,-1)*{\circ}="4";
\ar @{.} "1r";"4" <0pt>
\ar @{.} "1r";"3" <0pt>
\ar @{.} "1d";"4" <0pt>
\ar @{.} "1d";"2" <0pt>
\endxy} \Ea \ \simeq
\Ba{c}\resizebox{11.5mm}{!}{
\xy
(0,8)*{\circ}="3";
(-5,-1)*{\circ}="2";
(5,-1)*{\circ}="4";
\ar @{-} "4";"2" <0pt>
\ar @{-} "3";"4" <0pt>
\endxy
}\Ea.
$$
In this paper we adopt the following ``pictorial" convention: flags of any hyperedge with valency of its $*$-vertex greater than 2 are shown as {\em dotted edges}, while bivalent hyperedges (i.e. edges)
are always shown as {\em solid edges}, precisely as in the picture just above and the following one
$$
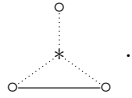
\Ba{c}\resizebox{15mm}{!}{ \xy
(-7,0)*{\circ}="0";
(7,0)*{\circ}="1";
(0,12)*{\circ}="2";
  (0,5)*{*}="s";
 \ar @{.} "s";"2" <0pt>
 \ar @{.} "s";"0" <0pt>
 \ar @{.} "s";"1" <0pt>
 \ar @{-} "0";"1" <0pt>
\endxy} \Ea.
$$

 \sip

Thus hypergraphs can be understood as graphs with two types of vertices whose ``edges" always connect vertices of different types. However this interpretation can be very misleading as
vertices and $*$-vertices are treated in our applications in radically different ways --- the vertices stand for points in configuration spaces while edges and hyperedges stand for propagators, that is, semialgebraic differential forms associated with the vertices to which they are attached.

\subsection{\bf  A cyclic operad of hypergraphs $f\HGrac_d$}\label{4: subsec on HGra}
For any integer $d\in \Z$ there is a cyclic operad
 $\HGrac_d:=\{\HGrac_d((n+1))\}$ where
the $\bS_{(n+1)}$-module $\HGrac_d((n+1))$ is generated by hypergraphs $\Ga$ with labelled $n+1$
vertices and {\it trivalent}\, hyperedges of cohomological degree $1-d$; more precisely, every $*$-vertex
is assigned the cohomological degree $1+2d$ while its three flags are assigned the degree $-d$ each, e.g.
 $$
  \Ba{c}\resizebox{11mm}{!}{ \xy
(-5,0)*+{_0}*\frm{o}="0";
 (5,0)*+{_1}*\frm{o}="1";
 (0,8)*+{_2}*\frm{o}="2";
\endxy} \Ea\in \HGrac_d((2+1)), \
 \Ba{c}\resizebox{15mm}{!}{ \xy
(-7,0)*+{_0}*\frm{o}="0";
 (7,0)*+{_1}*\frm{o}="1";
 (0,10)*+{_2}*\frm{o}="2";
  (0,3.5)*{*}="s";
 \ar @{.} "s";"2" <0pt>
 \ar @{.} "s";"0" <0pt>
 \ar @{.} "s";"1" <0pt>
\endxy} \Ea \in \HGrac_d((2+1)), \ \ \
 \Ba{c}\resizebox{20mm}{!}{ \xy
(0,5)*+{_0}*\frm{o}="0";
 (0,-5)*+{_1}*\frm{o}="1";
 (-8,0)*+{_2}*\frm{o}="2";
  (-3,0)*{*}="s1";
  (8,0)*+{_3}*\frm{o}="3";
  (3,0)*{*}="s2";
 \ar @{.} "s1";"2" <0pt>
 \ar @{.} "s1";"0" <0pt>
 \ar @{.} "s1";"1" <0pt>
  \ar @{.} "s2";"3" <0pt>
 \ar @{.} "s2";"0" <0pt>
 \ar @{.} "s2";"1" <0pt>
\endxy} \Ea\in \HGrac_d((3+1)).
 $$
 Note that for $d$ odd (which is the most important for us case) these hypergraphs can not have  parallel flags; moreover, in this
  case each hypergraph comes equipped with an orientation which is, by definition, an ordering of
  all its $*$-vertices and flags (up to their permutation and multiplying by the sign of that permuation).
 The permutation group $\bS_{(n+1)}$ acts on hypergraphs by relabeling the vertices. The operadic
 composition of such hypergraphs
$$
\Ba{rccc}
_i\circ_j: & \HGrac_d((n+1))\ot \HGrac_d((m+1)) & \lon & \HGrac_d((m+n))\\
           &  \Ga' \ot \Ga''&\lon &  \Ga' \mbox{$_i\circ_j$} \Ga''
\Ea
$$
is defined as a linear combination of hypergraphs (cf.\ \S{\ref{3: subsec on Gra_cyc}})
$$
\Ga'  \mbox{$_i\circ_j$} \Ga'' \sum_{f': F(i) \rar V(\Ga'')\atop
f'': F(j)\rar  V(\Ga')} \Ga' \mbox{$_{f'}\circ_{f''}$} \Ga''
$$
parameterized by all possible maps 
$$
f': F(i)\rar  V(\Ga'') \ \text{and} \  f'': F(j)\rar  V(\Ga')
$$
where
\Bi
\item[(1)] $F(i)$ (resp., $F(j)$) is the set of flags attached to the vertex $v_i\in \Ga'$
labelled by $i$ (resp., the vertex $v_j\in V(\Ga'')$ labelled by $j$), and 
\item[(2)]
for any given pair $(f',f'')$
the graph $\Ga'_{f'}\circ_{f''} \Ga''$ is obtained from $\Ga'$ and $\Ga''$ by erasing vertices
$i$ and, respectively, $j$ and reattaching the resulting dangling flags to the vertices of $\Ga''$
and, respectively, of $\Ga'$ in accordance with the maps  $f'$ and,
respectively, $f''$.
\Ei
The induced orientation in each summand is equal to $or_{\Ga'}\wedge or_{\Ga''}$.

\sip

 Here is an example of the cyclic operadic composition
$$
  \Ba{c}\resizebox{11mm}{!}{ \xy
(-5,0)*+{_0}*\frm{o}="0";
 (5,0)*+{_1}*\frm{o}="1";
 (0,8)*+{_2}*\frm{o}="2";
\endxy} \Ea \mbox{$_2\circ_0$}
 \Ba{c}\resizebox{15mm}{!}{ \xy
(-7,0)*+{_0}*\frm{o}="0";
 (7,0)*+{_1}*\frm{o}="1";
 (0,10)*+{_2}*\frm{o}="2";
  (0,3.5)*{*}="s";
 \ar @{.} "s";"2" <0pt>
 \ar @{.} "s";"0" <0pt>
 \ar @{.} "s";"1" <0pt>
\endxy} \Ea
=
 \Ba{c}\resizebox{13mm}{!}{ \xy
(-5,0)*+{_0}*\frm{o}="0";
 (5,0)*+{_1}*\frm{o}="1";
 (-5,10)*+{_2}*\frm{o}="2";
  (5,10)*+{_3}*\frm{o}="3";
  (0,5)*{*}="s";
 \ar @{.} "s";"2" <0pt>
 \ar @{.} "s";"0" <0pt>
  \ar @{.} "s";"3" <0pt>
\endxy} \Ea
+
 \Ba{c}\resizebox{13mm}{!}{ \xy
(-5,0)*+{_0}*\frm{o}="0";
 (5,0)*+{_1}*\frm{o}="1";
 (-5,10)*+{_2}*\frm{o}="2";
  (5,10)*+{_3}*\frm{o}="3";
  (0,5)*{*}="s";
 \ar @{.} "s";"2" <0pt>
 \ar @{.} "s";"1" <0pt>
  \ar @{.} "s";"3" <0pt>
\endxy} \Ea
$$

\subsubsection{\bf Lemma}\label{3: f from Lie_d to Gra_d} {\em For any odd integer $d$ there exists a morphism
of cyclic operads
$$
f: \Lie_d \lon \HGrac_d
$$
given on the generator by}
$$
\Ba{c}\begin{xy}
 <0mm,0.66mm>*{};<0mm,3mm>*{}**@{-},
 <0.39mm,-0.39mm>*{};<2.2mm,-2.2mm>*{}**@{-},
 <-0.35mm,-0.35mm>*{};<-2.2mm,-2.2mm>*{}**@{-},
 <0mm,0mm>*{\bu};<0mm,0mm>*{}**@{},
   <0mm,0.66mm>*{};<0mm,3.4mm>*{^{_0}}**@{},
   <0.39mm,-0.39mm>*{};<2.9mm,-4mm>*{^{_2}}**@{},
   <-0.35mm,-0.35mm>*{};<-2.8mm,-4mm>*{^{_1}}**@{},
\end{xy}\Ea
\lon
 \Ba{c}\resizebox{13mm}{!}{ \xy
(-6,0)*+{_1}*\frm{o}="1";
 (6,0)*+{_2}*\frm{o}="2";
 (0,10)*+{_0}*\frm{o}="0";
  (0,4)*{*}="s";
 \ar @{-} "s";"2" <0pt>
 \ar @{-} "s";"0" <0pt>
 \ar @{-} "s";"1" <0pt>
\endxy} \Ea
$$
\begin{proof}
 We have
\Beqrn
f\left(
\Ba{c}\resizebox{9mm}{!}{ \begin{xy}
 <0mm,0mm>*{\bu};<0mm,0mm>*{}**@{},
 <0mm,0.69mm>*{};<0mm,4.3mm>*{^0}**@{-},
 <-0.39mm,-0.39mm>*{};<-2.4mm,-2.4mm>*{}**@{-},
 <0.35mm,-0.35mm>*{};<1.9mm,-1.9mm>*{}**@{-},
 <2.4mm,-2.4mm>*{\bu};<2.4mm,-2.4mm>*{}**@{},
 <2.0mm,-2.8mm>*{};<0mm,-4.9mm>*{}**@{-},
 <2.8mm,-2.9mm>*{};<4.7mm,-4.9mm>*{}**@{-},
    <-0.39mm,-0.39mm>*{};<-3.3mm,-4.0mm>*{_1}**@{},
    <2.0mm,-2.8mm>*{};<-0.5mm,-6.7mm>*{_2}**@{},
    <2.8mm,-2.9mm>*{};<5.2mm,-6.7mm>*{_3}**@{},
 \end{xy}}\Ea
 \right)
&=&  \Ba{c}\resizebox{12mm}{!}{ \xy
(-6,0)*+{_1}*\frm{o}="1";
 (6,0)*+{_2}*\frm{o}="2";
 (0,10)*+{_0}*\frm{o}="0";
  (0,4)*{*}="s";
 \ar @{-} "s";"2" <0pt>
 \ar @{-} "s";"0" <0pt>
 \ar @{-} "s";"1" <0pt>
\endxy} \Ea \mbox{$_2\circ_0$}
 \Ba{c}\resizebox{12mm}{!}{ \xy
(-6,0)*+{_0}*\frm{o}="1";
 (6,0)*+{_2}*\frm{o}="2";
 (0,10)*+{_1}*\frm{o}="0";
  (0,4)*{*}="s";
 \ar @{-} "s";"2" <0pt>
 \ar @{-} "s";"0" <0pt>
 \ar @{-} "s";"1" <0pt>
\endxy} \Ea
 \\
&=& \Ba{c}\resizebox{13mm}{!}{ \xy
(-7,0)*+{_0}*\frm{o}="0";
 (7,0)*+{_2}*\frm{o}="2";
 (-7,10)*+{_1}*\frm{o}="1";
  (7,10)*+{_3}*\frm{o}="3";
  (-2,5)*{*}="s1";
  (2,5)*{*}="s2";
 \ar @{-} "s1";"1" <0pt>
 \ar @{-} "s1";"0" <0pt>
  \ar @{-} "s1";"3" <0pt>
\ar @{-} "s2";"2" <0pt>
\ar @{-} "s2";"3" <0pt>
\ar @{-} "s2";"1" <0pt>
\endxy} \Ea
+
\Ba{c}\resizebox{13mm}{!}{ \xy
(-7,0)*+{_0}*\frm{o}="0";
 (7,0)*+{_2}*\frm{o}="2";
 (-7,10)*+{_1}*\frm{o}="1";
  (7,10)*+{_3}*\frm{o}="3";
  (-2,5)*{*}="s1";
  (2,5)*{*}="s2";
 \ar @{-} "s1";"2" <0pt>
 \ar @{-} "s1";"0" <0pt>
  \ar @{-} "s1";"1" <0pt>
\ar @{-} "s2";"2" <0pt>
\ar @{-} "s2";"3" <0pt>
\ar @{-} "s2";"1" <0pt>
\endxy} \Ea
+
\Ba{c}\resizebox{13mm}{!}{ \xy
(-7,0)*+{_0}*\frm{o}="0";
 (7,0)*+{_2}*\frm{o}="2";
 (-7,10)*+{_1}*\frm{o}="1";
  (7,10)*+{_3}*\frm{o}="3";
  (-2,5)*{*}="s1";
  (2,5)*{*}="s2";
 \ar @{-} "s1";"1" <0pt>
 \ar @{-} "s1";"0" <0pt>
  \ar @{-} "s1";"3" <0pt>
\ar @{-} "s2";"2" <0pt>
\ar @{-} "s2";"3" <0pt>
\ar @{-} "s2";"0" <0pt>
\endxy} \Ea
+
\Ba{c}\resizebox{13mm}{!}{ \xy
(-7,0)*+{_0}*\frm{o}="0";
 (7,0)*+{_2}*\frm{o}="2";
 (-7,10)*+{_1}*\frm{o}="1";
  (7,10)*+{_3}*\frm{o}="3";
  (-2,5)*{*}="s1";
  (2,5)*{*}="s2";
 \ar @{-} "s1";"1" <0pt>
 \ar @{-} "s1";"0" <0pt>
  \ar @{-} "s1";"2" <0pt>
\ar @{-} "s2";"2" <0pt>
\ar @{-} "s2";"3" <0pt>
\ar @{-} "s2";"0" <0pt>
\endxy} \Ea
\Eeqrn
Next applying to the r.h.s.\ the element $\Id + (123) + (123)^2$ of the group algebra $\K[\bS_{(3+1)}])$
and taking into account orientations of the hypergraphs, one obtains zero.
\end{proof}

From now on we assume that the integer parameter $d$
  takes only odd values.

\subsection{A hypergraph complex} We define {\em a (full) hypergraph complex}\,  as the deformation complex of the above
  morphism of cyclic operads, 
$$
\fHGCc_d:=\Def(\Lie_d \stackrel{f}{\rar} \HGrac_d)=\prod_{n\geq 1}
\HGrac_d((n+1))\ot_{\bS_{n+1}} \sgn_{n+1}[2d-d(n+1)]
$$
This complex plays a purely auxiliary role in this paper: (a) we use it in twisting of $\HGrac$ discussed below (which in turn is used in the  construction of a dg cyclic operad $f\BVHgraphs$ in \S 6), and (b) we  need a  purely technical fact  about $\fHGCc_d$ (proven below) that it contains a subcomplex  $\HGCc_d$
generated by hypergraphs with all vertices {\it at least trivalent}. 

\sip

The dg Lie algebra $\fHGCc_d$  is generated by hypergraphs $\Ga$ whose vertices are unlabelled, e.g. 
\Beq\label{4: examples fHGC_d}
 \Ba{c}\resizebox{12mm}{!}{ \xy
(-6,0)*{\bu}="1";
 (6,0)*{\bu}="2";
 (0,10)*{\bu}="0";
  (0,4)*{*}="s";
 \ar @{.} "s";"2" <0pt>
 \ar @{.} "s";"0" <0pt>
 \ar @{.} "s";"1" <0pt>
\endxy} \Ea, \ \
\Ba{c}\resizebox{12mm}{!}{ \xy
(-6,0)*{\bu}="1";
 (6,0)*{\bu}="2";
 (0,10)*{\bu}="0";
  (-3,5)*{*}="s1";
  (3,5)*{*}="s2";
   (0,0)*{*}="s3";
 \ar @{.} "s1";"2" <0pt>
 \ar @{.} "s1";"0" <0pt>
 \ar @{.} "s1";"1" <0pt>
 \ar @{.} "s2";"2" <0pt>
 \ar @{.} "s2";"0" <0pt>
 \ar @{.} "s2";"1" <0pt>
  \ar @{.} "s3";"2" <0pt>
 \ar @{.} "s3";"0" <0pt>
 \ar @{.} "s3";"1" <0pt>
\endxy} \Ea \in \fHGCc_d.
\Eeq
The cohomological degree of  $\Ga\in \fHGCc_d$  is given by
$$
|\Ga|=d\# \sV_{\bu}(\Ga) + (1+2d)\# \sV_{*}(\Ga) - d\# \sF(\Ga) - 2d.
$$
Each such a hypergraph $\Ga$ is equipped
with an orientation given by the ordering of its vertices, flags and $*$-vertices up to a permutation
and the multiplication by the sign of that permutation.
For example, the hypergraph
$$
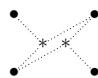
\Ba{c}\resizebox{13mm}{!}{ \xy
(-7,0)*{\bu}="0";
(7,0)*{\bu}="2";
(-7,10)*{\bu}="1";
(7,10)*{\bu}="3";
  (-2,5)*{*}="s1";
  (2,5)*{*}="s2";
 \ar @{.} "s1";"1" <0pt>
 \ar @{.} "s1";"0" <0pt>
  \ar @{.} "s1";"3" <0pt>
\ar @{.} "s2";"2" <0pt>
\ar @{.} "s2";"3" <0pt>
\ar @{.} "s2";"0" <0pt>
\endxy} \Ea
$$
vanishes identically in $\fHGCc_d$ as it admits an automorphism which reverses its orientation.

\sip

The Lie bracket in $\fHGCc_d$ is given by cyclic compositions,
$$
[\Ga', \Ga'']=\sum_{v'\in V_\bu(\Ga')\atop
v''\in V_\bu(\Ga'')} \Ga' \mbox{$_{v'}\circ_{v''}$} \Ga''
$$
The element
$$
\ga:= \Ba{c}\resizebox{12mm}{!}{ \xy
(-6,0)*{\bu}="1";
 (6,0)*{\bu}="2";
 (0,10)*{\bu}="0";
  (0,4)*{*}="s";
 \ar @{.} "s";"2" <0pt>
 \ar @{.} "s";"0" <0pt>
 \ar @{.} "s";"1" <0pt>
\endxy} \Ea
$$
is a Maurer-Cartan element corresponding to the morphism $f$ in the Lemma above. It makes $\fHGCc_d$ into a dg Lie algebra with the differential given by splitting the black vertices

$$
\delta\Ga=[\ga, \Ga]=\sum_{v\in \sV_\bu(\Ga)}\delta_v\Ga
$$
where $\delta_v \Ga$ is obtained from $\Ga$ by substituting into the vertex $v$ the following hypergraph,
\Beq\label{4: d on black in hGC}
\delta_v: \bu \ \lon  \  \frac{1}{2}
\Ba{c}\resizebox{9mm}{!}{ \xy
(0,-4)*{\bu}="1";
 (0,4)*{\bu}="2";
 (8,0)*{}="0";
  (3,0)*{*}="s";
 \ar @{.} "s";"2" <0pt>
 \ar @{.>} "s";"0" <0pt>
 \ar @{.} "s";"1" <0pt>
\endxy} \Ea
\Eeq
and taking two summations,
\Bi
\item[(a)] one summation goes over all possible reattachments of flags connected earlier to $v$ among
the two new black vertices, i.e.\ it goes over all possible splittings of $F(v)$ into the disjoint
union of two subsets,
\item[(b)] the second summation goes over attachments of the hanging flag (depicted as an arrow in the
picture above) to all other black vertices of $\Ga$, i.e.\ it a sum over the set $\sV(\Ga)\setminus v$.
\Ei
Here is an elementary pictorial check that $\delta$ is indeed a differential:
$$
\delta^2 \Ga= \sum_{v\in \sV_\bu(\Ga)}(\delta)^2_v\Ga
$$
where the sum of two operations
$$
(\delta)^2_v: \bu \lon
\Ba{c}\resizebox{9mm}{!}{ \xy
(0,-8)*{\bu}="1";
 (0,0)*{\bu}="2";
  (0,8)*{\bu}="3";
 (7,-4)*{}="0'";
  (3,-4)*{*}="s'";
  (7,4)*{}="0''";
  (3,4)*{*}="s''";
 \ar @{.} "s'";"2" <0pt>
 \ar @{.>} "s'";"0'" <0pt>
 \ar @{.} "s'";"1" <0pt>
  \ar @{.} "s''";"2" <0pt>
 \ar @{.>} "s''";"0''" <0pt>
 \ar @{.} "s''";"3" <0pt>
\endxy} \Ea
+
\Ba{c}\resizebox{18mm}{!}{ \xy
(0,-4)*{\bu}="1";
 (0,4)*{\bu}="2";
  (-8,0)*{\bu}="3";
 (8,0)*{}="0'";
  (3,0)*{*}="s'";
  (-3,0)*{*}="s''";
 \ar @{.} "s'";"2" <0pt>
 \ar @{.>} "s'";"0'" <0pt>
 \ar @{.} "s'";"1" <0pt>
  \ar @{.} "s''";"2" <0pt>
 \ar @{.} "s''";"1" <0pt>
 \ar @{.} "s''";"3" <0pt>
\endxy} \Ea
$$
vanishes identically as each operation creates a sum of  pairs of graphs with opposite orientations
which cancel each other
(arrows mean summations over attaching that edges to all other vertices of the graph, i.e.
excluding the ones shown in the picture).

\subsection{\bf At least trivalent hypergraph complex $\HGCc_d$} Let $\HGCc_{d}^{\geq 2}$
be a Lie subalgebra of $\fHGCc_d$ generated by {\it connected}\, hypergraphs with every black vertex having valency $\geq 2$.
Note that in contrast to the {\it full}\, Kontsevich graph complex $\fGC_d$ (the one in which no restructions of valencies of graph vertices are put, see \cite{Wi}), the subspace of $\fHGCc_d$ spanned by
connected hypergraphs with at least one {\it univalent}\,  black vertex is  {\it not}\,  a subcomplex
as the differential can increase the valency of black vertices by attaching a new flag.

\subsubsection{\bf Lemma}\label{4: Lemma on univalent in hGC} {\em $\HGCc_{d}^{\geq 2}$ is a dg Lie subalgebra of $\fcHGCc_d$.}

\begin{proof} It is obvious that this subspace is a Lie subalgebra.
It is less obvious that $\HGCc_{d}^{\geq 2}$ is a subcomplex, i.e. that the terms with univalent vertices in $\delta \Ga$ cancel out for any $\Ga\in \HGCc_{d}^{\geq 2}$. Indeed, the terms
with univalent vertices may appear in $\delta \Ga$  via attachments to every ordered pair
$(v,w)\in \sV(\Ga)\times \sV(\Ga)$ of different vertices $v\neq w$ a hyperedge with one univalent vertex
$$
\sum_{(v,w)\in \sV(\Ga)\times \sV(\Ga)\atop v\neq w}\Ba{c}\resizebox{12mm}{!}{ \xy
(-6,0)*{\bu}="1";
 (6,0)*{\bu}="2";
 (-6,-2)*{_v};
 (6,-2)*{_w};
 (0,10)*{\bu}="0";
  (0,4)*{*}="s";
 \ar @{.} "s";"2" <0pt>
 \ar @{.} "s";"0" <0pt>
 \ar @{.} "s";"1" <0pt>
\endxy} \Ea.
$$
Every such a summand occurs twice, once when $\delta$ splits $v$ and once when $\delta$ splits $w$;
both such summands cancel each other as they come with opposite orientations.
\end{proof}


\sip

Let $\HGCc_d$ be a linear subspace of $\HGCc_{d}^{\geq 2}$  spanned by graphs with all vertices at least trivalent.

\subsubsection{\bf Theorem}\label{4: trivalency theorem in hGC} {\em $\HGCc_{d}$ is a dg Lie subalgebra of $\HGCc_d^{\geq 2}$.}
\begin{proof}  The closure of $\HGCc_{d}$ under the Lie bracket is obvious.
Let us call a bivalent black vertex $v$ of a graph $\Ga\in \HGCc_d^{\geq 2}$  {\em bad}\, if
the two hyperedges
attached to $v$ share one more common black vertex (of any valency $\geq 2$); bad bivalent vertices fit ``squares"
of the form
$$
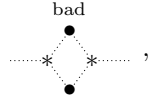
\Ba{c}\resizebox{17mm}{!}{ \xy
(0,-4)*{\bu}="1";
 (0,4)*{\bu}="2";
  (-8,0)*{}="3";
 (8,0)*{}="0'";
  (3,0)*{*}="s'";
  (0,6.4)*{^\text{bad}}="0''";
  (-3,0)*{*}="s''";
 \ar @{.} "s'";"2" <0pt>
 \ar @{.} "s'";"0'" <0pt>
 \ar @{.} "s'";"1" <0pt>
  \ar @{.} "s''";"2" <0pt>
 \ar @{.} "s''";"1" <0pt>
 \ar @{.} "s''";"3" <0pt>
\endxy} \Ea,
$$
where the lower black vertex can be of any valency $\geq 2$.
If $\Ga$ 
has no bivalent vertices, then the terms in $\delta\Ga$ which contain 
non-bad bivalent vertices, i.e. the ones of the form 
$$
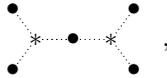
\Ba{c}\resizebox{19mm}{!}{ \xy
(-8,-4)*{\bu}="1";
 (-8,4)*{\bu}="2";
 (8,-4)*{\bu}="3";
 (8,4)*{\bu}="4";
 (0,0)*{\bu}="0";
  (-5,0)*{*}="s";
  (5,0)*{*}="s2";
 \ar @{.} "s";"2" <0pt>
 \ar @{.} "s";"0" <0pt>
 \ar @{.} "s";"1" <0pt>
  \ar @{.} "s2";"3" <0pt>
 \ar @{.} "s2";"s" <0pt>
 \ar @{.} "s2";"4" <0pt>
\endxy} \Ea\ ,
$$
 cancel out for the same reason as in the case of Kontsevich's  graph complex $\GC_d$ --- every such non-bad bivalent vertex appears twice in the sum with opposite signs.

\sip

Let us next show that for any $\Ga\in \HGCc_{d}$ the terms in $\delta\Ga$ with bivalent bad vertices cancel out as well.
Here the argument is a bit more complicated as we have to
study not pairs (as in the proof of Lemma~{\ref{4: Lemma on univalent in hGC}}), but triples
$(v_1,v_2,v_3)$ of black vertices of $\Ga$ such that they all are connected by one hyperedge in $\Ga$. The differential
$\delta$ applied to each vertex of such a triple gives us the following linear combination of hypergraphs,
each hypergraph having a newly created bivalent black vertex,
$$
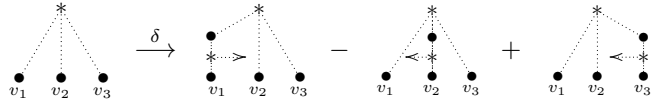
\Ba{c}\resizebox{15mm}{!}{ \xy
(-6,0)*{\bu}="1";
 (-0,0)*{\bu}="2";
 (6,0)*{\bu}="3";
 (-6,-2)*{_{v_1}};
  (0,-2)*{_{v_2}};
 (6,-2)*{_{v_3}};
  (0,10)*{*}="s";
 \ar @{.} "s";"2" <0pt>
 \ar @{.} "s";"3" <0pt>
 \ar @{.} "s";"1" <0pt>
\endxy} \Ea
\stackrel{\delta}{\lon}
\Ba{c}\resizebox{15mm}{!}{ \xy
(-7,0)*{\bu}="1";
(-7,6)*{\bu}="0";
(-7,3)*{*}="s0";
(-2,3)*{}="00";
 (-0,0)*{\bu}="2";
 (6,0)*{\bu}="3";
 (-6,-2)*{_{v_1}};
  (0,-2)*{_{v_2}};
 (6,-2)*{_{v_3}};
  (0,10)*{*}="s";
 \ar @{.} "s";"2" <0pt>
 \ar @{.} "s";"3" <0pt>
 \ar @{.} "s";"0" <0pt>
 \ar @{.} "s0";"1" <0pt>
 \ar @{.} "s0";"0" <0pt>
  \ar @{.>} "s0";"00" <0pt>
\endxy} \Ea
-
\Ba{c}\resizebox{15mm}{!}{ \xy
(-7,0)*{\bu}="1";
(0,6)*{\bu}="0";
(0,3)*{*}="s0";
(-4,3)*{}="00";
 (-0,0)*{\bu}="2";
 (6,0)*{\bu}="3";
 (-7,-2)*{_{v_1}};
  (0,-2)*{_{v_2}};
 (6,-2)*{_{v_3}};
  (0,10)*{*}="s";
 \ar @{.} "s";"2" <0pt>
 \ar @{.} "s";"3" <0pt>
 \ar @{.} "s";"0" <0pt>
 \ar @{.} "s";"1" <0pt>
 \ar @{.} "s0";"2" <0pt>
 \ar @{.} "s0";"0" <0pt>
  \ar @{.>} "s0";"00" <0pt>
\endxy} \Ea
+
\Ba{c}\resizebox{15mm}{!}{ \xy
(7,0)*{\bu}="1";
(7,6)*{\bu}="0";
(7,3)*{*}="s0";
(2,3)*{}="00";
 (-0,0)*{\bu}="2";
 (-6,0)*{\bu}="3";
 (7,-2)*{_{v_3}};
  (0,-2)*{_{v_2}};
 (-6,-2)*{_{v_1}};
  (0,10)*{*}="s";
 \ar @{.} "s";"2" <0pt>
 \ar @{.} "s";"3" <0pt>
 \ar @{.} "s";"0" <0pt>
 \ar @{.} "s0";"1" <0pt>
 \ar @{.} "s0";"0" <0pt>
  \ar @{.>} "s0";"00" <0pt>
\endxy} \Ea
$$
Bad bivalent vertices appear once the hanging edge at the vertex $v_i$ hits any one of the two
remaining vertices. This gives us six summands with bad bivalent vertices, and a straightforward
inspection of their induced orientations shows that they all cancel out.
\end{proof}

The above Lemma implies that, for example, the hypergraph
$$
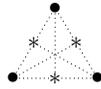
\Ba{c}\resizebox{14mm}{!}{ \xy
(-6,0)*{\bu}="1";
 (6,0)*{\bu}="2";
 (0,10)*{\bu}="0";
  (-3,5)*{*}="s1";
  (3,5)*{*}="s2";
   (0,0)*{*}="s3";
 \ar @{.} "s1";"2" <0pt>
 \ar @{.} "s1";"0" <0pt>
 \ar @{.} "s1";"1" <0pt>
 \ar @{.} "s2";"2" <0pt>
 \ar @{.} "s2";"0" <0pt>
 \ar @{.} "s2";"1" <0pt>
  \ar @{.} "s3";"2" <0pt>
 \ar @{.} "s3";"0" <0pt>
 \ar @{.} "s3";"1" <0pt>
\endxy} \Ea
$$
represents a non-trivial cohomology class in $H^{-d}(\HGCc_d)$ as the differential can not
create bivalent vertices and hence vanishes on this hypergraph.

\subsubsection{\bf Remark}\label{4: subsec map H from GC to HGC}
Consider an odd graph complex $\GC_d$, $d\in 2\Z +1$, and define a degree $-d$ linear map
$$
\Ba{rccc}
\caH: & \GC_d^{\geq 2} & \lon & \HGCc_d\\
   & \Ga            & \lon & \caH(\Ga)
\Ea
$$
by replacing each (directed up to the flip and the sign) edge of $\Ga$ by a hyperedge
with one ``hanging" flag,
$$
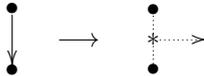
\Ba{c}\resizebox{3mm}{!}{ \xy
(0,-4)*{\bu}="1";
 (0,4)*{\bu}="2";
 \ar @{<-} "1";"2" <0pt>
\endxy} \Ea
\ \ \lon\ \ \Ba{c}\resizebox{9mm}{!}{ \xy
(0,-4)*{\bu}="1";
 (0,4)*{\bu}="2";
 (7,0)*{}="0";
  (0,0)*{*}="s";
 \ar @{.} "s";"2" <0pt>
 \ar @{.>} "s";"0" <0pt>
 \ar @{.} "s";"1" <0pt>
\endxy} \Ea
$$
and then summing up over all possible attachments of that newly created ``hanging" flag to
vertices of $\Ga$, and setting to zero graphs with at least one bivalent vertex.
For example
$$
\caH: \Ga=\Ba{c}\resizebox{9mm}{!}{\xy
 (0,0)*{\bullet}="a",
(8,0)*{\bu}="b",
(4,7)*{\bu}="c",
\ar @{->} "a";"b" <0pt>
\ar @{->} "a";"c" <0pt>
\ar @{->} "c";"b" <0pt>
\endxy}\Ea
\lon
\caH(\Ga)=\Ba{c}\resizebox{12mm}{!}{ \xy
(-6,0)*{\bu}="1";
 (6,0)*{\bu}="2";
 (0,10)*{\bu}="0";
  (-3,5)*{*}="s1";
  (3,5)*{*}="s2";
   (0,0)*{*}="s3";
 \ar @{.} "s1";"2" <0pt>
 \ar @{.} "s1";"0" <0pt>
 \ar @{.} "s1";"1" <0pt>
 \ar @{.} "s2";"2" <0pt>
 \ar @{.} "s2";"0" <0pt>
 \ar @{.} "s2";"1" <0pt>
  \ar @{.} "s3";"2" <0pt>
 \ar @{.} "s3";"0" <0pt>
 \ar @{.} "s3";"1" <0pt>
\endxy} \Ea
$$

\sip

The orientation of $\caH(\Ga)$ is uniquely determined by the orientation of $\Ga$ as a choice of a
direction on each edge $e$ of $\Ga$ gets encoded into the ordering of the two half-edges into which
the newly created $*$-vertex splits $e$. Note the map preserves the cohomological degrees of vertices and
(hyper)edges; the difference $|\caH(\Ga)| - |\Ga|=-d$ comes from the total degree shift $-2d$
in the case of hypergraphs from $\HGCc_d$ and the total shift $-d$ in the case of ordinary graphs from
$\GC_d^{\geq 2}$.

\sip

It is easy to see that the linear map $\caH: \GC_d \rar \HGCc_d$
is a morphism of complexes (which is not a quasi-isomorphism). This observation plays no role in the rest of this paper.

\subsection{Twisted extension of $\HGrac_d$} \label{3: subsec on twsisting of HGra}
The morphism $f$ in Lemma~{\ref{3: f from Lie_d to Gra_d}} can be used to {\it twist}\, the cyclic operad $\HGrac_d$  as follows (cf.\ Appendix I in \cite{Wi}).
For any odd integer $d$
let $f\HGraphs_d=\{f\HGraphs_d((n+1))\}_{n\geq 1}$ be a cyclic operad generated by hypergraphs with hyperedges as in $\HGrac_d$ but with
{\it  two}\, types of vertices,
labelled ones (depicted in white color), and unlabelled ones (depicted in black color), e.g.
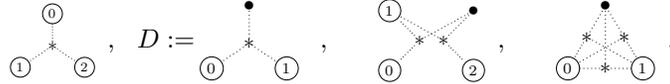
$$
\Ba{c}\resizebox{12mm}{!}{ \xy
(-6,0)*+{_1}*\frm{o}="1";
 (6,0)*+{_2}*\frm{o}="2";
 (0,10)*+{_0}*\frm{o}="0";
  (0,4)*{*}="s";
 \ar @{.} "s";"2" <0pt>
 \ar @{.} "s";"0" <0pt>
 \ar @{.} "s";"1" <0pt>
\endxy} \Ea, \ \ 
D:=\hspace{-3mm} \Ba{c}\resizebox{14mm}{!}{ \xy
(-6,0)*+{_0}*\frm{o}="1";
 (6,0)*+{_1}*\frm{o}="2";
 (0,10)*{\bu}="0";
  (0,4)*{*}="s";
 \ar @{.} "s";"2" <0pt>
 \ar @{.} "s";"0" <0pt>
 \ar @{.} "s";"1" <0pt>
\endxy} \Ea
\ , \ \ \
 \Ba{c}\resizebox{15mm}{!}{ \xy
(-7,0)*+{_0}*\frm{o}="0";
 (7,0)*+{_2}*\frm{o}="2";
 (-7,10)*+{_1}*\frm{o}="1";
  (7,10)*{\bu}="3";
  (-2,5)*{*}="s1";
  (2,5)*{*}="s2";
 \ar @{.} "s1";"1" <0pt>
 \ar @{.} "s1";"0" <0pt>
  \ar @{.} "s1";"3" <0pt>
\ar @{.} "s2";"2" <0pt>
\ar @{.} "s2";"3" <0pt>
\ar @{.} "s2";"1" <0pt>
\endxy} \Ea
, \ \ \
\Ba{c}\resizebox{14mm}{!}{ \xy
(-6,0)*+{_0}*\frm{o}="1";
 (6,0)*+{_1}*\frm{o}="2";
 (0,10)*{\bu}="0";
  (-3,5)*{*}="s1";
  (3,5)*{*}="s2";
   (0,0)*{*}="s3";
 \ar @{.} "s1";"2" <0pt>
 \ar @{.} "s1";"0" <0pt>
 \ar @{.} "s1";"1" <0pt>
 \ar @{.} "s2";"2" <0pt>
 \ar @{.} "s2";"0" <0pt>
 \ar @{.} "s2";"1" <0pt>
  \ar @{.} "s3";"2" <0pt>
 \ar @{.} "s3";"0" <0pt>
 \ar @{.} "s3";"1" <0pt>
\endxy} \Ea.
$$
The second hypergraph $D$ in the above list is skew-symmetric over $(01)$, and hence it defines a derivation of the cyclic operad
$f\HGraphs_d$. Let us consider a differential $\delta_\bu$ on $f\HGraphs_d$ which acts only on black
vertices of hypergraphs
exactly as in formula (\ref{4: d on black in hGC}).

\subsubsection{\bf Lemma}\label{4: Lemma on MC property of D} {\em The hypergraph
$
D=\hspace{-2mm} \Ba{c}\resizebox{13mm}{!}{ \xy
(-6,0)*+{_0}*\frm{o}="1";
 (6,0)*+{_1}*\frm{o}="2";
 (0,10)*{\bu}="0";
  (0,4)*{*}="s";
 \ar @{.} "s";"2" <0pt>
 \ar @{.} "s";"0" <0pt>
 \ar @{.} "s";"1" <0pt>
\endxy} \Ea
$
 satisfies the equation}\
 $
\delta_{\bu} D + D\ \mbox{${_1\circ_0}$}\ D=0
$.

\begin{proof}
One has
$$
\delta_\bu
\Ba{c}\resizebox{13mm}{!}{ \xy
(-6,0)*+{_0}*\frm{o}="1";
 (6,0)*+{_1}*\frm{o}="2";
 (0,10)*{\bu}="0";
  (0,4)*{*}="s";
 \ar @{.} "s";"2" <0pt>
 \ar @{.} "s";"0" <0pt>
 \ar @{.} "s";"1" <0pt>
\endxy} \Ea
=
\Ba{c}\resizebox{15mm}{!}{ \xy
(-7,0)*+{_0}*\frm{o}="0";
 (7,0)*+{_1}*\frm{o}="1";
 (0,6)*{\bu}="b1";
  (0,0)*{*}="s1";
   (0,18)*{\bu}="b1";
  (0,12)*{*}="s2";
 \ar @{.} "s1";"0" <0pt>
 \ar @{.} "s1";"1" <0pt>
 \ar @{.} "s1";"b1" <0pt>
 \ar @{.} "s2";"1" <0pt>
 \ar @{.} "s2";"b1" <0pt>
\endxy} \Ea
+
\Ba{c}\resizebox{15mm}{!}{ \xy
(-7,0)*+{_0}*\frm{o}="0";
 (7,0)*+{_1}*\frm{o}="1";
 (0,6)*{\bu}="b1";
  (0,0)*{*}="s1";
   (0,18)*{\bu}="b1";
  (0,12)*{*}="s2";
 \ar @{.} "s1";"0" <0pt>
 \ar @{.} "s1";"1" <0pt>
 \ar @{.} "s1";"b1" <0pt>
 \ar @{.} "s2";"0" <0pt>
 \ar @{.} "s2";"b1" <0pt>
\endxy} \Ea, \
\Ba{c}\resizebox{13mm}{!}{ \xy
(-7,0)*+{_0}*\frm{o}="0";
(7,0)*+{_1}*\frm{o}="1";
(0,12)*{\bu}="2";
  (0,5)*{*}="s";
 \ar @{.} "s";"2" <0pt>
 \ar @{.} "s";"0" <0pt>
 \ar @{.} "s";"1" <0pt>
\endxy} \Ea
\mbox{${_1\circ_0}$}
\Ba{c}\resizebox{13mm}{!}{ \xy
(-7,0)*+{_0}*\frm{o}="0";
(7,0)*+{_1}*\frm{o}="1";
(0,12)*{\bu}="2";
  (0,5)*{*}="s";
 \ar @{.} "s";"2" <0pt>
 \ar @{.} "s";"0" <0pt>
 \ar @{.} "s";"1" <0pt>
\endxy} \Ea
=
-
\Ba{c}\resizebox{15mm}{!}{ \xy
(-7,0)*+{_0}*\frm{o}="0";
 (7,0)*+{_1}*\frm{o}="1";
 (0,6)*{\bu}="b1";
  (0,0)*{*}="s1";
   (0,18)*{\bu}="b1";
  (0,12)*{*}="s2";
 \ar @{.} "s1";"0" <0pt>
 \ar @{.} "s1";"1" <0pt>
 \ar @{.} "s1";"b1" <0pt>
 \ar @{.} "s2";"1" <0pt>
 \ar @{.} "s2";"b1" <0pt>
\endxy} \Ea
-
\Ba{c}\resizebox{15mm}{!}{ \xy
(-7,0)*+{_0}*\frm{o}="0";
 (7,0)*+{_1}*\frm{o}="1";
 (0,6)*{\bu}="b1";
  (0,0)*{*}="s1";
   (0,18)*{\bu}="b1";
  (0,12)*{*}="s2";
 \ar @{.} "s1";"0" <0pt>
 \ar @{.} "s1";"1" <0pt>
 \ar @{.} "s1";"b1" <0pt>
 \ar @{.} "s2";"0" <0pt>
 \ar @{.} "s2";"b1" <0pt>
\endxy} \Ea
$$
Hence the claim.
\end{proof}

\sip

The above Lemma implies that the degree 1 derivation on the cyclic operad $f\HGraphs_d$ given by
\Beq\label{4: differential on HGrapgs}
\p \Ga= \delta_{\bu} \Ga + D(\Ga)
\Eeq
is a differential. It acts on black vertices by splitting them as in (\ref{4: d on black in hGC}) and
on white vertices by splitting them as follows,
$$
\p:  \Ba{c}\resizebox{4.5mm}{!}{ \xy
(0,0)*+{_i}*\frm{o}="0";
\endxy}\Ea
\lon
\Ba{c}\resizebox{9mm}{!}{ \xy
(0,-5)*+{_i}*\frm{o}="1";
 (0,4)*{\bu}="2";
 (8,0)*{}="0";
  (3,0)*{*}="s";
 \ar @{.} "s";"2" <0pt>
 \ar @{.>} "s";"0" <0pt>
 \ar @{.} "s";"1" <0pt>
\endxy} \Ea
$$
Let $\HGraphs_d$ be a dg cyclic sub-operad of $f\HGraphs_d$ spanned by graphs whose every black vertex is
at least trivalent and is connected to at least two white vertices by a connected path of flags. The following hypergraphs 
$$
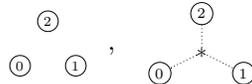
\Ba{c}\resizebox{11mm}{!}{ \xy
(-5,0)*+{_0}*\frm{o}="0";
 (5,0)*+{_1}*\frm{o}="1";
 (0,8)*+{_2}*\frm{o}="2";
\endxy} \Ea \ , \
 \Ba{c}\resizebox{15mm}{!}{ \xy
(-7,0)*+{_0}*\frm{o}="0";
 (7,0)*+{_1}*\frm{o}="1";
 (0,10)*+{_2}*\frm{o}="2";
  (0,3.5)*{*}="s";
 \ar @{.} "s";"2" <0pt>
 \ar @{.} "s";"0" <0pt>
 \ar @{.} "s";"1" <0pt>
\endxy} \Ea
$$
represent non-trivial cohomology classes
in $\HGraphs_d$.
The dg cyclic operad $\HGraphs_d$ is used below as one of the main building blocks
in the construction of a combinatorial hypergraph model of $\BVc$.



\bip

{\large
\section{\bf  Towards a hypergraph model for the cyclic operad $\BVc$}
}

\sip

\subsection{A geometric model for $\BVc$}
\label{5: subsec on fM and alpha-forms}
Let $\FFM_2=\{\FFM_2((n+1))\}_{n\geq 1}$ be the topological operad of compactified
moduli spaces of genus zero algebraic curves with labelled framed points; its homology operad
$H_\bu(\FFM_2)$ is precisely the cyclic
operad $\BVc$ of Batalin-Vilkovisky algebras \cite{KSV}.
The associated dual cyclic cooperad
$$
(\BVc)^*\simeq H^\bu(\FFM_2)=\{H^\bu(\FFM_2((n+1)))\}_{n\geq 1}
$$
is a Hopf cooperad, that is, a cooperad in the category of graded commutative algebras.
The cohomology groups $H^\bu(\FFM_2((n+1)))$ for each $n\geq 1$ are  generated as graded commutative algebras by  the
cohomology classes $\bar{\al}_{ij}$ of the closed 1-forms
$$
\al_{ij}=\al_{ji}:=\pi_{ij}^*(\text{Vol}_{S^1}), \ \forall i\neq j, \ i,j\in \{0,1,\ldots,n\},
$$
 obtained from  the standard homogeneous volume form on
$S^1\simeq \FFM_2((2))$ via the standard forgetful maps $\pi_{ij}: \FFM_2((n+1))^0\rar \FFM_2((2))$.
Let us consider the following 2-forms on $\FFM_2((n+1))$, $n\geq 1$,
\Beq\label{5: def of Omega_ijk}
\Omega_{ijk}:={\al}_{ij}{\al}_{jk} + {\al}_{ki}{\al}_{ij} + {\al}_{jk}{\al}_{ki}.
\Eeq
and let
\Beq\label{5: def of barOmega_ijk}
\bar{\Omega}_{ijk}:=\bar{\al}_{ij}\bar{\al}_{jk} + \bar{\al}_{ki}\bar{\al}_{ij} + \bar{\al}_{jk}\bar{\al}_{ki}.
\Eeq
stand for the associated cohomology class in $H^\bu(\FFM_2((n+1)))$. There is an isomorphism of graded commutative algebras \cite{GS}
$$
(\BVc)^*((n+1))\simeq H^\bu(\FFM_2((n+1)))=\frac{\R [\bar{\al}_{ij}]_{0\leq i\neq j\leq n}}{I},
$$
where $I$ is the ideal in the free graded commutative algebra $\R [\bar{\al}_{ij}]$ generated by the following {\em cyclic Arnold relations},
\Beq\label{5: Cyclic Arnold for Omega_ijk}
\bar{\Omega}_{jkl} - \bar{\Omega}_{ikl} + \bar{\Omega}_{ijl} -\bar{\Omega}_{ijk} \equiv 0  \ \ \ \forall i,j,k,l\in \{0,1,\ldots,n\}, \
\# \{i,j,k,l\}=4.
\Eeq


The cooperad structure in $(\BVc)^*$ is
given by the co-compositions,
$$
\Delta_{I,J}: (\BVc)^*(I) \lon  (\BVc)^*(I'\sqcup x') \ot (\BVc)^*(x''\sqcup I''),
$$
associated with any decomposition of the set $I=\{0,1,\ldots, n+1\}$ into the disjoint union of non-empty subsets
$I=I'\sqcup I''$,  which in turn are
 completely determined by their values on the co-generators given as follows (see \S 4.4 in \cite{GS})
\Beq\label{5: Delta on al}
\Delta_{I',I''}(\bar{\al}_{ij})=\left\{
\Ba{ll}
\bar{\al}_{ij}\ot 1 & \text{if}\ i,j\in I'\\
1\ot \bar{\al}_{ij} & \text{if}\ i,j\in I''\\
 \bar{\al}_{ix'}\ot 1 + 1\ot \bar{\al}_{x''j}  & \text{if}\ i\in I', j\in I''\\
\Ea
\right.
\Eeq
Note that these co-compositions induce the following
cocompositions of the 2-forms  (\ref{5: def of barOmega_ijk})
\Beq\label{5: Delta on Omega}
\Delta_{I',I''}(\bar{\Omega}_{ijk})=\left\{
\Ba{ll}
\bar{\Omega}_{ijk}\ot 1 & \text{if}\ i,j,k\in I'\\
1\ot \bar{\Omega}_{ijk} & \text{if}\ i,j,k\in I''\\
\bar{\Omega}_{ijx'}\ot 1 - (\bar{\al}_{ix'} - \bar{\al}_{jx'})\ot \bar{\al}_{x''k}  & \text{if}\ i,j\in I', k\in I''\\
1\ot \bar{\Omega}_{x''jk} + \bar{\al}_{ix'}\ot(\bar{\al}_{x''j} - \bar{\al}_{x''k})  & \text{if}\ i\in I', j,k\in I''\\
\Ea
\right.
\Eeq

\subsection{The dual cooperad $(\Grac_2)^*$}  In \S 2 we considered a cyclic operad $\Grac_2$ which 
comes equipped with an explicit morphism of cyclic operads (\ref{2: BVc to Grac})-(\ref{2: BVc to Grac formulae}). 
Its dual 
$$
(\Grac_2)^*=\left\{\left(\Grac_2((n+1))\right)^*\right\}_{n\geq 1}  
$$ 
is a Hopf cooperad such that $(\Grac_2)^*((n+1))$ is the free graded
commutative algebra $\K[[\theta_{ij}]]$ generated over a field $\K$ by degree 1 formal variables
\Beq\label{5: generators theta}
\theta_{ij}=\theta_{ji} \ \ \ \  \forall i\neq j, \ i,j\in \{0,1,\ldots,n\}.
\Eeq
Note that we do {\em not}\, assume here that $\theta_{ij}$ satisfy the cyclic Arnold relations
(\ref{5: Cyclic Arnold for Omega_ijk})-(\ref{5: def of Omega_ijk}). The co-compositions in $(\Grac_2)^*$ 
are given on the generators by (cf.\ (\ref{5: Delta on al}))
\Beq\label{Delta on theta small}
\Delta_{I',I''}(\theta_{ij})=\left\{
\Ba{ll}
\theta_{ij}\ot 1 & \text{if}\ i,j\in I'\\
1\ot \theta_{ij} & \text{if}\ i,j\in I''\\
 \theta_{ix'}\ot 1 + 1\ot \theta_{x''j}  & \text{if}\ i\in I', j\in I''\\
\Ea
\right.
\Eeq
Next we consider its hypergraph extension.

\subsection{An auxiliary cyclic Hopf cooperad $co\BVGrac$}  Consider a collection of graded commutative algebras
$$
co\BVGrac=\{co\BVGrac(n+1)\}_{n\geq 1},
$$
where $co\BVGrac(n+1)$ is the free graded commutative algebra $\K[[\theta_{ij}, \Theta_{ijk}]]$
generated by degree 1 formal variables (\ref{Delta on theta small})
 and by degree 2 formal variables
$\Theta_{ijk}$ which satisfy the following skewsymmetry conditions
\Beq\label{5: generators Theta}
\Theta_{ijk}=-\Theta_{jik}=-\Theta_{ikj} \ \ \ \  \forall i, j,k\in \{0,1,\ldots,n\}.
\Eeq


 Note
that we do {\em not}\, assume any relations between the generators $\theta_{ij}$ and $\Theta_{ijk}$
(as in (\ref{5: def of Omega_ijk})).

\subsubsection{\bf Lemma} {\em The $\bS$-module $co\BVGrac$ is a
cyclic Hopf cooperad with respect to the
co-compositions given on the generators by  (\ref{Delta on theta small}) and (cf.\ (\ref{5: Delta on Omega}))}
\Beq\label{5: Delta on Theta large}
\Delta_{I',I''}(\Theta_{ijk})=\left\{
\Ba{ll}
\Theta_{ijk}\ot 1 & \text{if}\ i,j,k\in I'\\
1\ot \Theta_{ijk} & \text{if}\ i,j,k\in I''\\
\Theta_{ijx'}\ot 1 - ({\theta}_{ix'} -
{\theta}_{jx'})\ot \theta_{x''k}  & \text{if}\ i,j\in I', k\in I''\\
1\ot \Theta_{x''jk} + {\theta}_{ix'}\ot({\theta}_{x''j} - \theta_{x''k})
 & \text{if}\ i\in I', j,k\in I''\\
\Ea
\right.
\Eeq

\begin{proof} We have to check the co-associativity identities for the generators $\Theta_{ijk}$, i.e.\ the axioms (\ref{3: Phi^R}) and (\ref{3: Phi^L}) when applied to the generators $\Theta_{ijk}$  for any decomposition of a finite set $I$ into the disjoint union of three subsets
$$
I=I'\sqcup I''\sqcup I''', \ \ \  \# I'\geq 2,\ \# I''\geq 2,\ \# I'''\geq 2.
$$
It is immediate to see that
$$
\Phi^L_{I',I'',I''}(\Theta_{ijk})=\Phi^R_{I',I'',I''}(\Theta_{ijk})=
\left\{
\Ba{ll}
\Theta_{ijk}\ot 1\ot 1  &\text{if}\ i,j,k\in I',\\
1\ot \Theta_{ijk} \ot 1  &\text{if}\ i,j,k\in I''\\\
1\ot 1\ot  \Theta_{ijk} &\text{if}\ i,j,k\in I'''.
\Ea
\right.
$$

If
$i,j\in I'$, $k\in I''$ we again have the equality as
$$
\Phi^L_{I',I'',I'''}(\Theta_{ijk})=
\left(\Delta_{I',I''\sqcup y'}\ot \Id\right) (\Theta_{ijk}\ot 1)=(\Theta_{ijx'}\ot 1 - ({\theta}_{ix'} -
{\theta}_{jx'})\ot \theta_{x''k})\ot 1
$$
and
\Beqrn
\Phi^R_{I',I'',I'''}(\Theta_{ijk})&=&
\left(1\ot \Delta_{I''\sqcup y', I'''}\right)(\Theta_{ijx'}\ot 1 - ({\theta}_{ix'} -
{\theta}_{jx'})\ot \theta_{x''k})\\
&=&
(\Theta_{ijx'}\ot 1 - ({\theta}_{ix'} -
{\theta}_{jx'})\ot \theta_{x''k})\ot 1.
\Eeqrn

If
$i,j\in I'$, $k\in I'''$, one has
\Beqrn
\Phi^L_{I',I'',I'''}(\Theta_{ijk})&=&
\left(\Delta_{I',I''\sqcup y'}\ot \Id\right)(\Theta_{ijy'}\ot 1 - ({\theta}_{iy'} -
{\theta}_{jy'})\ot \theta_{y''k})\\
&=&\Theta_{ijx'}\ot 1 \ot 1 -  ({\theta}_{ix'} -
{\theta}_{jx'})\ot \theta_{x''y'}\ot 1  \\
&& -({\theta}_{ix'}\ot 1 + 1\ot \theta_{x''y'} - {\theta}_{jx'}\ot 1 - 1\ot \theta_{x''y'})\ot \theta_{y''k}\\
&=&\Theta_{ijx'}\ot 1 \ot 1 -  ({\theta}_{ix'} - {\theta}_{jx'})\ot(\theta_{x''y'}\ot 1 +  1 \ot \theta_{y''k}),
\Eeqrn
while
\Beqrn
\Phi^R_{I',I'',I'''}(\Theta_{ijk})&=&
\left(1\ot \Delta_{I''\sqcup y', I'''}\right)(\Theta_{ijx'}\ot 1 - ({\theta}_{ix'} -
{\theta}_{jx'})\ot \theta_{x''k})\\
&=& \Theta_{ijx'}\ot 1 \ot 1 -  ({\theta}_{ix'} - {\theta}_{jx'})\ot(\theta_{x''y'}\ot 1 +  1 \ot \theta_{y''k})
\Eeqrn
i.e.\ the co-associativity condition is satisfied in this case as well.  Similarly one checks the conditions for $j,k\in I'''$ and $i$ in $I'$ or in $I''$.

\sip

It remains to check the condition for the most non-trivial case when $i\in I'$, $j\in I''$ and
$k\in I'''$. We have
\Beqrn
\Phi^L_{I',I'',I'''}(\Theta_{ijk})&=&
\left(\Delta_{I',I''\sqcup y'}\ot \Id\right)(\Theta_{ijy'}\ot 1 - ({\theta}_{iy'}-
{\theta}_{jy'})\ot \theta_{y''k})\\
&=& 1\ot \Theta_{x''jy'}\ot 1 + {\theta}_{ix'}\ot({\theta}_{x''j} - \theta_{x''y'})\ot 1
 - ({\theta}_{ix'}\ot 1 + 1\ot \theta_{x''y'} - 1\ot {\theta}_{jy'})\ot \theta_{y''k}
\Eeqrn
which gives precisely the same result as the following calculation,
\Beqrn
\Phi^R_{I',I'',I'''}(\Theta_{ijk})&=&
\left(1\ot \Delta_{I''\sqcup y', I'''}\right)(1\ot \Theta_{x''jk} + {\theta}_{ix'}\ot ({\theta}_{x''j} - \theta_{x''k})\\
&=& 1\ot \Theta_{x''jy'}\ot 1 - 1\ot ({\theta}_{x''y'} -
{\theta}_{jy'})\ot \theta_{y''k}
+ {\theta}_{ix'}\ot( {\theta}_{x''j} \ot 1 - \theta_{x''y'} \ot 1 -
1\ot \theta_{y''k})\\
&=&  \Phi^L_{I',I'',I'''}(\Theta_{ijk}).
\Eeqrn
The proof is completed.
\end{proof}


\subsection{A cyclic operad of hypergraphs}\label{4: Subsec on BVGrac}
The cyclic operad dual to the above cooperad,
$$
\BVGrac:=(co\BVGrac)^*,
$$
admits a nice combinatorial description in terms of hypergraphs with labelled white vertices with  hyperedges of two types,
$$
\text{the dual of $\Theta_{ijk}$} \rar
\Ba{c}\resizebox{12mm}{!}{ \xy
(-7,0)*+{_i}*\frm{o}="0";
(7,0)*+{_j}*\frm{o}="1";
(0,12)*+{_k}*\frm{o}="2";
  (0,5)*{*}="s";
 \ar @{.} "s";"2" <0pt>
 \ar @{.} "s";"0" <0pt>
 \ar @{.} "s";"1" <0pt>
\endxy} \Ea
, \ \
\text{the dual of $\theta_{ij}$} \rar
 \Ba{c}\resizebox{11mm}{!}{ \xy
(-5,0)*+{_i}*\frm{o}="0";
 (5,0)*+{_j}*\frm{o}="1";
 \ar @{-} "1";"0" <0pt>
\endxy} \Ea,
$$
 trivalent and  bivalent ones. Following our convention in \S {\ref{2: subsec on hypergraphs}}, we call bivalent hyperedges simply {\em edges}\, and
show them in put pictures as {\em solid}\, intervals; the flags of trivalent hyperedges are shown as {\em dotted}\, intervals. For example,
$$
\Ba{c}\resizebox{11mm}{!}{ \xy
(-5,0)*+{_0}*\frm{o}="0";
 (5,0)*+{_1}*\frm{o}="1";
 \ar @{-} "1";"0" <0pt>
\endxy} \Ea
\in \BVGrac((2)),
 \ \ \ \ \
  \Ba{c}\resizebox{10mm}{!}{ \xy
(-4,0)*+{_0}*\frm{o}="0";
 (4,0)*+{_1}*\frm{o}="1";
 (0,5.8)*+{_2}*\frm{o}="2";
\endxy} \Ea
\ , \
\Ba{c}\resizebox{12mm}{!}{ \xy
(-7,0)*+{_0}*\frm{o}="0";
(7,0)*+{_1}*\frm{o}="1";
(0,12)*+{_2}*\frm{o}="2";
  (0,5)*{*}="s";
 \ar @{.} "s";"2" <0pt>
 \ar @{.} "s";"0" <0pt>
 \ar @{.} "s";"1" <0pt>
 \ar @{-} "0";"1" <0pt>
\endxy} \Ea
, \
 \in \BVGrac((3)).
$$
Each flag of a hyperedge is assigned the cohomological degree $-1$, while its $*$-vertex is assigned degree $+1$; each solid edge has degree $-1$.
Every hypergraph in $\BVGrac$ is equipped with an {\em orientation}\, defined as an ordering of its edges, flags, and $*$-vertices (up to a permutation and the corresponding sign). The set of edges (resp., flags) attached to a vertex $\ \wi\ $ of hypergraph $\Ga$ is denoted by $E_i(\Ga)$ (resp., $F_i(\Ga)$).

\sip

The operadic
 composition of such hypergraphs
$$
\Ba{rccc}
_i\circ_j: & \BVGrac((n+1))\ot \BVGrac((m+1)) & \lon & \BVGrac((m+n))\\
           &  \Ga_1 \ot \Ga_2&\lon &  \Ga_1\ \mbox{$_i\circ_j$}\ \Ga_2
\Ea
$$
is defined (cf.\ \S {\ref{3: subsec on Gra_cyc}} and {\ref{4: subsec on HGra}}) by erasing the white vertices $\ \wi\in \Ga_1$ and $\ \wj\in \Ga_2$, and then
(i) reattaching their flags and edges among the vertices of the opponent and (ii) gluing some edges attached to the deleted vertices into new hyperedges. More precisely,
$$
\Ga_1\ \mbox{$_i\circ_j$}\ \Ga_2= \sum_{z,s_1, s_2} \sum_{I'\sqcup I''}\sum_{J'\sqcup J''} \sum_{f_i: F_i(\Ga)\rar (V(\Ga_2)\setminus \resizebox{1.8mm}{!}{\wj})\atop
f_j: F_j(\Ga)\rar( V(\Ga_1)\setminus \resizebox{1.8mm}{!}{\wi})}  \Ga(z,s_1,s_2, I'\sqcup I'',J'\sqcup J'', f_i,f_j)
$$
where 
\Bi
\item[(i)] the first sum from the right runs over all possible maps  $f_i: F_i(\Ga)\rar (V(\Ga_2)\setminus \resizebox{3mm}{!}{\wj})$ and
$f_j: F_j(\Ga)\rar( V(\Ga_1)\setminus {\resizebox{3mm}{!}{\wi}})$,

\item[(ii)] the sums $\sum_{I'\sqcup I''}$ and $\sum_{J'\sqcup J''}$ run over all possible decompositions of the sets of edges attached to $\ \wi\ $ and to $\ \wj\ $,
$$
E_i
=I'\sqcup I'',\ \ \  
E_j
=J'\sqcup J'' \ \text{with}\ \ \# I''=\# J''
$$

\item[(iii)] the sum  $\sum_{z,s_i, s_j}$ runs over all possible bijections
$
z: I''\rar J''$
and all possible maps $s_i:I' \rar (V(\Ga_2)\setminus \resizebox{1.8mm}{!}{\wj})$ and $s_j: J'\rar ( V(\Ga_1)\setminus \resizebox{1.8mm}{!}{\wi})$.
\Ei
The hypergraph  
$\Ga(z,s_1,s_2, I'\sqcup I'',J'\sqcup J'', f_i,f_j)$ is obtained from $\Ga_1$ and $\Ga_2$ by erasing their white vertices $\ \wi\in \Ga_1$ and $\ \wj\in \Ga_2$
and
\Bi
\item[(i)] attaching flags from $F_i$ (resp.\ $F_j$) to the white vertices of $\Ga_2$ (resp.\ $\Ga_1$) as prescribed by the map $f_i$ (resp.\ $f_j$),  
\item[(ii)] attaching edges  from $I'$ (resp.\ $J'$) to the white vertices of $\Ga_2$ (resp.\ $\Ga_1$) as prescribed by the map $s_i$ (resp.\  $s_j$)
\item[(iii)] gluing the dangling endpoint of the edge $e\in I''$ to the dangling end point of corresponding edge $z(e)\in J''$  and making them into two flags 
 of a new hyperedge,
$$
 \Ba{c}\resizebox{21mm}{!}{ \xy
(-7,0)*+{}="1";
 (8,0)*+{}="2";
 (-4,3)*{_e};
 (6,3)*{^{z(e)}};
 (0,10)*{}="0";
  (0,4)*{*}="s";
 \ar @{.} "s";"2" <0pt>
 \ar @{.>} "s";"0" <0pt>
 \ar @{.} "s";"1" <0pt>
\endxy} \Ea
\ \ \ \
\text{orientated as}\ e\wedge z(e) \wedge \text{new flag}\wedge  \text{$*$-vertex} ,
$$
and then taking a sum over attachments of the third ``hanging" flag to the remaining white vertices of $\Ga_1$ and $\Ga_2$ in all possible ways.
\Ei

For example,
$$
  \Ba{c}\resizebox{11mm}{!}{ \xy
(-5,0)*+{_1}*\frm{o}="1";
 (5,0)*+{_2}*\frm{o}="2";
 (0,8)*+{_0}*\frm{o}="0";
  \ar @{-} "1";"2" <0pt>
\endxy} \Ea
\ \mbox{$_2\circ_0$} \
 \Ba{c}\resizebox{11mm}{!}{ \xy
(-5,0)*+{_0}*\frm{o}="0";
 (5,0)*+{_1}*\frm{o}="1";
 (0,8)*+{_2}*\frm{o}="2";
  \ar @{-} "1";"2" <0pt>
  \ar @{-} "0";"1" <0pt>
\endxy} \Ea
=
 \Ba{c}\resizebox{12mm}{!}{ \xy
(-5,0)*+{_1}*\frm{o}="1";
 (5,0)*+{_2}*\frm{o}="2";
 (-5,10)*+{_0}*\frm{o}="0";
  (5,10)*+{_3}*\frm{o}="3";
 \ar @{-} "1";"3" <0pt>
  \ar @{-} "0";"2" <0pt>
 \ar @{-} "3";"2" <0pt>
\endxy} \Ea
+
 \Ba{c}\resizebox{12mm}{!}{ \xy
(-5,0)*+{_1}*\frm{o}="1";
 (5,0)*+{_2}*\frm{o}="2";
 (-5,10)*+{_0}*\frm{o}="0";
  (5,10)*+{_3}*\frm{o}="3";
  (0,5)*{*}="s";
 \ar @{.} "s";"0" <0pt>
 \ar @{.} "s";"1" <0pt>
  \ar @{.} "s";"2" <0pt>
  \ar @{-} "3";"2" <0pt>
\endxy} \Ea
+
 \Ba{c}\resizebox{12mm}{!}{ \xy
(-5,0)*+{_1}*\frm{o}="1";
 (5,0)*+{_2}*\frm{o}="2";
 (-5,10)*+{_0}*\frm{o}="0";
  (5,10)*+{_3}*\frm{o}="3";
  (0,5)*{*}="s";
 \ar @{.} "s";"3" <0pt>
 \ar @{.} "s";"1" <0pt>
  \ar @{.} "s";"2" <0pt>
  \ar @{-} "3";"2" <0pt>
\endxy} \Ea.
$$

\sip

The operad $\Grac_2$ is not a suboperad of $\BVGrac$; nevertheless we still have a well-defined morphism of
cyclic operads
\Beq\label{5: f from BV to BVgrac}
f: \BVc \lon \BVGrac
\Eeq
given on the generators by the same formulae as before,
\Beq\label{5: f from BV to BVgrac formulae}
\Ba{c}\resizebox{9mm}{!}{\xy
 (0,1)*{_{_0}}="a",
(7,1)*{_{_1}}="b",
(3.5,0.6)*{^{_\Delta}},
\ar @{-} "a";"b" <0pt>
\endxy
}\Ea \rar  \resizebox{11mm}{!}{ \xy
(-5,1)*+{_0}*\frm{o}="0";
 (5,1)*+{_1}*\frm{o}="1";
 \ar @{-} "1";"0" <0pt>
\endxy}  , \ \ \
\Ba{c}\begin{xy}
 <0mm,0.66mm>*{};<0mm,3mm>*{}**@{-},
 <0.39mm,-0.39mm>*{};<2.2mm,-2.2mm>*{}**@{-},
 <-0.35mm,-0.35mm>*{};<-2.2mm,-2.2mm>*{}**@{-},
 <0mm,0mm>*{\circ};<0mm,0mm>*{}**@{},
   <0mm,0.66mm>*{};<0mm,3.4mm>*{^{_0}}**@{},
   <0.39mm,-0.39mm>*{};<2.9mm,-4mm>*{^{_2}}**@{},
   <-0.35mm,-0.35mm>*{};<-2.8mm,-4mm>*{^{_1}}**@{},
\end{xy}\Ea
 \rar  \Ba{c}\resizebox{10mm}{!}{ \xy
(-4,0)*+{_0}*\frm{o}="0";
 (4,0)*+{_1}*\frm{o}="1";
 (0,5.8)*+{_2}*\frm{o}="2";
\endxy} \Ea.
\Eeq

\subsection{A twisted extension of  $\BVGrac$}
The operad $\BVGrac$ comes equipped with a morphism 
of cyclic operads
$$
f: \Lie_3 \lon \BVGrac
$$
given on the generator by (cf.\ Lemma {\ref{3: f from Lie_d to Gra_d}}) 
$$
\Ba{c}\begin{xy}
 <0mm,0.66mm>*{};<0mm,3mm>*{}**@{-},
 <0.39mm,-0.39mm>*{};<2.2mm,-2.2mm>*{}**@{-},
 <-0.35mm,-0.35mm>*{};<-2.2mm,-2.2mm>*{}**@{-},
 <0mm,0mm>*{\bu};<0mm,0mm>*{}**@{},
   <0mm,0.66mm>*{};<0mm,3.4mm>*{^{_0}}**@{},
   <0.39mm,-0.39mm>*{};<2.9mm,-4mm>*{^{_2}}**@{},
   <-0.35mm,-0.35mm>*{};<-2.8mm,-4mm>*{^{_1}}**@{},
\end{xy}\Ea
\lon
 \Ba{c}\resizebox{13mm}{!}{ \xy
(-6,0)*+{_1}*\frm{o}="1";
 (6,0)*+{_2}*\frm{o}="2";
 (0,10)*+{_0}*\frm{o}="0";
  (0,4)*{*}="s";
 \ar @{-} "s";"2" <0pt>
 \ar @{-} "s";"0" <0pt>
 \ar @{-} "s";"1" <0pt>
\endxy} \Ea.
$$
Hence the operad $\BVGrac$ can be twisted (in a close analogy to the twisting of $\HGra$ in \S {\ref{3: subsec on twsisting of HGra}})  into a dg cyclic operad
$$
f\BVGcyc=\left\{f\BVGcyc((n+1))\right\}_{n\geq 1}
$$
using a cyclic version of Thomas Willwacher's  twisting endofunctor \cite{Wi} .

\sip

 By definition, $f\BVGcyc((n+1))$
is an $\bS_{(n+1)}$-module spanned over a field $\K$ by the so called {\it admissible}\, hypergraphs $\Ga$  which have two types of hyperedges (dotted and solid ones) precisely as in the case of $\BVGrac$ and  also two types of vertices (black and white ones); here are a few examples,
$$
\Ba{c}\resizebox{13mm}{!}{ \xy
(-5,0)*+{_0}*\frm{o}="0";
 (5,0)*+{_1}*\frm{o}="1";
 \ar @{-} "1";"0" <0pt>
\endxy} \Ea
, \
\Ba{c}\resizebox{16mm}{!}{ \xy
(-7,0)*+{_0}*\frm{o}="0";
(7,0)*+{_1}*\frm{o}="1";
(0,6.5)*{\bu}="2";
  (0,0)*{*}="s";
 \ar @{.} "s";"2" <0pt>
 \ar @{.} "s";"0" <0pt>
 \ar @{.} "s";"1" <0pt>
\endxy} \Ea
\in f\BVGcyc(2),
 \ 
  \Ba{c}\resizebox{13mm}{!}{ \xy
(-4,0)*+{_0}*\frm{o}="0";
 (4,0)*+{_1}*\frm{o}="1";
 (0,5.8)*+{_2}*\frm{o}="2";
\endxy} \Ea
\ , \
\Ba{c}\resizebox{16mm}{!}{ \xy
(-7,0)*+{_0}*\frm{o}="0";
(7,0)*+{_1}*\frm{o}="1";
(0,12)*+{_2}*\frm{o}="2";
  (0,5)*{*}="s";
 \ar @{.} "s";"2" <0pt>
 \ar @{.} "s";"0" <0pt>
 \ar @{.} "s";"1" <0pt>
 \ar @{-} "0";"1" <0pt>
\endxy} \Ea
, \
\Ba{c}\resizebox{19mm}{!}{ \xy
(-7,0)*+{_0}*\frm{o}="0";
 (9,0)*+{_2}*\frm{o}="2";
 (-7,10)*+{_1}*\frm{o}="1";
  (7,10)*{\bu}="3";
  (-3.5,5)*{*}="s1";
  (3,3)*{*}="s2";
 \ar @{.} "s1";"1" <0pt>
 \ar @{.} "s1";"0" <0pt>
  \ar @{.} "s1";"3" <0pt>
\ar @{.} "s2";"2" <0pt>
\ar @{.} "s2";"3" <0pt>
\ar @{.} "s2";"1" <0pt>
\ar @{-} "3";"2" <0pt>
\ar @{-} "2";"0" <0pt>
\endxy} \Ea \in f\BVGcyc(3)
$$

 More precisely, a hypergraph $\Ga$ is called {\em admissible}\, if
\Bi
\item[(i)] $\Ga$ has precisely $n+1$ white vertices labelled by integers from
the set $((n+1))$; they are assigned the cohomological degree zero.
\item[(ii)] $\Ga$ may have black vertices which are unlabelled; they are assigned the cohomological degree
$3$ and their set is denoted by $V_\bu(\Ga)$.
\item[(iii)] $\Ga$ may have {\it solid edges}
 connecting pairs of different white vertices, or pairs consisting of a white vertex and a black connecting a white vertex; {\em no edges between black vertices are allowed}\footnote{The condition ``no edges between black vertices" is better understood as follows: consider
for a moment a larger cyclic operad $ff\BVGcyc$ in which such edges are allowed, and then notice that hypergraphs with edges
between black vertices form an operadic ideal. The quotient by this ideal is precisely $f\BVGcyc$.}; the set of edges
of $\Ga$ is denoted by $E(\Ga)$; each edge has degree $-1$.
\item[(iv)] $\Ga$ may have hyperedges of exactly the same type as in the case of $\BVGrac$,
$$
\Ba{c}\resizebox{10mm}{!}{ \xy
(-4,0)*{}="0";
 (4,0)*{}="1";
 (0,6)*{}="2";
  (0,2)*{*}="s";
 \ar @{.} "s";"2" <0pt>
 \ar @{.} "s";"0" <0pt>
 \ar @{.} "s";"1" <0pt>
\endxy} \Ea\ .
$$
Every hyperedge is required to have at least one flag connected to a white vertex,  i.e.\  hyperedges of the form
\Beq\label{7: wrong hyperedges}
 \Ba{c}\resizebox{12mm}{!}{ \xy
(-6,0)*{\bu}="1";
 (6,0)*{\bu}="2";
 (0,10)*{\bu}="0";
  (0,4)*{*}="s";
 \ar @{.} "s";"2" <0pt>
 \ar @{.} "s";"0" <0pt>
 \ar @{.} "s";"1" <0pt>
\endxy}\Ea
\Eeq
are not allowed (it is better to think that we take a quotient by the operadic ideal generated by hypergraphs with at least one hyperedge of this form).
The set of hyperedges of $\Ga$  is denoted by $V_*(\Ga)$; the star vertex has degree +1, while its flags have degree $-1$ each.

\item[(v)] Every black vertex of $\Ga$ (if any) can be connected to at least two white vertices by a continuous  path made of edges or flags.
\item[(vi)] Every admissible hypergraph $\Ga$ is equipped with an {\em orientation}\, $or(\Ga)$ defined as an ordering of its edges, flags, $*$-vertices and black vertices up to a permutation and the corresponding sign, i.e.
    $$
    or(\Ga)\in \det E(\Ga) \ot \det(Fl(\Ga)) \ot \det V_*(\Ga) \ot \det V_\bu(\Ga).
    $$

\Ei
 Thus the cohomological degree of $\Ga$ is given by
$$
|\Ga|=3V_{\bu}(\Ga) - 2V_*(\Ga) - E(\Ga).
$$
Note that we put no restrictions on the valencies of vertices of admissible hypergraphs $\Ga$ yet.

\sip

The cyclic operadic structure in $f\BVGcyc$ is an obvious extension of the one in $\BVGrac$;
the only difference is that under the compositions $\Ga_1\ \mbox{$_i\circ_j$}\ \Ga_2$ the flags attached to,
say, white vertex $\wi \in \Ga_1$ can now hit not only white vertices in $\Ga_2$ but also allowed black vertices  (i.e.\ the ones which do not belong to prohibited subgraphs (\ref{7: wrong hyperedges})) 
of $\Ga_2$ as well; similarly, the hanging flag arising under gluing of edges
from $\Ga_1$ and $\Ga_2$ into a new hyperedge  hits now all allowed vertices, not only the white ones.


\subsection{The differential}
The cyclic operad $f\BVGcyc$ contains both $\BVGrac$ and $ \HGraphs_3$ as cyclic suboperads. Moreover,
the differential  $\p=\delta_\bu + D$ in $\HGraphs_3$ given by (\ref{4: differential on HGrapgs})
extends naturally to a differential $\p$  in $f\BVGcyc$. However this differential does not respect the  inclusion
$\BVc\rar f\BVGcyc$ as
$$
\p \Ba{c}\resizebox{11mm}{!}{ \xy
(-5,1)*+{_0}*\frm{o}="0";
 (5,1)*+{_1}*\frm{o}="1";
 \ar @{-} "1";"0" <0pt>
\endxy}
 \Ea
 =
 \Ba{c}\resizebox{16mm}{!}{ \xy
(-9,1)*+{_0}*\frm{o}="0";
(9,1)*+{_1}*\frm{o}="1";
(0,7)*{\bu}="b-up";
(0,0)*{*}="s";
 \ar @{-} "0";"b-up" <0pt>
 \ar @{-} "1";"b-up" <0pt>
  \ar @{.} "s";"b-up" <0pt>
 \ar @{.} "s";"1" <0pt>
 \ar @{.} "s";"0" <0pt>
\endxy} \Ea
  \neq 0.
$$
We twist the differential $\p$ by changing the hypergraph $D$ from Lemma {\ref{4: Lemma on MC property of D}} as follows,
$$
D \rar \hat{D}:=
\Ba{c}\resizebox{13mm}{!}{ \xy
(-7,1)*+{_0}*\frm{o}="0";
(-3.5,-1)*{_1};
(3.5,-1)*{_2};
(0,3)*{^3};
(7,1)*+{_1}*\frm{o}="1";
(0,1)*{\bu};
 \ar @{-} "0";"1" <0pt>
\endxy} \Ea
+\hspace{-2mm}
\Ba{c}\resizebox{13mm}{!}{ \xy
(-7,0)*+{_0}*\frm{o}="0";
(7,0)*+{_1}*\frm{o}="1";
(0,12)*{\bu}="2";
  (0,5)*{*}="s";
(-3.3,4.3)*{_1};
(3.3,4.3)*{_2};
(1.4,7.5)*{^3};
(0,2.6)*{_5};
(0,14)*{^4};
 \ar @{.} "s";"2" <0pt>
 \ar @{.} "s";"0" <0pt>
 \ar @{.} "s";"1" <0pt>
\endxy} \Ea
$$
The orientations are defined as $1\wedge 2\wedge 3$ in the first summand and as
$1\wedge 2\wedge 3\wedge 4 \wedge 5$ in the second summand. This hypergraph satsifies $(01)\hat{D}=-\hat{D}$
and hence defines a derivation of $f\BVGcyc$  indeed.

\subsubsection{\bf Lemma}\label{5: Lemma on MC property of hat D} {\em The hypergraph $\hat{D}\in f\BVGcyc$ satisfies the equation}
$$
\delta_{\bu} \hat{D} + \hat{D}\ \mbox{${_1\circ_0}$}\ \hat{D}=0
$$
\begin{proof} In view of Lemma {\ref{4: Lemma on MC property of D}}, it is enough to show that
$$
\delta_\bu
\Ba{c}\resizebox{13mm}{!}{ \xy
(-7,1)*+{_0}*\frm{o}="0";
(7,1)*+{_1}*\frm{o}="1";
(0,1)*{\bu};
 \ar @{-} "0";"1" <0pt>
\endxy} \Ea
+
\Ba{c}\resizebox{13mm}{!}{ \xy
(-7,1)*+{_0}*\frm{o}="0";
(7,1)*+{_1}*\frm{o}="1";
(0,1)*{\bu};
 \ar @{-} "0";"1" <0pt>
\endxy} \Ea
\ \mbox{${_1\circ_0}$}\
\Ba{c}\resizebox{13mm}{!}{ \xy
(-7,1)*+{_0}*\frm{o}="0";
(7,1)*+{_1}*\frm{o}="1";
(0,1)*{\bu};
 \ar @{-} "0";"1" <0pt>
\endxy} \Ea
+
\left(\Id - (01)\right)\left(
\Ba{c}\resizebox{13mm}{!}{ \xy
(-7,1)*+{_0}*\frm{o}="0";
(7,1)*+{_1}*\frm{o}="1";
(0,1)*{\bu};
 \ar @{-} "0";"1" <0pt>
\endxy} \Ea
\ \mbox{${_1\circ_0}$}\
\hspace{-2mm}
\Ba{c}\resizebox{13mm}{!}{ \xy
(-7,0)*+{_0}*\frm{o}="0";
(7,0)*+{_1}*\frm{o}="1";
(0,12)*{\bu}="2";
  (0,5)*{*}="s";
 \ar @{.} "s";"2" <0pt>
 \ar @{.} "s";"0" <0pt>
 \ar @{.} "s";"1" <0pt>
\endxy} \Ea
\right)
=0.
$$
where the symbol (01) stands for the permutation acting on the white vertices of the hypergraphs standing to the right of it.
We have, using hypergraphs with identical orientations,
$$
\delta_\bu\Ba{c}\resizebox{13mm}{!}{ \xy
(-7,1)*+{_0}*\frm{o}="0";
(7,1)*+{_1}*\frm{o}="1";
(0,1)*{\bu};
 \ar @{-} "0";"1" <0pt>
\endxy} \Ea
=
\Ba{c}\resizebox{16mm}{!}{ \xy
(-9,1)*+{_0}*\frm{o}="0";
(9,1)*+{_1}*\frm{o}="1";
(0,7)*{\bu}="b-up";
(0,-7)*{\bu}="b-d";
(0,0)*{*}="s";
 \ar @{-} "0";"b-up" <0pt>
 \ar @{-} "1";"b-d" <0pt>
  \ar @{.} "s";"b-up" <0pt>
 \ar @{.} "s";"b-d" <0pt>
 \ar @{.} "s";"0" <0pt>
\endxy} \Ea
+
\Ba{c}\resizebox{16mm}{!}{ \xy
(-9,1)*+{_0}*\frm{o}="0";
(9,1)*+{_1}*\frm{o}="1";
(0,7)*{\bu}="b-up";
(0,-7)*{\bu}="b-d";
(0,0)*{*}="s";
 \ar @{-} "0";"b-up" <0pt>
 \ar @{-} "1";"b-d" <0pt>
  \ar @{.} "s";"b-up" <0pt>
 \ar @{.} "s";"b-d" <0pt>
 \ar @{.} "s";"1" <0pt>
\endxy} \Ea,
$$
$$
\Ba{c}\resizebox{13mm}{!}{ \xy
(-7,1)*+{_0}*\frm{o}="0";
(7,1)*+{_1}*\frm{o}="1";
(0,1)*{\bu};
 \ar @{-} "0";"1" <0pt>
\endxy} \Ea
\ \mbox{${_1\circ_0}$}\
\Ba{c}\resizebox{13mm}{!}{ \xy
(-7,1)*+{_0}*\frm{o}="0";
(7,1)*+{_1}*\frm{o}="1";
(0,1)*{\bu};
 \ar @{-} "0";"1" <0pt>
\endxy} \Ea=
\underbrace{
\Ba{c}\resizebox{13mm}{!}{ \xy
(-7,0)*+{_0}*\frm{o}="0";
(7,0)*+{_1}*\frm{o}="1";
(0,5)*{\bu}="u";
(0,-5)*{\bu}="d";
 \ar @{-} "0";"u" <0pt>
 \ar @{-} "0";"d" <0pt>
 \ar @{-} "1";"u" <0pt>
 \ar @{-} "1";"d" <0pt>
\endxy} \Ea}_{\equiv 0}
  -
\Ba{c}\resizebox{16mm}{!}{ \xy
(-9,1)*+{_0}*\frm{o}="0";
(9,1)*+{_1}*\frm{o}="1";
(0,7)*{\bu}="b-up";
(0,-7)*{\bu}="b-d";
(0,0)*{*}="s";
 \ar @{-} "0";"b-up" <0pt>
 \ar @{-} "1";"b-d" <0pt>
  \ar @{.} "s";"b-up" <0pt>
 \ar @{.} "s";"b-d" <0pt>
 \ar @{.} "s";"0" <0pt>
\endxy} \Ea
-
\Ba{c}\resizebox{16mm}{!}{ \xy
(-9,1)*+{_0}*\frm{o}="0";
(9,1)*+{_1}*\frm{o}="1";
(0,7)*{\bu}="b-up";
(0,-7)*{\bu}="b-d";
(0,0)*{*}="s";
 \ar @{-} "0";"b-up" <0pt>
 \ar @{-} "1";"b-d" <0pt>
  \ar @{.} "s";"b-up" <0pt>
 \ar @{.} "s";"b-d" <0pt>
 \ar @{.} "s";"1" <0pt>
\endxy} \Ea,
$$
$$
(\Id - (01))\left(\Ba{c}\resizebox{13mm}{!}{ \xy
(-7,0)*+{_0}*\frm{o}="0";
(7,0)*+{_1}*\frm{o}="1";
(0,12)*{\bu}="2";
  (0,5)*{*}="s";
 \ar @{.} "s";"2" <0pt>
 \ar @{.} "s";"0" <0pt>
 \ar @{.} "s";"1" <0pt>
\endxy} \Ea
\ \mbox{${_1\circ_0}$}\
\Ba{c}\resizebox{13mm}{!}{ \xy
(-7,1)*+{_0}*\frm{o}="0";
(7,1)*+{_1}*\frm{o}="1";
(0,1)*{\bu};
 \ar @{-} "0";"1" <0pt>
\endxy} \Ea\right)
=
(\Id - (01))\Ba{c}\resizebox{14mm}{!}{ \xy
(-7,0)*+{_0}*\frm{o}="0";
(7,0)*+{_1}*\frm{o}="1";
(0,12)*{\bu}="2";
  (0,6)*{*}="s";
(0,0)*{\bu};
 \ar @{.} "s";"2" <0pt>
 \ar @{.} "s";"0" <0pt>
 \ar @{.} "s";"1" <0pt>
 \ar @{-} "0";"1" <0pt>
\endxy} \Ea
\equiv 0,
$$

 In last three lines we used the assumption that the admissible hypergraphs have no solid edges between black vertices (and this is the main reason to make that assumption). The above formulae imply the Proposition.
\end{proof}

Hence the degree 1 derivation on the cyclic operad $f\BVGcyc$ given by
\Beq\label{5: full differential as a sum in BVGraphs}
\delta \Ga= \delta_{\bu} \Ga + \hat{D}(\Ga)
\Eeq
is a differential in $f\BVGcyc$. Let us describe it in more detail:  $\delta$ is fully determined by its action on black and white vertices of a hypergraph $\Ga\in f\BVGcyc((n+1))$,
\Beq\label{full delta in BVGraphs}
\delta\Ga=\underbrace{\sum_{v\in \sV_\bu(\Ga)}\delta_{v}\Ga}_{\delta_\bu \Ga} + \underbrace{\sum_{i=0}^n \delta_i^*\Ga}_{\delta_* \Ga}  + \underbrace{\sum_{i=0}^n \delta_i^{(1)}\Ga}_{\delta^{(1)}\Ga}  +  \underbrace{\sum_{i=0}^n \delta_i^{(2)})\Ga}_{\delta^{(2)}\Ga}
\Eeq
where
\Bi
\item[(i)] $\delta_{v}$ acts on the vertex $v$ in the complete analogy to (\ref{4: d on black in hGC}),
i.e.\ it first splits $v$ as follows
\Beq\label{5: d on black in BVGraphs}
\delta_v: \bu \ \lon  \  \frac{1}{2}
\Ba{c}\resizebox{10mm}{!}{ \xy
(0,-4)*{\bu}="1";
 (0,4)*{\bu}="2";
 (8,0)*{}="0";
  (3,0)*{\ast}="s";
 \ar @{.} "s";"2" <0pt>
 \ar @{.>} "s";"0" <0pt>
 \ar @{.} "s";"1" <0pt>
\endxy} \Ea,
\Eeq
and then redistributes the flags and edges attached earlier to $v$ among the newly created two black
vertices in all possible ways\footnote{In this and other pictures  the newly created ``hanging" flag or the ``hanging" edge means a sum over all its possible legitimite attachments to other vertices of $\Ga$.}. 

\item[(ii)]
The operation $\delta_i^*$ splits the white vertex  $\ \wi$
$$
\delta_i^*:\  \wi \lon
\Ba{c}\resizebox{10mm}{!}{
\xy
(0,-4.3)*+{_i}*\frm{o}="1";
 (0,4)*{\bu}="2";
 (8,0)*{}="0";
  (3,0)*{\ast}="s";
 \ar @{.} "s";"2" <0pt>
 \ar @{.>} "s";"0" <0pt>
 \ar @{.} "s";"1" <0pt>
\endxy} \Ea
$$
and redistributes flags and edges attached earlier to $\ \wi\ $ in all possible  ways among the
 among $\ \wi\ $ and the newly created black vertex.

\item[(iii)]
The operation $\delta_i^{(2)}$ (the superscript $(2)$ indicates that it creates two solid edges) acts  on the white vertex  $\ \wi\ $
by splitting  $\ \wi\ $ as follows
$$
\delta_i^{(2)}:\ \wi \lon
\Ba{c}\resizebox{9mm}{!}{ \xy
(0,-3)*+{_i}*\frm{o}="0";
 (0,3)*{\bu}="1";
 (6,3)*{}="2";
 \ar @{-} "1";"0" <0pt>
 \ar @{->} "1";"2" <0pt>
\endxy}\Ea.
$$
and then runs a summation over all possible legitimate reattachments  of edges and flags attached earlier to $\ \wi\ $  among $\ \wi\ $ and the newly created black vertex.

\item[(iv)] the operation $\delta_i^{(1)}$ acts non-trivially on $\wi$ if only if $\wi$ has at least one solid edge attached; if $E_i$ stands for the set of solid edges attached to  $\ \wi\ $, we have
$$
\delta_i^{(1)}=\sum_{e\in E_i} \delta_{i,e}^{(1)},
$$
where the operator $\delta_{i,e}^{(1)}$ is given pictorially by
$$
\delta^{(1)}_{i,e}:\
\resizebox{4.5mm}{!}{ \xy
(0,-3)*+{_i}*\frm{o}="0";
 (0,3)*{\odot}="1";
 (0,4.7)*{^a}="2";
 (-1,0.5)*{_e};
 \ar @{-} "1";"0" <0pt>
\endxy} \ \ 
\lon \ \ 
\resizebox{15mm}{!}{ \xy
(0,-3)*+{_i}*\frm{o}="0";
 (0,3)*{\bu}="1";
 (6,3)*{*}="2";
 (6,9)*{}="u";
 (6,1.4)*{_e};
 (12,3)*{\odot}="3";
 (12,5)*{^a};
 \ar @{-} "1";"0" <0pt>
 \ar @{.} "1";"2" <0pt>
 \ar @{.} "1";"3" <0pt>
 \ar @{.>} "2";"u" <0pt>
\endxy}\ ;
$$
it acts in more details as follows: 
\Bi
\item[(i)] $\delta^{(1)}_{i,e}$ substitutes into the white vertex\  $\wi$ \ the graph $\resizebox{4mm}{!}{ \xy
(0,1)*+{_i}*\frm{o}="0";
 (0,6)*{\bu}="1";
 \ar @{-} "1";"0" <0pt>
\endxy}$\ \ and takes the sum over all possible  reattachments  of all edges (except $e$!) and flags attached earlier to $\ \wi\ $ among $\ \wi\ $ and the {\em newly}\, created black vertex;
\item[(ii)] attaches $e$, for a moment, to the newly created black vertex and then changes it immediately to a hyperedge with one ``hanging" flag,
\item[(iii)] takes the sum over attachments of  that hanging flag (shown in the picture above as a dotted arrow)
to all allowed vertices of $\Ga$ (the vertex $\ \wi\ $ including). 

\Ei

\Ei

\sip

For example,
\Beqrn
\delta \Ba{c}\resizebox{11mm}{!}{ \xy
(-5,1)*+{_0}*\frm{o}="0";
 (5,1)*+{_1}*\frm{o}="1";
 \ar @{-} "1";"0" <0pt>
\endxy}
 \Ea
 &=&
 (\Id + (01))\left(
 \Ba{c}\resizebox{16mm}{!}{ \xy
(-9,1)*+{_0}*\frm{o}="0";
(9,1)*+{_1}*\frm{o}="1";
(0,7)*{\bu}="b-up";
 \ar @{-} "0";"b-up" <0pt>
 \ar @{-} "1";"b-up" <0pt>
  \ar @{-} "0";"1" <0pt>
\endxy} \Ea
+ 
 \Ba{c}\resizebox{16mm}{!}{ \xy
(-9,1)*+{_0}*\frm{o}="0";
(9,1)*+{_1}*\frm{o}="1";
(0,7)*{\bu}="b-up";
(0,0)*{*}="s";
 \ar @{-} "1";"b-up" <0pt>
  \ar @{.} "s";"b-up" <0pt>
 \ar @{.} "s";"1" <0pt>
 \ar @{.} "s";"0" <0pt>
\endxy} \Ea
-
 \Ba{c}\resizebox{16mm}{!}{ \xy
(-9,1)*+{_0}*\frm{o}="0";
(9,1)*+{_1}*\frm{o}="1";
(0,7)*{\bu}="b-up";
(0,0)*{*}="s";
 \ar @{-} "0";"b-up" <0pt>
  \ar @{.} "s";"b-up" <0pt>
 \ar @{.} "s";"1" <0pt>
 \ar @{.} "s";"0" <0pt>
\endxy} \Ea\right)
= 0,\\
\delta \left( \Ba{c}\resizebox{9mm}{!}{ \xy
(-4,0)*+{_0}*\frm{o}="0";
 (4,0)*+{_1}*\frm{o}="1";
 (0,5.8)*+{_2}*\frm{o}="2";
\endxy} \Ea\right) &= & 0,
\Eeqrn
implying the following

\subsection{\bf Proposition} {\em There is a morphism of dg cyclic operads
\Beq\label{5: Main Morphism F}
F: \BVc \lon f\BVGcyc
\Eeq
given on the generators by the formulae (\ref{5: f from BV to BVgrac formulae}). }

\mip

This morphism is not a quasi-isomorphism, the right hand side has lots of ``wrong" cohomology classes which we discuss in detail below. 

\subsection{On valencies of black vertices} It is obvious that the subspace
of $f\BVGcyc$ generated by hypergraphs with all black vertices having valency $\geq 3$ forms a cyclic sub-operad. A remarkable fact is that it closed 
under the differential.

\subsubsection{\bf Lemma} {\it  The subspace
of $f\BVGcyc$ generated by hypergraphs with all black vertices (if any) at least trivalent is closed under the differential $\delta$.}

\begin{proof} Assume  $\Ga\in \BVGcyc$ has no black vertices of valency $<3$. Let us prove that $\delta\Ga$ also has no black vertices of valency $<3$.

\sip
 Arguing as in Lemma {\ref{4: Lemma on univalent in hGC}} it is easy to see that  $\delta\Ga$ can not contain univalent black vertices.

\sip

Due to a calculation made in the proof of Theorem {\ref{4: trivalency theorem in hGC}}, $\delta\Ga$ can not contain  bivalent black vertices
with with no edges attached (such vertices can only come in pairs which cancel each other). Similarly, $\delta\Ga$ can not not contain 2-valent black vertices with {\it only}\, edges attached --- such vertices can only be of the form
$$
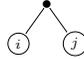
\Ba{c}\resizebox{11mm}{!}{ \xy
(-5,0)*+{_i}*\frm{o}="0";
(5,0)*+{_j}*\frm{o}="1";
 (0,7)*{\bu}="2";
 \ar @{-} "2";"0" <0pt>
 \ar @{-} "1";"2" <0pt>
\endxy}\Ea
$$
and hence they can only  originate from the part $\delta^{(2)}$
of the full differential; hence such bivalent vertices  always come  in pairs, one from $\delta^{(2)}_i \, \wi$\ and one from $\delta^{(2)}_j\,\wj$\ , which cancel each other.

\sip

It remains to check $\delta\Ga$ does not contain bivalent black vertices $v$
with precisely one flag and precisely one edge attached. The edge connects $v$ to some white vertex,
say $\wi$\ . Such vertices can appear as summands in
$$
\delta^{(1)}_{i,e}:\
\resizebox{4.5mm}{!}{ \xy
(0,-3)*+{_i}*\frm{o}="0";
 (0,3)*{\odot}="1";
 (0,4.7)*{^a}="2";
 (-1,0.5)*{_e};
 \ar @{-} "1";"0" <0pt>
\endxy} \ \ 
\lon \ \ 
\resizebox{15mm}{!}{ \xy
(0,-3)*+{_i}*\frm{o}="0";
 (0,3)*{\bu}="1";
 (6,3)*{*}="2";
 (6,9)*{}="u";
 (6,1.4)*{_e};
 (12,3)*{\odot}="3";
 (12,5)*{^a};
 \ar @{-} "1";"0" <0pt>
 \ar @{.} "1";"2" <0pt>
 \ar @{.} "1";"3" <0pt>
 \ar @{.>} "2";"u" <0pt>
\endxy}\ .
$$
If the vertex $a$ is  white, say $a=\wj$\ , then such terms come always in pairs,
$$
(\delta^{(2)}_i + \delta^*_j) (...\Ba{c}\resizebox{11mm}{!}{ \xy
(-5,1)*+{_i}*\frm{o}="0";
 (5,1)*+{_j}*\frm{o}="1";
 \ar @{-} "1";"0" <0pt>
\endxy}\Ea ...)
=
\ldots + \Ba{c}\resizebox{13mm}{!}{ \xy
(0,-3)*+{_i}*\frm{o}="0";
(12,-3)*+{_j}*\frm{o}="j";
 (0,3)*{\bu}="1";
 (6,3)*{*}="2";
  (12,3)*{}="3";
  (6,9)*{}="4";
 \ar @{-} "1";"0" <0pt>
 \ar @{.} "1";"2" <0pt>
  \ar @{.} "2";"j" <0pt>
  \ar @{.} "2";"4" <0pt>
\endxy}\Ea
-
\Ba{c}\resizebox{13mm}{!}{ \xy
(0,-3)*+{_i}*\frm{o}="0";
(12,-3)*+{_j}*\frm{o}="j";
 (0,3)*{\bu}="1";
 (6,3)*{*}="2";
  (12,3)*{}="3";
  (6,9)*{}="4";
 \ar @{-} "1";"0" <0pt>
 \ar @{.} "1";"2" <0pt>
  \ar @{.} "2";"j" <0pt>
  \ar @{.} "2";"4" <0pt>
\endxy}\Ea + \ldots
$$
which  cancel each other. If the vertex $a$ is black, $a=\bu$, then such bivalent vertices also cancel out  with  similar terms coming from $\delta_\bu$ applied to $a$.
 The proof is completed.
\end{proof}

\sip

In view of the above Lemma we assume from now on that the dg cyclic operad  $f\BVGcyc$
is generated by hypergraphs {\it with no univalent or bivalent black vertices}.


\bip

{\large
\section{\bf  Dg cooperad $f\BVGcyc^*$  and cyclic Arnold relations }
}

\sip

\subsection{Dual cooperad} Let
$$
f\BVGcyc^*=\{f\BVGcyc((n+1))^*\}_{n\geq 1}
$$
be the dg cyclic cooperad dual to the operad $\BVGcyc$ introduced in the previous section.
Thus $f\BVGcyc((n+1))^*$ is a dg vector space generated by hypergraphs with labeled $n+1$ white vertices, some number of unlabeled black vertices,
such that (i) there are no edges  and no hyperedges connecting {\it solely}\,  black vertices,  (ii) every black vertex is at least trivalent and  can be connected by a continuous path of edges or flags to at least two white vertices. The differential $d$ in $f\BVGcyc((n+1))^*$ consists of four parts (cf.\ (\ref{5: full differential as a sum in BVGraphs}))
\Beq\label{6: d on CoBVGraphs as three parts}
d=d_{\bu\bu} + d_{\bu\circ}  + d_{\bu\circ}^{(1)}  + d_{\bu\circ}^{(2)}
\Eeq
where
\Bi
\item[(i)]
$d_{\bu\bu}$  acts on unordered pairs $(v',v'')$ of black vertices connected by a hyperedge $*$
 by erasing from $*$ its unique flag connected to some other vertex and then
collapsing the remaining two flags together with their vertices  $v'$ and $v''$ into one black vertex,
$$
d_{\bu\bu}:
\Ba{c}\resizebox{22mm}{!}
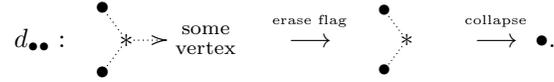
{ \xy
(0,-4)*{\bu}="1";
 (0,4)*{\bu}="2";
 (8,0)*{}="0";
 (13,0)*{\Ba{c}_\text{some}\\ ^\text{vertex}\Ea};
  (3,0)*{\ast}="s";
 \ar @{.} "s";"2" <0pt>
 \ar @{.>} "s";"0" <0pt>
 \ar @{.} "s";"1" <0pt>
\endxy} \Ea
\stackrel{^{\text{erase flag}}}{\lon}
\Ba{c}\resizebox{10mm}{!}{ \xy
(0,-4)*{\bu}="1";
 (0,4)*{\bu}="2";
 (8,0)*{}="0";
  (3,0)*{\ast}="s";
 \ar @{.} "s";"2" <0pt>
 \ar @{.} "s";"1" <0pt>
\endxy} \Ea
\stackrel{^{\text{collapse}}}{\lon} \bu.
$$
More precisely, $d_{\bu\bu}(\Ga)$ is a sum over such triples $(v',v'',*)$.

\mip

\item[(ii)]
$d_{\bu\circ}$ acts on triples $(\ \wi\ ,\bu, *)$ consisting of a white vertex  $\ \wi\ $ and a black vertex $\bu$ which are connected by an hyperedge $*$ as follows: erase the unique flag of $*$ which is not connected to  $\ \wi\ $  or $\bu$ and then
collapse the remaining two flags together with the black vertex $\bu$ onto the white vertex \ $\wi$\ ,
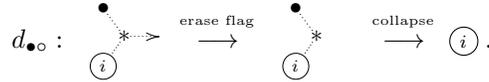
$$
d_{\bu\circ}:
\Ba{c}\resizebox{10mm}{!}{
\xy
(0,-4.3)*+{_i}*\frm{o}="1";
 (0,4)*{\bu}="2";
 (8,0)*{}="0";
  (3,0)*{\ast}="s";
 \ar @{.} "s";"2" <0pt>
 \ar @{.>} "s";"0" <0pt>
 \ar @{.} "s";"1" <0pt>
\endxy} \Ea
\stackrel{^{\text{erase flag}}}{\lon}
\Ba{c}\resizebox{10mm}{!}{
\xy
(0,-4.3)*+{_i}*\frm{o}="1";
 (0,4)*{\bu}="2";
 (8,0)*{}="0";
  (3,0)*{\ast}="s";
 \ar @{.} "s";"2" <0pt>
 \ar @{.} "s";"1" <0pt>
\endxy} \Ea
\stackrel{^{\text{collapse}}}{\lon}
\ \wi\ .
$$
More precisely, $d_{\bu\circ}(\Ga)$ is a sum over such triples $(\ \wi\ ,\bu,*)$.
Note that the action of $d_{\bu\circ}$ on a triple $(\ \wi\ ,\bu,*)$ is set to be zero
if the vertices $\ \wi\ $ and $\bu$ are also connected by an edge (as we do not allow tadpoles
in our hypergraphs).

\mip

\item[(iii)] $d_{\circ\bu}^{(2)}$ acts on triples $(\ \wi\ ,\bu, e)$ consisting of a white vertex and a black vertex $\bu$ connected by an edge $e$ such that the set $E(\bu)$
    of all edges connected to $\bu$ has cardinality $\geq2$. Then 
    $$
    d_{\circ\bu}^{(2)}=\sum_{e'\in E(\bu)\atop{e'\neq e}} d_{\circ\bu}^{e'}
    $$
    where $d_{\circ\bu}^{e'}$ acts on the hypergraph $\Ga$ by erasing the edge $e'$ and then collapsing $\bu$ and $e$ onto the white vertex    $\ \wi\ $,
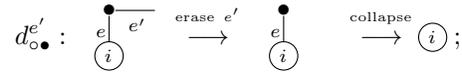
    $$
d_{\circ\bu}^{e'}:
\Ba{c}\resizebox{9mm}{!}{ \xy
(0,-3)*+{_i}*\frm{o}="0";
 (0,3)*{\bu}="1";
 (6,3)*{}="2";
(3.8,1.1)*{_{e'}};
(-1,0)*{_{e}};
 \ar @{-} "1";"0" <0pt>
 \ar @{-} "1";"2" <0pt>
\endxy}\Ea
\stackrel{^{\text{erase}\ e'}}{\lon}
\Ba{c}\resizebox{9mm}{!}{ \xy
(0,-3)*+{_i}*\frm{o}="0";
 (0,3)*{\bu}="1";
 (6,3)*{}="2";
(-1,0)*{_{e}};
 \ar @{-} "1";"0" <0pt>
\endxy}\Ea
\stackrel{^{\text{collapse}}}{\lon} \wi\ ;
$$
More precisely, $d_{\bu\circ}^{(2)}(\Ga)$ is a sum over such triples $(\ \wi\ ,\bu,e)$.

\mip

\item[(iv)]   $d_{\circ\bu}^{(1)}$ acts on 4-tuples $(\ \wi\ ,\bu, e,f)$ consisting of a white vertex connected by an edge $e$ to
 a black vertex $\bu$ ,  and some flag $f$
attached to $\bu$. Let $f'$ and $f''$ be the two other flags of the hyperedge $*$ owning the flag $f$, and  let $w'$ and $w''$ stand for the black or white vertices of $\Ga$ to which flags $f'$ and $f''$ are attached (they can not be both black).
 Then $d_{\circ\bu}^{(1)}$ acts on such a 4-tuple $(\ \wi\ ,\bu, e,f)$ by first collapsing the vertex $\bu$ and the edge $e$ into $\ \wi\ $,
\Beq\label{6: d'' in cooperad, step 1}
    \Ba{c}\resizebox{15mm}{!}{ \xy
(0,-3)*+{_i}*\frm{o}="0";
(12,3)*{\circ}="bu";
(12,4.6)*{^{w'}};
(9,0.9)*{_{f'}};
(3,0.9)*{_{f}};
(8,6.5)*{^{f''}};
(6,12)*{^{w''}};
(-2,0.3)*{^{e}};
 (0,3)*{\bu}="1";
 (6,3)*{*}="2";
  (12,3)*{}="3";
  (6,10)*{\circ}="4";
 \ar @{-} "1";"0" <0pt>
 \ar @{.} "1";"2" <0pt>
  \ar @{.} "2";"bu" <0pt>
  \ar @{.} "2";"4" <0pt>
\endxy}\Ea
\stackrel{^{\text{collaps}\ \bu \&\ e}}{\lon}
   \Ba{c}\resizebox{15mm}{!}{ \xy
(0,-3)*+{_i}*\frm{o}="0";
(-4,10)*{^{w'}};
(-4,5)*{_{f'}};
(-1,0.4)*{_{f}};
(4,5)*{^{f''}};
(4.7,10)*{^{w''}};
 (0,3)*{*}="s";
  (4,8)*{\circ}="l";
  (-4,8)*{\circ}="r";
 \ar @{.} "0";"s" <0pt>
 \ar @{.} "s";"l" <0pt>
  \ar @{.} "s";"r" <0pt>
\endxy}\Ea,
    \Eeq
and then producing two hypergraphs, one by erasing
 $f'$ and another one by  erasing $f''$, and finally replacing
  the remaining dotted interval in each summand by a solid edge as in the picture,
 \Beq\label{6: d'' in cooperad, step 2}
  \Ba{c}\resizebox{15mm}{!}{ \xy
(0,-3)*+{_i}*\frm{o}="0";
(-4,10)*{^{w'}};
(-4,5)*{_{f'}};
(-1,0.4)*{_{f}};
(4,5)*{^{f''}};
(4.7,10)*{^{w''}};
 (0,3)*{*}="s";
  (4,8)*{\circ}="l";
  (-4,8)*{\circ}="r";
 \ar @{.} "0";"s" <0pt>
 \ar @{.} "s";"l" <0pt>
  \ar @{.} "s";"r" <0pt>
\endxy}\Ea
\stackrel{^{\text{erase}\, f'' \pm \text{erase}\, f'}}{\lon}
\Ba{c}\resizebox{4mm}{!}{ \xy
(0,0)*+{_i}*\frm{o}="0";
(0,9.9)*{^{w'}};
(0,4.5)*{*}="s";
 (0,8)*{\circ}="1";
 \ar @{.} "1";"0" <0pt>
\endxy}\Ea
\pm
\Ba{c}\resizebox{5mm}{!}{ \xy
(0,0)*+{_i}*\frm{o}="0";
(0,9.9)*{^{w''}};
(0,4.5)*{*}="s";
 (0,8)*{\circ}="1";
 \ar @{.} "1";"0" <0pt>
\endxy}\Ea
\stackrel{^{\text{replace}}}{\lon}
\Ba{c}\resizebox{4mm}{!}{ \xy
(0,-3)*+{_i}*\frm{o}="0";
(0,4.6)*{^{w'}};
 (0,3)*{\circ}="1";
 \ar @{-} "1";"0" <0pt>
\endxy}\Ea
\pm
\Ba{c}\resizebox{4.7mm}{!}{ \xy
(0,-3)*+{_i}*\frm{o}="0";
(0,4.6)*{^{w''}};
 (0,3)*{\circ}="1";
 \ar @{-} "1";"0" <0pt>
\endxy}\Ea
\Eeq
\Ei

\sip

It is best to understand the action of $d$ on concrete examples.

\subsubsection{\bf Example}\label{5: Ex 1} We have
$$
d\Ba{c}\resizebox{17mm}{!}{ \xy
(-10,0)*+{i}*\frm{o}="0";
(0,0)*+{j}*\frm{o}="1";
(10,0)*+{k}*\frm{o}="2";
  (5,10)*{\bu}="bu";
 (-5,10)*{*}="s";
 \ar @{-} "bu";"2" <0pt>
 \ar @{-} "bu";"0" <0pt>
 \ar @{.} "s";"bu" <0pt>
  \ar @{.} "s";"0" <0pt>
 \ar @{.} "s";"1" <0pt>
\endxy} \Ea
=
\Ba{c}\resizebox{12mm}{!}{ \xy
(-7,0)*+{_i}*\frm{o}="0";
(7,0)*+{_j}*\frm{o}="1";
(0,12)*+{_k}*\frm{o}="2";
  (0,5)*{*}="s";
 \ar @{.} "s";"2" <0pt>
 \ar @{.} "s";"0" <0pt>
 \ar @{.} "s";"1" <0pt>
\endxy} \Ea
-
 \Ba{c}\resizebox{13mm}{!}{ \xy
(-6,0)*+{_i}*\frm{o}="0";
 (6,0)*+{_j}*\frm{o}="1";
 (0,9)*+{_k}*\frm{o}="2";
 \ar @{-} "0";"2" <0pt>
 \ar @{-} "0";"1" <0pt>
\endxy}
\Ea
-
 \Ba{c}\resizebox{13mm}{!}{ \xy
(-6,0)*+{_i}*\frm{o}="0";
 (6,0)*+{_j}*\frm{o}="1";
 (0,9)*+{_k}*\frm{o}="2";
 \ar @{-} "1";"2" <0pt>
 \ar @{-} "1";"0" <0pt>
\endxy}
\Ea
-
 \Ba{c}\resizebox{13mm}{!}{ \xy
(-6,0)*+{_i}*\frm{o}="0";
 (6,0)*+{_j}*\frm{o}="1";
 (0,9)*+{_k}*\frm{o}="2";
 \ar @{-} "2";"0" <0pt>
 \ar @{-} "2";"1" <0pt>
\endxy}
\Ea
$$
Up to $d$-exact term on the l.h.s., this is precisely the combinatorial incarnation od the relation (\ref{5: def of Omega_ijk}) which we implicitly had in mind when defining
the cyclic operad structure on hypergraphs (it comes with the correct signs, i.e.\ the with correct induced orientations of hypergraphs which are skipped in the picture).

\subsubsection{\bf Example}\label{5: Ex 2} We have
\Beqrn
d\Ba{c}\resizebox{19mm}{!}{ \xy
(0,0)*+{i}*\frm{o}="0";
(8,0)*+{j}*\frm{o}="1";
(16,0)*+{k}*\frm{o}="2";
(24,0)*+{l}*\frm{o}="3";
  (18,10)*{\bu}="bu";
 (6,10)*{*}="s";
 \ar @{-} "bu";"2" <0pt>
 \ar @{-} "bu";"3" <0pt>
 \ar @{.} "s";"bu" <0pt>
  \ar @{.} "s";"0" <0pt>
 \ar @{.} "s";"1" <0pt>
\endxy} \Ea
&=&
\underbrace{\Ba{c}\resizebox{11mm}{!}{ \xy
(-5,0)*+{_j}*\frm{o}="1";
 (5,0)*+{_k}*\frm{o}="2";
 (-5,10)*+{_i}*\frm{o}="0";
  (5,10)*+{_l}*\frm{o}="3";
  (0,5)*{*}="s";
 \ar @{.} "s";"0" <0pt>
 \ar @{.} "s";"1" <0pt>
  \ar @{.} "s";"2" <0pt>
\endxy} \Ea
-
 \Ba{c}\resizebox{11mm}{!}{ \xy
(-5,0)*+{_j}*\frm{o}="1";
 (5,0)*+{_k}*\frm{o}="2";
 (-5,10)*+{_i}*\frm{o}="0";
  (5,10)*+{_l}*\frm{o}="3";
  (0,5)*{*}="s";
 \ar @{.} "s";"0" <0pt>
 \ar @{.} "s";"1" <0pt>
  \ar @{.} "s";"3" <0pt>
\endxy} \Ea}_{d^{(2)}_{\bu\circ}}
+
\underbrace{\Ba{c}\resizebox{11mm}{!}{ \xy
(-5,0)*+{_j}*\frm{o}="1";
 (5,0)*+{_k}*\frm{o}="2";
 (-5,10)*+{_i}*\frm{o}="0";
  (5,10)*+{_l}*\frm{o}="3";
 \ar @{-} "0";"2" <0pt>
 \ar @{-} "0";"3" <0pt>
\endxy} \Ea
+
\Ba{c}\resizebox{11mm}{!}{ \xy
(-5,0)*+{_j}*\frm{o}="1";
 (5,0)*+{_k}*\frm{o}="2";
 (-5,10)*+{_i}*\frm{o}="0";
  (5,10)*+{_l}*\frm{o}="3";
 \ar @{-} "1";"2" <0pt>
 \ar @{-} "1";"3" <0pt>
\endxy} \Ea}_{d_{\bu\circ}}\\
&&
\underbrace{\Ba{c}\resizebox{11mm}{!}{ \xy
(-5,0)*+{_j}*\frm{o}="1";
 (5,0)*+{_k}*\frm{o}="2";
 (-5,10)*+{_i}*\frm{o}="0";
  (5,10)*+{_l}*\frm{o}="3";
 \ar @{-} "0";"3" <0pt>
 \ar @{-} "2";"3" <0pt>
\endxy} \Ea
+
\Ba{c}\resizebox{11mm}{!}{ \xy
(-5,0)*+{_j}*\frm{o}="1";
 (5,0)*+{_k}*\frm{o}="2";
 (-5,10)*+{_i}*\frm{o}="0";
  (5,10)*+{_l}*\frm{o}="3";
 \ar @{-} "1";"3" <0pt>
 \ar @{-} "2";"3" <0pt>
\endxy} \Ea
+
\Ba{c}\resizebox{11mm}{!}{ \xy
(-5,0)*+{_j}*\frm{o}="1";
 (5,0)*+{_k}*\frm{o}="2";
 (-5,10)*+{_i}*\frm{o}="0";
  (5,10)*+{_l}*\frm{o}="3";
 \ar @{-} "0";"2" <0pt>
 \ar @{-} "2";"3" <0pt>
\endxy} \Ea
+
\Ba{c}\resizebox{11mm}{!}{ \xy
(-5,0)*+{_j}*\frm{o}="1";
 (5,0)*+{_k}*\frm{o}="2";
 (-5,10)*+{_i}*\frm{o}="0";
  (5,10)*+{_l}*\frm{o}="3";
 \ar @{-} "0";"2" <0pt>
 \ar @{-} "2";"3" <0pt>
\endxy} \Ea
}_{d_{\bu\circ}^{(1)}}\\
\Eeqrn

\subsubsection{\bf Cyclic Arnold relations}
Combining results of the previous two examples we get the equality (we skip showing orientations
of hypergraphs, so the formula below is shown up to signs)
$$
d\left(\Ba{c}\resizebox{19mm}{!}{ \xy
(0,0)*+{i}*\frm{o}="0";
(8,0)*+{j}*\frm{o}="1";
(16,0)*+{k}*\frm{o}="2";
(24,0)*+{l}*\frm{o}="3";
  (18,10)*{\bu}="bu";
 (6,10)*{*}="s";
 \ar @{-} "bu";"2" <0pt>
 \ar @{-} "bu";"3" <0pt>
 \ar @{.} "s";"bu" <0pt>
  \ar @{.} "s";"0" <0pt>
 \ar @{.} "s";"1" <0pt>
\endxy} \Ea
+
\Ba{c}\resizebox{21mm}{!}{ \xy
(-10,0)*+{i}*\frm{o}="0";
(0,0)*+{j}*\frm{o}="1";
(10,0)*+{k}*\frm{o}="2";
(18,0)*+{l}*\frm{o}="3";
  (5,10)*{\bu}="bu";
 (-5,10)*{*}="s";
 \ar @{-} "bu";"2" <0pt>
 \ar @{-} "bu";"0" <0pt>
 \ar @{.} "s";"bu" <0pt>
  \ar @{.} "s";"0" <0pt>
 \ar @{.} "s";"1" <0pt>
\endxy} \Ea
-
\Ba{c}\resizebox{21mm}{!}{ \xy
(-10,0)*+{i}*\frm{o}="0";
(0,0)*+{j}*\frm{o}="1";
(10,0)*+{k}*\frm{o}="2";
(18,0)*+{l}*\frm{o}="3";
  (5,10)*{\bu}="bu";
 (-5,10)*{*}="s";
 \ar @{-} "bu";"3" <0pt>
 \ar @{-} "bu";"0" <0pt>
 \ar @{.} "s";"bu" <0pt>
  \ar @{.} "s";"0" <0pt>
 \ar @{.} "s";"1" <0pt>
\endxy} \Ea
\right)=
$$
$$
+
\Ba{c}\resizebox{11mm}{!}{ \xy
(-5,0)*+{_j}*\frm{o}="1";
 (5,0)*+{_k}*\frm{o}="2";
 (-5,10)*+{_i}*\frm{o}="0";
  (5,10)*+{_l}*\frm{o}="3";
 \ar @{-} "0";"2" <0pt>
 \ar @{-} "0";"3" <0pt>
\endxy} \Ea
+
\Ba{c}\resizebox{11mm}{!}{ \xy
(-5,0)*+{_j}*\frm{o}="1";
 (5,0)*+{_k}*\frm{o}="2";
 (-5,10)*+{_i}*\frm{o}="0";
  (5,10)*+{_l}*\frm{o}="3";
 \ar @{-} "1";"2" <0pt>
 \ar @{-} "1";"3" <0pt>
\endxy} \Ea
+
\Ba{c}\resizebox{11mm}{!}{ \xy
(-5,0)*+{_j}*\frm{o}="1";
 (5,0)*+{_k}*\frm{o}="2";
 (-5,10)*+{_i}*\frm{o}="0";
  (5,10)*+{_l}*\frm{o}="3";
 \ar @{-} "0";"3" <0pt>
 \ar @{-} "2";"3" <0pt>
\endxy} \Ea
+
\Ba{c}\resizebox{11mm}{!}{ \xy
(-5,0)*+{_j}*\frm{o}="1";
 (5,0)*+{_k}*\frm{o}="2";
 (-5,10)*+{_i}*\frm{o}="0";
  (5,10)*+{_l}*\frm{o}="3";
 \ar @{-} "1";"3" <0pt>
 \ar @{-} "2";"3" <0pt>
\endxy} \Ea
+
\Ba{c}\resizebox{11mm}{!}{ \xy
(-5,0)*+{_j}*\frm{o}="1";
 (5,0)*+{_k}*\frm{o}="2";
 (-5,10)*+{_i}*\frm{o}="0";
  (5,10)*+{_l}*\frm{o}="3";
 \ar @{-} "0";"2" <0pt>
 \ar @{-} "2";"3" <0pt>
\endxy} \Ea
+
\Ba{c}\resizebox{11mm}{!}{ \xy
(-5,0)*+{_j}*\frm{o}="1";
 (5,0)*+{_k}*\frm{o}="2";
 (-5,10)*+{_i}*\frm{o}="0";
  (5,10)*+{_l}*\frm{o}="3";
 \ar @{-} "0";"2" <0pt>
 \ar @{-} "2";"3" <0pt>
\endxy} \Ea
$$
$$
+
\Ba{c}\resizebox{11mm}{!}{ \xy
(-5,0)*+{_j}*\frm{o}="1";
 (5,0)*+{_k}*\frm{o}="2";
 (-5,10)*+{_i}*\frm{o}="0";
  (5,10)*+{_l}*\frm{o}="3";
 \ar @{-} "0";"1" <0pt>
 \ar @{-} "0";"2" <0pt>
\endxy} \Ea
+
\Ba{c}\resizebox{11mm}{!}{ \xy
(-5,0)*+{_j}*\frm{o}="1";
 (5,0)*+{_k}*\frm{o}="2";
 (-5,10)*+{_i}*\frm{o}="0";
  (5,10)*+{_l}*\frm{o}="3";
 \ar @{-} "1";"0" <0pt>
 \ar @{-} "1";"2" <0pt>
\endxy} \Ea
+
\Ba{c}\resizebox{11mm}{!}{ \xy
(-5,0)*+{_j}*\frm{o}="1";
 (5,0)*+{_k}*\frm{o}="2";
 (-5,10)*+{_i}*\frm{o}="0";
  (5,10)*+{_l}*\frm{o}="3";
 \ar @{-} "2";"0" <0pt>
 \ar @{-} "2";"1" <0pt>
\endxy} \Ea
+
\Ba{c}\resizebox{11mm}{!}{ \xy
(-5,0)*+{_j}*\frm{o}="1";
 (5,0)*+{_k}*\frm{o}="2";
 (-5,10)*+{_i}*\frm{o}="0";
  (5,10)*+{_l}*\frm{o}="3";
 \ar @{-} "0";"1" <0pt>
 \ar @{-} "0";"3" <0pt>
\endxy} \Ea
+
\Ba{c}\resizebox{11mm}{!}{ \xy
(-5,0)*+{_j}*\frm{o}="1";
 (5,0)*+{_k}*\frm{o}="2";
 (-5,10)*+{_i}*\frm{o}="0";
  (5,10)*+{_l}*\frm{o}="3";
 \ar @{-} "1";"0" <0pt>
 \ar @{-} "1";"3" <0pt>
\endxy} \Ea
+
\Ba{c}\resizebox{11mm}{!}{ \xy
(-5,0)*+{_j}*\frm{o}="1";
 (5,0)*+{_k}*\frm{o}="2";
 (-5,10)*+{_i}*\frm{o}="0";
  (5,10)*+{_l}*\frm{o}="3";
 \ar @{-} "3";"0" <0pt>
 \ar @{-} "3";"1" <0pt>
\endxy} \Ea
$$
which at the cohomology level gives us precisely the cyclic Arnold relations
(\ref{5: Cyclic Arnold for Omega_ijk})-(\ref{5: def of Omega_ijk}) for the combinatorial
incarnations of the generators $\bar{\al}_{ij}$ of $H^\bu(\FFM_2)$.

\bip

{\large 
\section{\bf From $f\BVGcyc^*$ to the semialgebraic De Rham algebra on $\FFM_2$}
}
\sip

\subsection{De Rham algebra on $\FFM_2$}
Let
$$
\Omega_{\FFM_2}=\{ \Omega_{\FFM_2((n+1))}, d_{dR}\}_{n\geq 1}
$$
be the $\bS$-module of dg algebras of semialgebraic differential forms on the compactified
moduli spaces $\FFM_2((n+1))$ \cite{KS, HLTV, GS}; its cohomology $H^\bu(\Omega_{\FFM_2})$  is precisely the cyclic  cooperad
 $(\BVc)^*$.

\sip

For any $i,j\in ((n+1))$ there is a projection
$$
p_{ij}: \FFM_2((n+1))^0\lon \FFM_2((2))
$$
from the big open cell  $\FFM_2((n+1))^0$  of the compactified space
$\FFM_2((n+1))$
which forgets all framed points except the ones labelled by $i$ and $j$. Since $\FFM_2((2))$
can be identified with the circle $S^1$, one can define a semialgebraic 1-form
$$
\al_{ij}=\al_{ji}:=\pi_{ij}^*(\text{Vol}_{S^1}), \ \forall i\neq j, \ i,j\in \{0,1,\ldots,n\},
$$
on $\FFM_2((n+1))^0$ as a pull-back the normalized homogeneous volume form $\text{Vol}_{S^1}$ on $S^1$,
$\int_{S^1} \text{Vol}_{S^1}=1$. This form extends nicely to the reak compactification $\FFM_2((n+1))$
of $\FFM_2((n+1))^0$, and we use it to construct the $S^1$-basic 2-form (\ref{5: def of Omega_ijk})
on $\FFM_2((3))\simeq S^1\times S^1 \times S^1$ and hence on any $\FFM_2((n+1))$ with $n\geq 2$.
These two differential forms, $\al_{ij}$ and $\Omega_{ijk}$, are the propagators which we use
below to build a suitable  morphism from $(f\BVGcyc)^*$ to $\Omega_{\FFM_2}$.

\subsection{Propagators in an affine coordinate system} Using the action of the projective group
$PGL(2,\C)$ on framed points in $\P^1\simeq \C\sqcup {\infty}$ one can fix the point labelled,
say, by $0$ at $\infty$ and fix the tangent ray at that point to be tangent to the real big
circle $\R$ through $\infty$. This trick identifies
\Beq\label{6: fcM via FM}
\FFM_2((n+1))\simeq \mathsf{FM}_2(n)\times (S^1)^n
\Eeq
where $\mathsf{FM}_2(n)$ is the standard Fulton-MacPherson
compacification of the configuration space
$$
\text{Conf}_n(\C):=\frac{i: [n]\hook \C}{\R^+\ltimes \R^2}
$$
of $n$ different points in the complex plane  $\C$ modulo the affine group action. Using this identification
one gets \cite{GS} the following affine representations of the propagators \cite{GS}, 
$$
\al_{ij}= \eta_i + \eta_j - 2\omega_{ij}
$$
and hence
\Beq\label{6: Omega_ijk in affine coordiantes}
\Omega_{ijk}=\eta_i\eta_j +\eta_k\eta_i + \eta_j\eta_k + 2\omega_{ij}(\eta_i -\eta_j)
+   2\omega_{ki}(\eta_k -\eta_i)  + 2\omega_{jk}(\eta_j -\eta_k) + 4(\om_{ij}\om_{jk}
+ \om_{ki}\om_{ij} + \om_{jk}\om_{ki})
\Eeq
where $\eta_i$ is the normalized to 1 homogeneous volume form on the $i$-th copy of $S^1$
in the identification (\ref{6: fcM via FM}), and
$$
\om_{ij}:=\frac{1}{2\pi}dArg(z_i - z_j),
$$
a generator of $H^\bu(\mathsf{FM}_2(n)$. Note that according to the (non-cyclic)
Arnold relations the combination
$$
\om_{ij}\om_{jk}
+ \om_{ki}\om_{ij} + \om_{jk}\om_{ki}=d\beta_{ijk}
$$
is a $d$-exact 2-form on $\mathsf{FM}_2(n)$ and hence can be ignored upon integrations
over {\it boundaries}\, of closed semialgebraic varieties.

\subsection{Restrictions on the boundary strata}
Cyclic compositions in the topological operad $\FFM_2$,
$$
_i\circ_j: \FFM_2((n+1)) \times \FFM_2((m+1)) \lon \FFM_2((m+n))
$$
induce, for any decomposition $((m+n))=I'\sqcup I''$, the pull-back maps
\Beq\label{6: Delta on de Rham forms}
\Delta_{I',I''}: \Omega_{\fcM_{m+n}} \lon
\Omega_{\fcM_{I'\sqcup x'} \times \fcM_{I''\sqcup x''}}
\Eeq
These maps when restricted to  the propagators $\al_{ij}$ (and hence to $\Omega_{ijk}$) are given by the formulae (\ref{5: Delta on al})
(resp., by (\ref{5: Delta on Omega})) as they factor through the injections
$$
\Ba{rccc}
i: & \Omega_{\FFM_2((I'\sqcup x'))} \ot_\R \Omega_{\FFM_2((I''\sqcup x''))}  & \lon &
\Omega_{\FFM_2((I'\sqcup x')) \times \FFM_2((I''\sqcup x''))}\\
&  \om_1 \ot \om_2 &\lon & \pi_1^*(\om_1) \wedge \pi_2^*(\om_2),
\Ea\ 
$$
where $\pi_1$ and and $\pi_2$ are the canonical projections to the factors of the cartesian product
$\FFM_2((I'\sqcup x')) \times \FFM_2((I''\sqcup x''))$. It is proven in Lemma 7.1 in \cite{GS} that
\Beq\label{6: Delta on al}
\Delta_{I',I''}({\al}_{ij})=\left\{
\Ba{ll}
\pi_1^*({\al}_{ij}) & \text{if}\ i,j\in I'\\
\pi_2^*({\al}_{ij}) & \text{if}\ i,j\in I''\\
 \pi_1^*({\al}_{ix'})+ \pi_2^*({\al}_{x''j})  & \text{if}\ i\in I', j\in I''\\
\Ea
\right.
\Eeq
Hence one has for the degree 2 propagators  (cf.\ (\ref{5: Delta on Omega})),
\Beq\label{6: Delta on Omega}
\Delta_{I',I''}({\Omega}_{ijk})=\left\{
\Ba{ll}
\pi_1^*({\Omega}_{ijk}) & \text{if}\ i,j,k\in I'\\
\pi_2^*({\Omega}_{ijk}) & \text{if}\ i,j,k\in I''\\
\pi_1^*({\Omega}_{ijx'}) - (\pi_1^*({\al}_{ix'}) - \pi_1^*({\al}_{jx'}))\wedge \pi_2^*({\al}_{x''k})  & \text{if}\ i,j\in I', k\in I''\\
\pi_2^*({\Omega}_{x''jk}) + \pi_1^*({\al}_{ix'})\wedge
\left(\pi_2^*({\al}_{x''j}) - \pi_2^*({\al}_{x''k})\right)  & \text{if}\ i\in I', j,k\in I''\\
\Ea
\right.
\Eeq


\subsection{A morphism from $(\BVGcyc)^*$ to the de Rham algebra on $\FFM_2$}\label{6: subsec on morphism
CoBVGraphs to fM}

Given any $\Ga\in f\BVGcyc((n+1))^*$.
For any solid edge $e\in E(\Ga)$ we denote by $e',e''$ the two vertices it connects; similarly for a star
vertex $*\in V_*(\Ga)$ we denote by $v_*', v_*'',v_*'''$ the three vertices it connects.
Let us consider a linear map
\Beq\label{6: the map Phi}
\Ba{rccl}
\Phi: & (f\BVGcyc((n+1)))^*  & \lon &  \Omega_{\FFM_2((n+\# V_\bu(\Ga)+1))}\vspace{2mm} \\
& \Ga & \lon & \Phi_{n+1}(\Ga):=\displaystyle \bigwedge_{e\in E(\Ga)}\hspace{-2mm}
\al_{e'e''} \bigwedge_{*\in V_*(\Ga')}
\Omega_{v_*'v_*''v_*'''}
\Ea
\Eeq
The wedge products of 1-forms and the ordering of indices $v_*',v_*'',v_*'''$ in the
2-form  $\Omega_{v_*'v_*''v_*'''}$ are taken in accordance with the given orientation
of $\Ga$. There is  a semialgebraic fibration \cite{GS}
\Beq\label{6: fibration pi which forgets black points}
\pi^\Ga: \FFM_2((n+\# V_\bu(\Ga) +1)) \lon \FFM_2((n+1))
\Eeq
which forgets all black vertices.  It induces a well-defined push-down
morphism of graded vector spaces \cite{KS, HLTV, GS}
$$
\Ba{rccc}
\Omega_{n+1}: & f\BVGcyc((n+1))^*& \lon & \Omega_{\FFM_2((n+1))}\vspace{3mm}\\
& \Ga & \lon & \Omega_\Ga:= \pi_*^\Ga(\Phi(\Ga))
\Ea
$$
Thanks to formulae (\ref{6: Delta on al})-(\ref{6: Delta on Omega})
for the propagators, these maps are compatible with
of cooperadic structures on $(\BVGcyc)^*$ and $\FFM_2$ in the sense that  the
following diagram,
\Beq\label{6: diagram compatibility}
 \xymatrix{
  \ar[d]_{\Delta_{I',I''}}f\BVGcyc((n+m))^*\hspace{-8mm}
  \ \ \ \ \ \ \ \ \ \ \  \stackrel{\Omega_{n+m}}{\lon} \ \ \ \   
  \Omega_{\FFM_2((m+n))}\ar[r]^{\ \ \ \ \ \Delta_{I',I''} } &  \Omega_{\FFM_2((I'\sqcup x'))\times
  \FFM_2((I''\sqcup x''))}     \\
 f\BVGcyc((I'\sqcup x'))^*\ot f\BVGcyc((I''\sqcup x''))^* \ 
 \ar[r]_{\hspace{20mm}\Omega_{n +1}\ot \Omega_{m +1}} 
 & \ \ \ \ \  \Omega_{\FFM_2((I'\sqcup x'))}  \ot \Omega_{\FFM_2((I''\sqcup x''))}\ar[u]_i
}
\Eeq
commutes 
for any partition $((m+n))=I'\sqcup I''$ into non-empty subsets. 

\subsection{\bf Theorem} {\em The  morphism of graded $\bS$-modules,
\Beq\label{6: map Omega from coBVGraphs to Omega_fM}
\Ba{rccc}
\Omega: & f\BVGcyc^* & \lon &  \Omega_{\FFM_2}\vspace{3mm}\\
        & \Ga & \lon & \Omega(\Ga):=\pi_*^\Ga(\Phi(\Ga))
\Ea
\Eeq
 respects the differentials.}

\begin{proof} We have to show that for any hypergraph $\Ga\in (f\BVGcyc)^*$ one has
$$
\Omega(d\Ga)=d_{DR}\Omega(\Ga).
$$
Since the differential form $\Phi(\Ga)$ is closed for any $\Ga$, the fiberwise Stokes
 formula (see Proposition 9 in \cite{HLTV} and Theorem 6.2 in \cite{GS}) implies
$$
d_{DR}\Omega(\Ga)=\pi_*^\Ga(\Phi(\Ga))=\p \pi_*^\Ga(\Phi(\Ga)),
$$
where $\p \pi_*^\Ga(\Phi(\Ga))$ stands for the push-down along the boundary of the fibration
(\ref{6: fibration pi which forgets black points}). That boundary consists of limit
configurations of two types: ($B_1$) $p\geq 2$ ``black"  framed points (i.e.\ the one
which correspond to black vertices in $\Ga$) collapse into a ``black" point;
($B_2$) $q\geq 1$ ``black" framed points collapse into some ``white" point in the base
 $\FFM_2((n+1))$ (corresponding to some white vertex $\wi\in \Ga$). Since $\Phi(\Ga)$
 is eventually a wedge product of 1-forms $\al_{ij}$, one applies the Vanishing
 Lemma 6.5 from \cite{GS} (which in turn follows from  the Kontsevich Vanishing Lemma 6.13
 from \cite{Ko2}) to conclude that the only boundary strata which contributes to the
 push-down $\p \pi_*^\Ga(\Phi(\Ga))$ has $p=2$ in the case of the type $B_1$ boundary stratum and $q=1$ in
 the case of the boundary stratum of type $B_2$. In both cases one integrates in
$\p \pi_*^\Ga(\Phi(\Ga))$ over a 2-dimensional subspace in $\FFM_2((3))=S^1\times S^1\times S^1$,
where one factor corresponds to the tangent ray at the newly created node in the limit
configuration, and the other two factors correspond to the tangent rays of the collapsing
framed points; call them  $u$ and $v$ as in \cite{GS}. This boundary stratum can give a
non-zero contribution into the push-down $\p \pi_*(\Phi(\Ga))$ if and only if the
collapsing points $u$ and $v$ are connected either by an edge $e$ or by a hyperedge
in $\Ga$; moreover, it is easy to see that a non-zero contribution into the
boundary push-down can only come from the differential 2-form of the form
$\omega_{uv}\wedge \eta_v$ if it is present as a factor in $\Phi(\Ga)$: the angular
propagator $\om_{uv}$ takes care about integrating over the boundary factor $S^1$
corresponding to the node in the limit configuration, and $\eta_v$ contributes via
its integration over the boundary factor $S^1$ corresponding to real rays at $v$.
Let us consider next all possible limiting configurations of a pair of framed
points $u$ and $v$ (=vertices of $\Ga$)  which might contribute into $\p \pi_*^\Ga(\Phi(\Ga))$.

\sip

{\sc Option 1}: both collapsing points $u$ and $v$ correspond to black vertices in $\Ga$.
In this case $u$ and $v$ must be connected by a hyperedge $\Omega_{uvk}$ here the
flag $k$ is connected to some white white vertex $\wk$, and the affine representation
(\ref{6: Omega_ijk in affine coordiantes}) of the propagator $\Omega_{uvk}$ assures
us that such limit configurations do give us a non-zero contribution via the
summand $\omega_{uv}\wedge(\eta_v - \eta_u)$ in $\Omega_{uvk}$. All such
contributions into $\p \pi_*^\Ga(\Phi(\Ga))$ can hence be identified with
$\Omega(d_{\bu\bu}\Ga)$, i.e.\ they all correspond to part $d_{\bu\bu}$ of the full
differential (\ref{6: d on CoBVGraphs as three parts}) in $(f\BVGcyc)^*$.

\sip

{\sc Option 2}: one collapsing point, say $v$,  corresponds to a black vertex in $\Ga$,
another collapsing point $u$ corresponds to a white vertex $\wi$ in $\Ga$, and they are
 connected by a hyperedge. If these two vertices are also connected by an edge, the
 contribution is zero, otherwise the 2-form (\ref{6: Omega_ijk in affine coordiantes})
contains the required summand $\omega_{uv}\wedge \eta_v$ and hence integrates to 1 on the boundary
strata of this types giving rise to the contributions into $\p \pi_*^\Ga(\Phi(\Ga))$ which are
 equal precisely to $\Omega(d_{\bu\circ}\Ga)$.

\sip

{\sc Option 3}: one collapsing point, say $v$,  corresponds to a black vertex in $\Ga$,
another collapsing point, say $u$,  corresponds to a white vertex $\wi$ in $\Ga$, and
they are connected by an edge. This edge contributes the factor $\omega_{uv}$ into the
integral over the boundary strata. The other required factor $\eta_v$ can only originate
 either from flags or from other edges attached to the black vertex $v$. Both options
  can contribute, and we consider them separately.

\sip

{\em Option 3a}: $\eta_v$ originates from the propagator
 $\al_{vj}=\eta_v + \eta_j - 2\omega_{uj}$ corresponding to an edge which connects
 $v$ to some other white vertex $\wj$\ . The associated contributions into
the push-down  $\p \pi_*^\Ga(\Phi(\Ga))$ are equal precisely to
 $\Omega(d_{\bu\circ}^{(2)}\Ga)$.

\sip

{\em Option 3b}: $\eta_v$ originates from the propagator
$\Omega_{vxy}$ which connects $v$ to some other vertices $x,y$ in the limit $v\rar \wi$.
This is the most intriguing case. We can re-write
(\ref{6: Omega_ijk in affine coordiantes}) keeping only terms which have $\eta_v$ as a factor
$$
\Omega_{vxy}=\eta_v(\eta_x - \eta_y  - 2\omega_{vx} + 2\omega_{vy})+\ldots.
$$
After integrating out the standard form $\omega_{uv}\wedge \eta_v$   we are left in the
limit $v\rar \wi$ with the following linear combination of 1-forms,
$$
\al_{ax}-\al_{ay}
$$
which is precisely the result of the operations
(\ref{6: d'' in cooperad, step 1})-(\ref{6: d'' in cooperad, step 2})! Hence this class
of contributions into $\p \pi_*^\Ga(\Phi(\Ga))$ corresponds  precisely to
 $\Omega(d_{\bu\circ}^{(1)}\Ga)$. The proof is completed.
\end{proof}

\mip

The commutative diagram (\ref{6: diagram compatibility}) and the above Theorem imply the following result for the operad  $Chains(\FFM_2)$ of the semialgebraic 
chains  on $\FFM_2$.

\subsection{Proposition}\label{6: Prop on F from chains to fBVGcec} {\it There is a morphism of dg cyclic operads
\Beq\label{6: F chains to fBVHgraphs}
\cF: Chains(\FFM_2)\lon f\BVGcyc.
\Eeq
which sends a chain $g: \vartriangle \rar \FFM_2$ into the sum
$\sum_{\Ga} (\int_\vartriangle g^*(\Omega_\Ga))\Ga$ over the generators $\Ga$ of $f\BVGcyc$.}

\mip

As we show in the next section this morphism $\cF$ can not be a quasi-isomorphism unless we replace $f\BVGcyc$ by a ``smaller" cyclic operad.

\bip

{\large
\section{\bf A dg cyclic operad $\BVGcyc$ of forests}
}

\sip

\subsection{Hopf cooperad structure in $f\BVHgraphs$}

 For any 
 generators $\Ga_1,\Ga_2\in f\BVHgraphs^*(n)$, one can defined their product $\Ga_1\cdot \Ga_2\in f\BVHgraphs^*(n)$ as their disjoint union with subsequence identification of the corresponding white vertices. For example, the hypergraph
 $$
  \Ba{c}\resizebox{12mm}{!}
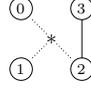
{ \xy
(-5,0)*+{_1}*\frm{o}="1";
 (5,0)*+{_2}*\frm{o}="2";
 (-5,10)*+{_0}*\frm{o}="0";
  (5,10)*+{_3}*\frm{o}="3";
  (0,5)*{*}="s";
 \ar @{.} "s";"0" <0pt>
 \ar @{.} "s";"1" <0pt>
  \ar @{.} "s";"2" <0pt>
  \ar @{-} "3";"2" <0pt>
\endxy} \Ea
 $$
 is the product   $\Ga_1\cdot \Ga_2$     of
 $$
 \Ga_1=
  \Ba{c}\resizebox{12mm}{!}
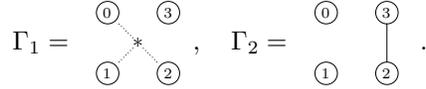
{ \xy
(-5,0)*+{_1}*\frm{o}="1";
 (5,0)*+{_2}*\frm{o}="2";
 (-5,10)*+{_0}*\frm{o}="0";
  (5,10)*+{_3}*\frm{o}="3";
  (0,5)*{*}="s";
 \ar @{.} "s";"0" <0pt>
 \ar @{.} "s";"1" <0pt>
  \ar @{.} "s";"2" <0pt>
\endxy} \Ea, \ \ \  
\Ga_2= \Ba{c}\resizebox{12mm}{!}{ \xy
(-5,0)*+{_1}*\frm{o}="1";
 (5,0)*+{_2}*\frm{o}="2";
 (-5,10)*+{_0}*\frm{o}="0";
  (5,10)*+{_3}*\frm{o}="3";
  \ar @{-} "3";"2" <0pt>
\endxy} \Ea\ .
 $$
 The differential and the co-operadic compositions respect this product so that 
 $f\BVHgraphs^*$ is a dg cyclic operad in the category of dg commutative associative algebras (cf.\ \cite{Ko,SW}).

\sip

A  hypergraph $\Ga\in f\BVGcyc^*$ (or in  $f\BVGcyc$) is called {\em internally connected}\,
if it becomes (or stays) connected after erasing all white vertices. For example, the first three hypergraphs below
are internally connected,
$$
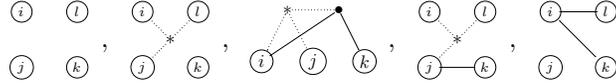
\Ba{c}\resizebox{11mm}{!}{ \xy
(-5,0)*+{_j}*\frm{o}="1";
 (5,0)*+{_k}*\frm{o}="2";
 (-5,10)*+{_i}*\frm{o}="0";
  (5,10)*+{_l}*\frm{o}="3";
\endxy} \Ea
,
 \Ba{c}\resizebox{11mm}{!}{ \xy
(-5,0)*+{_j}*\frm{o}="1";
 (5,0)*+{_k}*\frm{o}="2";
 (-5,10)*+{_i}*\frm{o}="0";
  (5,10)*+{_l}*\frm{o}="3";
  (0,5)*{*}="s";
 \ar @{.} "s";"0" <0pt>
 \ar @{.} "s";"1" <0pt>
  \ar @{.} "s";"3" <0pt>
\endxy} \Ea
,
\Ba{c}\resizebox{17mm}{!}{ \xy
(-10,0)*+{i}*\frm{o}="0";
(0,0)*+{j}*\frm{o}="1";
(10,0)*+{k}*\frm{o}="2";
  (5,10)*{\bu}="bu";
 (-5,10)*{*}="s";
 \ar @{-} "bu";"2" <0pt>
 \ar @{-} "bu";"0" <0pt>
 \ar @{.} "s";"bu" <0pt>
  \ar @{.} "s";"0" <0pt>
 \ar @{.} "s";"1" <0pt>
\endxy} \Ea
,
 \Ba{c}\resizebox{11mm}{!}{ \xy
(-5,0)*+{_j}*\frm{o}="1";
 (5,0)*+{_k}*\frm{o}="2";
 (-5,10)*+{_i}*\frm{o}="0";
  (5,10)*+{_l}*\frm{o}="3";
  (0,5)*{*}="s";
 \ar @{.} "s";"0" <0pt>
 \ar @{.} "s";"1" <0pt>
  \ar @{.} "s";"3" <0pt>
  \ar @{-} "1";"2" <0pt>
\endxy} \Ea
,
\Ba{c}\resizebox{11mm}{!}{ \xy
(-5,0)*+{_j}*\frm{o}="1";
 (5,0)*+{_k}*\frm{o}="2";
 (-5,10)*+{_i}*\frm{o}="0";
  (5,10)*+{_l}*\frm{o}="3";
 \ar @{-} "0";"2" <0pt>
 \ar @{-} "0";"3" <0pt>
\endxy} \Ea
$$
and the last two are not.

\sip

\subsection{\bf Classification of hyperedges} 
Every hyperedge (if any) of a generator $\Ga\in f\BVGcyc$  
has one of the following tree types:
\Bi
\item[{\it Type-0}]  {\it hyperedge} is, by definition, a hyperedge which is connected to three white vertices, 
$$
  \text{type-0 hyperedge}:\ \ \ \ 
  \Ba{c}\resizebox{30mm}{!}{ \xy
(-23,0)*+{_0}*\frm{o}="0";
(-7.5,0)*{...};
(+7.5,0)*{...};
(-17.5,0)*{...};
(-12,0)*+{_i}*\frm{o}="0";
(0,0)*+{_j}*\frm{o}="1";
(12,0)*+{_n}*\frm{o}="2";
  (0,9)*{*}="s";
 \ar @{.} "s";"2" <0pt>
 \ar @{.} "s";"0" <0pt>
 \ar @{.} "s";"1" <0pt>
\endxy} \Ea, \ \ \ 0\leq i<j<n
$$

\item[{\it Type-1}] {\it hyperedge} is, by definition, a hyperedge which is connected to two black vertices and to a unique white vertex $\ \wi\ $ with $i\in \{0,1,\dots, n\}$
    (small integers show a choice of an ordering of its flags)
    \Beq\label{7: type1 hyperedge}
    \text{type-1 hyperedge}:\ \ \ \ 
     \Ba{c}\resizebox{12mm}{!}
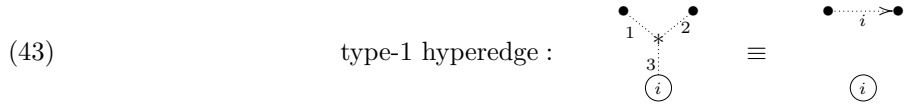
{ \xy
(0,-3)*+{_i}*\frm{o}="0";
(-4,5)*{_{1}};
(-1,0.4)*{_{3}};
(4,5)*{^{2}};
 (0,4)*{*}="s";
  (5,8)*{\bu}="r";
  (-5,8)*{\bu}="l";
 \ar @{.} "0";"s" <0pt>
 \ar @{.} "s";"l" <0pt>
  \ar @{.} "s";"r" <0pt>
\endxy}\Ea\ \ \ \equiv \ \ \   \Ba{c}\resizebox{12mm}{!}{ \xy
(0,-3)*+{_i}*\frm{o}="0";
(0,6)*{^{i}};
  (5,8)*{\bu}="r";
  (-5,8)*{\bu}="l";
 \ar @{.>} "l";"r" <0pt>
\endxy}\Ea
    \Eeq
Note that the parts $\delta_i^{(2)}$ and $\delta_i^*$ of the differential  act trivially on the flag connected to $\ \wi\ $  as we take the quotient by hyperedges of the form $
 \Ba{c}\resizebox{12mm}{!}{ \xy
(-6,0)*{\bu}="1";
 (6,0)*{\bu}="2";
 (0,10)*{\bu}="0";
  (0,4)*{*}="s";
 \ar @{.} "s";"2" <0pt>
 \ar @{.} "s";"0" <0pt>
 \ar @{.} "s";"1" <0pt>
\endxy} \Ea$. 
Hence we can identify such a hyperedge  with a {\it directed labelled by $i$  dotted  edge}\, as in the picture above. The direction of such a dotted edge can be flipped in accordance with the following sign rule
$$
\Ba{c}\resizebox{12mm}{!}{ \xy
(0,-2)*{^{i}};
  (5,0)*{\bu}="r";
  (-5,0)*{\bu}="l";
 \ar @{.>} "l";"r" <0pt>
\endxy}\Ea = - \Ba{c}\resizebox{12mm}{!}{ \xy
(0,-2)*{^{i}};
  (5,0)*{\bu}="r";
  (-5,0)*{\bu}="l";
 \ar @{.>} "r";"l" <0pt>
\endxy}\Ea
$$
so in pictures we do not show directions at all (assuming tacitly that some choice has been made). It is worth noting that in this notation the part $\delta_\bu\Ga=\sum_{v\in V_\bu(\Ga)} \delta_v \Ga$ of the full differential $\delta$ acts on black vertices as follows
$$
\delta_v: v=\bu\ \ \ \rar \ \  \sum_{i=0}^n\Ba{c}\resizebox{12mm}{!}{ \xy
(0,-2)*{^{i}};
  (5,0)*{\bu}="r";
  (-5,0)*{\bu}="l";
 \ar @{.>} "l";"r" <0pt>
\endxy}\Ea.
$$


\item[{\it Type-2}] {\it hyperedge}\, is, by definition, a hyperedge which is connected to two white vertices and to a unique black vertex,
$$
  \text{type-2 hyperedge}:\ \ \ \ \Ba{c}\resizebox{12mm}{!}
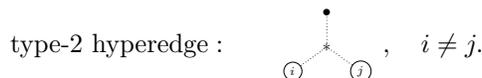
{ \xy
(-7,0)*+{_i}*\frm{o}="0";
(7,0)*+{_j}*\frm{o}="1";
(0,12)*{\bu}="2";
  (0,5)*{*}="s";
 \ar @{.} "s";"2" <0pt>
 \ar @{.} "s";"0" <0pt>
 \ar @{.} "s";"1" <0pt>
\endxy} \Ea  , \ \ \ i\neq j.
$$
The flags of a such hyperedge which are attached to white vertices are called  {\it type-2 flags}.

\Ei

\subsection{Proposition}\label{7: Prop on degrees of fBVHgraphs}  {\it The dg cyclic operad $f\BVHgraphs$ is non-positively graded, that is, for any $\Ga\in f\BVHgraphs$ one has $|\Ga|\leq 0$.}

\begin{proof} The internally connected graphs in  $f\BVHgraphs$ which have no hyperedges must have at most one single black vertex with $m\geq 3$ edges attached. Each such graph has the cohomological degree $3-m\leq 0$ so the Proposition holds true for such graphs.

\sip

Consider next a connected hypergraph  $\Ga\in f\BVHgraphs$ which has at least one hyperedge.
In general $\Ga$ can have (i) $n_k\geq 1$ black vertices which have precisely  $k$ flags attached, (ii) $p$ hyperedges of type-1, (iii) $q$ hyperedges of type-2 and at least $2n_1 + n_2$ edges as each black vertex is at least trivalent. We have
$$
2p+q=\sum_{k\geq 1} k n_k, \ \ \ \ \# E(\Ga)\geq 2n_1 + n_2.  
$$
and
\Beqrn
|\Ga| &=& 3 (\sum_{k\geq 1} n_k) - 2(p+q) - \# E(\Ga)\\
      &\leq & (\sum_{k\geq 1}(3-k) n_k) - q - 2n_1 - n_2 \\                                     &\leq & \sum_{k\geq 4}(3-k) n_k -q \leq  0.
\Eeqrn
The Proposition is proven.
\end{proof}

\subsection{A dg sub-cooperad of forests}
We define a dg cyclic sub-cooperad  $\wBVHgraphs^* \subset f\BVHgraphs^*$ as the one generated by hypergraphs $\Ga$ whose internally connected components are {\it trees}\, (i.e.\ they do not contain closed paths consisting solely of flags of hyperedges only). By the Proposition above,
 $\wBVHgraphs^*$ is concentrated in non-negative degrees. The same is true
for the the semialgebraic de Rham complex  $\Omega^\bu_{\FFM_2}$. 
The morphism (\ref{6: map Omega from coBVGraphs to Omega_fM})
restricts to the following one
\Beq\label{7: F wBVHgraphs to Omega_fM}
\Ba{rccc}
\Omega: & \wBVHgraphs^* & \lon &  \Omega_{\FFM_2}\\
        & \Ga & \lon & \Omega(\Ga):=\pi_*^\Ga(\Phi(\Ga)).
\Ea
\Eeq
This map is a morphism of dg graded commutative algebras,
$$
\Omega(\Ga_1\cdot \Ga_2)= \Omega(\Ga_1) \wedge \Omega(\Ga_2).
$$
However it is not a quasi-isomorphism because a certain class of hypergraphs (including the ones representing non-trivial cohomology classes) goes under the map $\Omega$ to zero.

\subsection{Theorem}\label{7: Theorem on Omega vanishing on degree zero} {\it If $\Ga\in \wBVHgraphs^*$ is an internally connected hypergraph which has cohomological degree zero and contains at least one black vertex, then $\Omega(\Ga)=0$.
}

\begin{proof}
In accordance with the morphism (\ref{6: the map Phi})  the differential form $\Phi_\Ga$ associated to a hypergraph $\Ga$ is obtained from the latter by assigning to each edge the propagator 
$$
\om_{ij}=\frac{1}{2\pi}dArg(z_i - z_j),
$$ 
and to each hyperedge the propagator 
$$
\Omega_{ijk}= \om_{ij}\om_{jk}
+ \om_{ki}\om_{ij} + \om_{jk}\om_{ki}
$$
In combinatorial terms the latter association  corresponds to the replacement
of each hyperedge in $\Ga$ by a linear combination of edges as in the following picture (we skip showing the ordering of edges)
\Beq\label{7: hyperedge goes to pairs of edges} 
\Ba{c}\resizebox{16mm}{!}{ \xy
(-7,0)*+{_i}="0";
(7,0)*+{_j}="1";
(0,13.5)*+{_k}="2";
  (0,5)*{*}="s";
 \ar @{.} "s";"2" <0pt>
 \ar @{.} "s";"0" <0pt>
 \ar @{.} "s";"1" <0pt>
\endxy} \Ea  \lon 
\Ba{c}\resizebox{16mm}{!}{ \xy
(-7,0)*+{_i}="i";
(7,0)*+{_j}="j";
(0,13.5)*+{_k}="k";
 \ar @{-} "i";"j" <0pt>
 \ar @{-} "j";"k" <0pt>
\endxy}\Ea
+
\Ba{c}\resizebox{16mm}{!}{ \xy
(-7,0)*+{_i}="i";
(7,0)*+{_j}="j";
(0,13.5)*+{_k}="k";
 \ar @{-} "i";"k" <0pt>
 \ar @{-} "i";"j" <0pt>
\endxy}\Ea
+
\Ba{c}\resizebox{16mm}{!}{ \xy
(-7,0)*+{_i}="i";
(7,0)*+{_j}="j";
(0,13.5)*+{_k}="k";
 \ar @{-} "i";"k" <0pt>
 \ar @{-} "k";"j" <0pt>
\endxy}\Ea
\Eeq

If $\Ga$ is any graph having a trivalent vertex of the form 
 $$
 \Ba{c}\resizebox{25mm}{!}
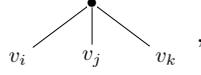
{ \xy
(-12,0)*+{v_i}="0";
(0,0)*+{v_j}="1";
(12,0)*+{v_k}="2";
  (0,9)*{\bu}="s";
 \ar @{-} "s";"2" <0pt>
 \ar @{-} "s";"0" <0pt>
 \ar @{-} "s";"1" <0pt>
\endxy} \Ea,
$$
then $\Omega(\Ga)=0$ by Lemma 6.3(ii) in \cite{GS}. The authors of \cite{GS} prove this claim by (a) fixing the three vertices $v_i,v_j,v_k$ 
 to which the trivalent black vertex is attached by solid edges at the positions $0$, $1$ and $\infty$ in the Riemann sphere, and then (b) by noticing that the complex conjugation is, on the one hand, an orientation preserving diffeomorphism of the tangent circle bundle, but on the other hand it reverses the sign of each propagator
$$
dArg(\overline{z_i - z_j})= - dArg(z_i - z_j).
$$
Hence the result of the integration $\pi_*^\Ga(\Phi(\Ga))$ for any graph containing a trivalent black vertex as in the picture above must be zero.

\sip

Any internally connected  tree $\Ga\in \wBVHgraphs^*$ with $|\Ga|=0$ must have black vertices precisely trivalent, and must have hyperedges (if any) of type-1 only; here is a typical example,
$$
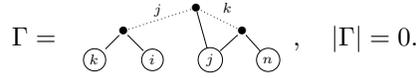
\Ga= \Ba{c}\resizebox{27mm}{!}{ \xy
(-10,1)*+{_k}*\frm{o}="0";
(0,1)*+{_i}*\frm{o}="1";
(10,1)*+{_j}*\frm{o}="2";
(20,1)*+{_n}*\frm{o}="n";
  (-5,6)*{\bu}="bu";
   (7.5,10)*{\bu}="r";
     (15.5,6)*{\bu}="rr";
(1,9.5)*{_j};
(13,10)*{_k};
 \ar @{-} "0";"bu" <0pt>
  \ar @{-} "1";"bu" <0pt>
   \ar @{.} "bu";"r" <0pt>
    \ar @{-} "rr";"2" <0pt>
     \ar @{-} "rr";"n" <0pt>
  \ar @{.} "rr";"r" <0pt>
 \ar @{-} "r";"2" <0pt>
\endxy} \Ea , \ \ \ |\Ga|=0.
$$
Under the map $\Omega$  one substitutes hyperedges in $\Ga$ as 
in the formula (\ref{7: hyperedge goes to pairs of edges}); it is easy to see that this substitution always creates at least one black vertex with precisely three solid edges. Hence $\Omega(\Ga)$ vanishes  by the above mentioned Lemma 6.3(ii) in \cite{GS}.
\end{proof}

\subsection{Unwanted cohomology classes in $\wBVHgraphs$}   Let  $\wBVHgraphs$ be the dg cyclic operad which is dual to dg cooperad $\wBVHgraphs^*$ defined above.
The operad  $\wBVHgraphs$ is the quotient of $f\BVHgraphs$ by the dg ideal generated by hypergraphs with
internal genus.
 By Proposition~{ \ref{7: Prop on degrees of fBVHgraphs}},   the dg cyclic
operad $\wBVHgraphs=\{\wBVHgraphs((n+1))\}_{n\geq 1}$ is concentrated in cohomological degrees $\leq 0$, and the morphism (\ref{5: Main Morphism F}) induces the following morphism of dg cyclic operads
\Beq\label{7: Main Morphism F from BVc to wHgraphs}
F: \BVc \lon \wBVHgraphs
\Eeq
given on the generators by the formulae (\ref{5: f from BV to BVgrac}).

\sip

Any graph of the form  $$
 \Ba{c}\resizebox{30mm}{!}
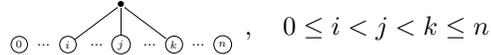
{ \xy
(-23,0)*+{_0}*\frm{o}="0";
(23,0)*+{_n}*\frm{o}="0";
(-5.5,0)*{...};
(+6.5,0)*{...};
(-17.5,0)*{...};
(17.5,0)*{...};
(-12,0)*+{_i}*\frm{o}="0";
(0,0)*+{_j}*\frm{o}="1";
(12,0)*+{_k}*\frm{o}="2";
  (0,9)*{\bu}="s";
 \ar @{-} "s";"2" <0pt>
 \ar @{-} "s";"0" <0pt>
 \ar @{-} "s";"1" <0pt>
\endxy} \Ea, \ \ \ 0\leq i<j<k \leq n
$$
is a {\it degree zero cohomology class}\, in the operad $\wBVHgraphs$. Similarly, some
  trivalent  hypergraphs  of the form (defined modulo IXH relations at each dotted edge)
$$
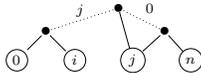
\Ba{c}\resizebox{27mm}{!}{ \xy
(-10,1)*+{_0}*\frm{o}="0";
(0,1)*+{_i}*\frm{o}="1";
(10,1)*+{_j}*\frm{o}="2";
(20,1)*+{_n}*\frm{o}="n";
  (-5,6)*{\bu}="bu";
   (7.5,10)*{\bu}="r";
     (15.5,6)*{\bu}="rr";
(1,9.5)*{_j};
(13,10)*{_0};
 \ar @{-} "0";"bu" <0pt>
  \ar @{-} "1";"bu" <0pt>
   \ar @{.} "bu";"r" <0pt>
    \ar @{-} "rr";"2" <0pt>
     \ar @{-} "rr";"n" <0pt>
  \ar @{.} "rr";"r" <0pt>
 \ar @{-} "r";"2" <0pt>
\endxy} \Ea 
$$
give us degree zero cohomology classes in $H^0(\wBVHgraphs)$. All such hypergraphs   vanish under  the morphism  $\Omega$. Hence to get the things right we have to truncate all the degree zero cohomology classes from $\wBVHgraphs$ which we do next.

\subsection{\bf A cyclic operad of $\Lie_\infty$-algebras $\wICH$} The dg cyclic operad $\wBVHgraphs$ is a Hopf operad, that is a cyclic operad in the category of coalgebras.
Let $\wICH((n+1))[1]$, $n\geq 1$, be the subspace of the complex $\wBVHgraphs((n+1))$ spanned by internally connected hypergraphs with at least one edge or hyperedge.
We decompose
\Beq\label{7: from wBVHgraphs to wICH}
\wBVHgraphs(n+1)=\left(\odot^\bu (\wICH((n+1))[1]), \delta=\sum_{k\geq 1}\delta_k\right)
\Eeq
where the constant summand $1$ in the symmetric tensor algebra corresponds to the unique graph
with $n+1$ white vertices and no edges and hyperedges, 
$$
  \Ba{c}\resizebox{30mm}{!}
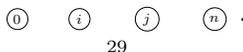
{ \xy
(-23,0)*+{_0}*\frm{o}="0";
(-12,0)*+{_i}*\frm{o}="0";
(0,0)*+{_j}*\frm{o}="1";
(12,0)*+{_n}*\frm{o}="2";
\endxy} \Ea.
$$
The summands 
\Beq\label{7: delta_k for wICH}
\delta_k: \odot^k (\wICH((n+1))[1]) \lon \wICH((n+1))[1], 
 \ \ \ \forall k\geq 1,
\Eeq
of the full differential represent the induced $\Lie_\infty$ compositions on  $\wICH((n+1))$. The collection
$$
\wICH=\{\wICH((n+1))\}_{n\geq 1}
$$
is a cyclic operad in the category of $\Lie_\infty$-algebras (cf.\ \cite{SW}). Due to the degree shift the graded spaces $\wICH((n+1))$ are all concentrated in degrees $\leq +1$, and the ``unwanted" cohomology classes all have  degree $+1$ in $\wICH((n+1))$ (cf.\ Theorem {\ref{7: Theorem on Omega vanishing on degree zero}}).

\subsection{Truncating $\wICH$} 
We define dg vector subspaces for any $n\geq 1$,
$$
\ICH((n+1))\equiv \bigoplus_{i=-\infty}^1\ICH^i((n+1)) 
\ \ \subset \ \  \wICH((n+1))\equiv \bigoplus_{i=-\infty}^1\wICH^i((n+1)) 
$$ 
by truncating the above trivalent cohomology classes as follows,
\Beq\label{7: def of ICH}
\ICH^i((n+1)):=\wICH^i((n+1)) \ \  \forall i\leq -1,  \ \ \ \ICH^0((n+1)):=\Ker \delta_1\subset \wICH  ,\ \ \    \ICH^1((n+1)):=0
\Eeq
Thus, by the very definition, we have
\Beq\label{8 H(ICH) versus H(wICH)}
H^i(\ICH((n+1)))=H^i((\wICH(n+1))) \  \forall i\leq 0,  \ \ H^1(\ICH((n+1)))=0.
\Eeq

The collection
$$
\ICH=\{\ICH((n+1))\}_{n\geq 1}
$$
is a dg cyclic suboperad of $\wICH$ in the category of $\Lie_\infty$-algebras.

\subsection{Main Definition} Applying the Chevalley symmetric monoidal functor to the operad  $\ICH$ of $\Lie_\infty$ algebras, we obtain
a dg cyclic suboperad  
$$
\BVHgraphs:= \odot^\bu(\ICH[1]) \subsetneq \odot^\bu(\wICH[1])=\wBVHgraphs
$$

The morphism of dg cyclic operads (\ref{7: Main Morphism F from BVc to wHgraphs}) factors through the inclusion $\BVHgraphs \subset \wBVHgraphs$.

\subsection{Main Theorem}\label{7: Main Theorem} {\it The morphism of dg cyclic operads
\Beq\label{F: BV to BVHgraphs}
\BVc\lon \BVHgraphs 
\Eeq
is a quasi-isomorphism}.

\sip

We prove the Main Theorem in several steps below in  \S\S\, 8-10  using the quasi-isomorphism (\ref{2: BVc to CE(T^c)}) and the fact that the cyclic operad $\ftt$ of infinitesimal ribbon braids contains free Lie algebras (as discussed in \S 2.8). Our approach is inspired by  Appendix B in   
$\check{\text{S}}$evera-Willwacher's paper \cite{SW} though the technical details in our situation are not, unfortunately, that short.

\sip

Let us formulate several immediate corollaries of the Main Theorem.

\subsection{Dual cooperad $\BVHgraphs^*$ and the formality of $\FFM_2$} There is a projection
$$
\pi: \wBVHgraphs^* \lon  \BVHgraphs^*
$$
of the dg cyclic cooperads. The Main Theorem together with Theorem {\ref{7: Theorem on Omega vanishing on degree zero}} imply the following

\subsubsection{\bf Corollary} {\it The map $\Omega: \wBVHgraphs^*\rar \Omega_{\FFM_2}$ in (\ref{6: map Omega from coBVGraphs to Omega_fM}) factors through the projection $\pi$ and induces a quasi-isomorphism of dg algebras
$$
\Omega: \BVHgraphs^* \lon \Omega_{\FFM_2},
$$
which fits the corresponding version of the commutative diagram (\ref{6: diagram compatibility}).}
\sip

This Corollary together with Proposition {\ref{6: Prop on F from chains to fBVGcec}} imply the following 

\subsection{Corollary} {\it There is a morphism of dg cyclic operads
\Beq\label{7: f chains to BVHgraphs}
\cF: Chains(\FFM_2)\lon \BVGcyc.
\Eeq
which is a quasi-isomorphism.}

\sip

The latter Corollary together with the Main Theorem {\ref{7: Main Theorem}}
imply the formality of $\FFM_2$ as a topological cyclic  operad.


\sip


\bip

{\large
\section{\bf Operad of internally connected hypergraphs}
}

\sip

\subsection{From infinitesimal ribbon braids to $\wICH$} Let
$$
\wICH=\{ \wICH((n+1))\}_{n\geq 1}
$$
be the cyclic operad of $\Lie_\infty$-algebras defined by the decomposition
(\ref{7: from wBVHgraphs to wICH}), and
$$
\ICH=\{ \wICH((n+1))\}_{n\geq 1}
$$
its cyclic suboperad defined by the truncation (\ref{7: def of ICH}).

\sip

Note that a hypergraph $\Ga\in\wICH((n+1)) $ has no solid edges between white vertices  unless $\Ga$ consists of white vertices only and has precisely one such edge; denote a one-edge graph by $\ga_{ij}$ if the edge
connects two different white vertices labelled by $i$ and $j$, e.g.
$$
\ga_{01}= \Ba{c}\resizebox{11mm}{!}{ \xy
(-5,0)*+{_1}*\frm{o}="1";
 (5,0)*+{_2}*\frm{o}="2";
 (-5,10)*+{_0}*\frm{o}="0";
  (5,10)*+{_3}*\frm{o}="3";
 \ar @{-} "0";"1" <0pt>
\endxy} \Ea
, \ \
 \ga_{12}= \Ba{c}\resizebox{11mm}{!}{ \xy
(-5,0)*+{_1}*\frm{o}="1";
 (5,0)*+{_2}*\frm{o}="2";
 (-5,10)*+{_0}*\frm{o}="0";
  (5,10)*+{_3}*\frm{o}="3";
 \ar @{-} "2";"1" <0pt>
 \endxy} \Ea
$$
 Note that all such graphs
$\ga_{ij}$ have cohomological degree zero in $\wICH(n+1)$ and represent cohomology classes in $H^0(\wICH)$,
$$
\delta \ga_{ij}=0, \ \ \forall \ i,j\in \{0,1,\ldots, n\}, \ i\neq j.
$$
In particular, all such graphs $\ga_{ij}$ belong to $\ICH\subset \wICH$.

\sip

Recall that $\ftt=\{\ftt((n+1))\}_{n\geq 1}$ stands for the cyclic operad of infinitesimal ribbon braids (see \S {\ref{2: subsec on T_cyc}}).

\subsubsection{\bf Proposition}\label{8: Prop on g from ftt to ICH}   {\em There is a morphism of cyclic operads
of $\Lie_\infty$-algebras
\Beq\label{8: morphism ftt to ICH}
g: \ftt \lon \ICH
\Eeq
given on the generators, for each $n\geq 1$, as follows}
\Beq\label{7: g_n+1: T_cyc to ICH}
\Ba{cccl}
g_{(n+1)}: & \ftt(({n+1})) & \lon & \wICH((n+1))\\
 &  T_{AB} & \lon & \ga_{AB}\ \ \forall A,B\in ((n+1)),\ A\neq B.
\Ea
\Eeq

\begin{proof} We have to check that for any $n\geq 1$ and any
$A,C,D\in ((n+1))$ one has
\Beq\label{8: condition on g}
g(\sum_{B=0}^n[T_{AB}, T_{CD}])=0.
\Eeq

{\sc Case 1}:  $\#\{A,C,D\}=3$.  Then
$$
0\equiv \sum_{B=0}^n[T_{AB}, T_{CD}]= [T_{AC}, T_{CD}] + [T_{AD}, T_{CD}]
$$
so that the required equality holds true if $g ([T_{AC}, T_{CD}])$ is skew-symmetric over indices $A,C,D$. This is indeed the case as
\Beqrn
g ([T_{AC}, T_{CD}])&=& \wICH\text{-part of }\
\delta\left(\Ba{c}\resizebox{17mm}{!}{ \xy
(-10,1)*+{_A}*\frm{o}="0";
(0,1)*+{_C}*\frm{o}="1";
(10,1)*+{_D}*\frm{o}="2";
 \ar @{-} "0";"1" <0pt>
  \ar @{-} "2";"1" <0pt>
\endxy} \Ea ...\right)\\
&=&
\underbrace{\Ba{c}\resizebox{17mm}{!}{ \xy
(-10,1)*+{_A}*\frm{o}="0";
(0,1)*+{_C}*\frm{o}="1";
(10,1)*+{_D}*\frm{o}="2";
  (0,10)*{\bu}="bu";
   (10,10)*{}="r";
%
 \ar @{-} "0";"bu" <0pt>
  \ar @{-} "1";"bu" <0pt>
  \ar @{-} "2";"bu" <0pt>
   \ar @{->} "bu";"r" <0pt>
\endxy} \Ea...}_{\delta^{(2)}_C} \ \ +
\underbrace{\Ba{c}\resizebox{17mm}{!}{ \xy
(-10,1)*+{_A}*\frm{o}="0";
(0,1)*+{_C}*\frm{o}="1";
(10,1)*+{_D}*\frm{o}="2";
  (0,14)*{\bu}="bu";
   (10,8)*{}="r";
 (0,8)*{*}="s";
 \ar @{.} "1";"bu" <0pt>
  \ar @{-} "0";"bu" <0pt>
  \ar @{-} "2";"bu" <0pt>
   \ar @{.>} "s";"r" <0pt>
\endxy} \Ea...}_{\delta_C^*} \ \
 \ \ +
\underbrace{\Ba{c}\resizebox{17mm}{!}{ \xy
(-10,1)*+{_C}*\frm{o}="0";
(0,1)*+{_A}*\frm{o}="1";
(10,1)*+{_D}*\frm{o}="2";
  (0,14)*{\bu}="bu";
   (10,8)*{}="r";
 (0,8)*{*}="s";
 \ar @{.} "1";"bu" <0pt>
  \ar @{-} "0";"bu" <0pt>
  \ar @{-} "2";"bu" <0pt>
   \ar @{.>} "s";"r" <0pt>
\endxy} \Ea...}_{\delta_A^{(1)}} \ \
 \ \ +
\underbrace{\Ba{c}\resizebox{17mm}{!}{ \xy
(-10,1)*+{_A}*\frm{o}="0";
(0,1)*+{_D}*\frm{o}="1";
(10,1)*+{_C}*\frm{o}="2";
  (0,14)*{\bu}="bu";
   (10,8)*{}="r";
 (0,8)*{*}="s";
 \ar @{.} "1";"bu" <0pt>
  \ar @{-} "0";"bu" <0pt>
  \ar @{-} "2";"bu" <0pt>
   \ar @{.>} "s";"r" <0pt>
\endxy} \Ea...}_{\delta_D^{(1)}} \ \
\Eeqrn
where dots stands for isolated white vertices (if any).
We conclude that $g ([T_{AC}, T_{CD}])$ is indeed skew-symmetric over $A,C,D$ so that
(\ref{8: condition on g}) holds true in Case 1 indeed.

\sip

{\sc Case 2}: $A=C$. Then
$$
g\left(\sum_{B=0}^n[T_{AB}, T_{AD}]\right)=
\Ba{c}\resizebox{17mm}{!}{ \xy
(-10,1)*+{_A}*\frm{o}="0";
(0,1)*+{...}="1";
(10,1)*+{_D}*\frm{o}="2";
  (0,10)*{\bu}="bu";
   (10,10)*{}="r";
%
 \ar @{-} "0";"bu" <0pt>
  \ar @{<-} "1";"bu" <0pt>
  \ar @{-} "2";"bu" <0pt>
   \ar @{->} "bu";"r" <0pt>
\endxy} \Ea \ \ +
\Ba{c}\resizebox{17mm}{!}{ \xy
(-10,1)*+{_A}*\frm{o}="0";
(0,1)*+{...}="1";
(10,1)*+{_D}*\frm{o}="2";
  (0,14)*{\bu}="bu";
   (10,8)*{}="r";
 (0,8)*{*}="s";
 \ar @{<.} "1";"bu" <0pt>
  \ar @{-} "0";"bu" <0pt>
  \ar @{-} "2";"bu" <0pt>
   \ar @{.>} "s";"r" <0pt>
\endxy} \Ea \ \
=0.
$$
Every term on the r.h.s.\ occurs twice with opposite signs so that the condition (\ref{8: condition on g}) holds true again.

The Proposition is proven. 
\end{proof}


\subsection{A plan to prove  the Main Theorem}

 Consider a cyclic operad of cochain complexes
$$
\odot^{\bu}(\ftt[1])= \{ \odot^{\bu}(\ftt((n+1))[1])\}_{n\geq 1}
$$
where  $\odot^{\bu}(\ftt((n+1))[1])$ is the Chevalley complex of the Lie algebra $\ftt((n+1))$.
It is well-known (see Lemma 1 in \cite{S}) that there is a quasi-isomorphism of cyclic operads
$$
G: \BVc \lon \odot^{\bu}(\ftt[1])
$$
which sends the multiplication generator of $\BVc$ into $1$ 
and the $BV$ operator $\Delta$ into the generator $T_{01}\in \ftt((2))[1]$ (see \S {\ref{2: subsec on T_cyc}}).

\sip

Composing $G$  with the morphism (\ref{8: morphism ftt to ICH}) we obtain another morphism of dg cyclic operads
$$
\BVc \stackrel{g\circ G}{\lon} \odot^{\bu}(\ICH[1])=\BVHgraphs
$$
We conclude that {\it if the morphism (\ref{8: morphism ftt to ICH}) is a quasi-isomorphism of $\Lie_\infty$-algebras, then the induced map $g\circ G: \BVc
\rar \BVHgraphs$ is a quasi-isomorphism as well, and the Main Theorem {\ref{7: Main Theorem}} follows}.

\sip

Hence to prove the Main theorem it is enough to prove that the morphism (\ref{8: morphism ftt to ICH}) is a quasi-isomorphism. It is more suitable to work with the $\Lie_\infty$-algebra $\wICH$ and the composition
\Beq\label{8: morphism ftt to wICH}
\hat{g}: \ftt \lon \ICH \hook \wICH
\Eeq
As 
$$
H^{\bu \leq 0}(\ICH)\equiv H^{\bu\leq 0}(\wICH)
$$
The morphism (\ref{8: morphism ftt to ICH}) is a quasi-isomorphism if and only the following Theorem holds true.

\subsection{Theorem}\label{8: qis on ftt and ICH} {\it 
The morphism (\ref{8: morphism ftt to wICH})  induces an isomorphism of cyclic operads of Lie algebras
$$
[\hat{g}]: \ftt \rar H^0(\wICH).
$$ 
Moreover $H^{\bu<0}(\wICH)=0$.} 

\bip

Here is the upshot of the above deliberations:
{\it once we prove Theorem {\ref{8: qis on ftt and ICH}}, the Main Theorem {\ref{7: Main Theorem}} follows immediately}.

\bip

{\large
\section{\bf Reducing the proof of Theorem {\ref{8: qis on ftt and ICH}}
to a simpler claim }
}

\sip

\subsection{``Zero vertex" filtrations of $\ftt((n+1))$ and $\wICH((n+1))$} We study
the following compatible filtrations  of both sides of the morphism of $\Lie_\infty$-algebras,
\Beq\label{9: g_n+1: T_cyc to wICH}
\Ba{cccl}
\hat{g}_{n+1}: & \ftt({n+1}) & \lon & \wICH(n+1)\\
 &  T_{AB} & \lon & \ga_{AB}=\Ba{c}\resizebox{19mm}{!}{ \xy
(-5,1)*+{_A}*\frm{o}="1";
 (5,1)*+{_B}*\frm{o}="2";
 (-15,1)*+{0}*\frm{o}="0";
  (15,1)*+{n}*\frm{o}="3";
  (-10,0)*{...};
  (10,0)*{...};
 \ar @{-} "2";"1" <0pt>
 \endxy} \Ea    \ \ \forall A,B\in ((n+1)),\ A\neq B,
\Ea
\Eeq

\Bi
\item[(i)] a filtration of the l.h.s.\ 
by the number of generators $T_{AB}$ with $A\neq 0$ and $B\neq 0$. The associated graded Lie algebra
is a direct sum of the usual Kohno-Drinfeld Lie algebra $t_n$ generated by $T_{ij}$, $i,j\in [n]$,
and the $n$-dimensional center spanned by $T_{0i}$, $i\in [n]$,
$$
gr \ftt({n+1}) = t_n\ \oplus\ \bigoplus_{i=1}^n\K\langle  T_{0i} \rangle.
$$

\item[(ii)] a decreasing and bounded filtration of the r.h.s.\ given by
$$
F_{-p}\wICH(n+1) :=\text{the subspace of $\ICH(n+1)$ generated by hypergraphs $\Ga$ with $||\Ga||\geq p$},
$$
where
 $$
||\Ga||= (\#\text{edges not connected to $\ \wo\ $}) + (\#\text{flags not connected to $\wo$}) -
\# V_\bu(\Ga) - \# V_*(\Ga).
$$
Denote by $(gr\wICH, d)$ be the associated graded complex.
\Ei

\sip

The morphism  (\ref{9: g_n+1: T_cyc to wICH}) respects both filtrations so that we get a morphism of associated graded complexes
\Beq\label{9: g_n+1: gr ftt(n+1) to gr wICH(n+1)} 
gr(\hat{g}_{(n+1)}):\ gr\ftt((n+1)) \lon  (gr\wICH((n+1)),d),
\Eeq
where the l.h.s.\ is assumed to have the trivial differential. We want to show by induction over $n$ that 
\Beq\label{9: H of gr wICH}
H^{\bu<0}((gr\wICH((n+1)),d)=0, \ \ \ H^{0}((gr\wICH((n+1)),d)\simeq gr\ftt((n+1)),
\Eeq
as this would imply Theorem {\ref{8: qis on ftt and ICH}}.

\sip

The induced  differential $d$ in  $gr\wICH(n+1)$ is the sum (cf.\ (\ref{5: full differential as a sum in BVGraphs}))
$$
d=d_\bu + d^{(2)}=\sum_{v \in V_\bu(\Ga)}d_v   + \sum_{i=1}^n d_i^{(2)}  
$$
where $d_v$ splits the black vertex  $v$ by creating a new dotted edge labelled by zero
(as the ``hanging" flag in the general formula (\ref{5: d on black in BVGraphs})
must be attached to $\ \wo\ $ in $gr\ICH(n+1)$ under $d_v$),
\Beq\label{9: d_bu in grICH}
d_v(\bu)= \sum \Ba{c}\resizebox{12mm}{!}{ \xy
(-4,1)*{\bu}="1";
 (4,1)*{\bu}="2";
  (0,3)*{_0};
 \ar @{.} "2";"1" <0pt>
 \endxy} \Ea.    
\Eeq
 The differential
$d_i^{(2)}$ acts on the white vertex $\ \wi\ $ non-trvially only if $\ \wi\ $ has at least one type-2 flag attached,
\Beq\label{7: d_0i}
d_i^{(2)}(\ \wi\ )=
\sum_{f\in Fl_i^2}
\resizebox{14mm}{!}{ \xy
(0,-3)*+{_0}*\frm{o}="0";
(10,-3)*+{_i}*\frm{o}="i";
 (5,1)*{\bu}="1";
 (5,7)*{}="l";
 (7,4)*{_f};
 \ar @{-} "1";"0" <0pt>
 \ar @{-} "1";"i" <0pt>
  \ar @{.} "l";"1" <0pt>
\endxy}\ ,
\Eeq
where $Fl_i^2$ is the set of all  type-2 flags attached to  $\ \wi\ $.

\sip

The parts $\delta^*$ and $\delta^{(1)}$ of the full differential (\ref{5: full differential as a sum in BVGraphs}) act trivially on internally connected hypergraphs as otherwise they would create hypergraphs with internal genus. 

\subsection{Inductive decompositions of  $gr\ft((n+1))$ and $gr\wICH((n+1))$}
Recall (see \S {\ref{2: subsec on T_cyc}}) that for each $n\geq 1$ the Lie algebra $\ft^c((n+1))$ fits a short exact sequence of Lie algebras
$$
0\lon \ft_n((n+1)) \lon \ftt((n+1)) \lon \ftt((n)) \lon 0
$$
where $\ft_n((n+1))$ is the subspace generated by elements $T_{An}$, $A\in {0,1,\dots,n-1}$. The Lie subalgebra $\ftt((n+1))$ can be identified with   the quotient of the free Lie algebra
$\mathfrak{freeLie}\langle T_{0n},T_{1,n}, \ldots, T_{n-1\,n}  \rangle$  modulo one relation saying that the sum
$$
T_{0n}+ T_{1n} + \ldots + T_{n-1\, n}
$$
is  central. Hence we can identify the associated graded with the following Lie algebra
$$
gr\ft_n(n+1)=\simeq  \mathfrak{freeLie} \langle T_{1n}, T_{2n}, \ldots, T_{n-1,n} \rangle
\ \ \oplus \ \ \K\langle T_{0n}\rangle, \ \ \ \forall n\geq 1,
$$
the element $T_{0n}$ being central.

\sip

Similarly, let $\wICH_n(n+1)$ be a subcomplex spanned by hypergraphs having at least one edge or flag attached to the vertex $\ \wn\ $. We have a short exact sequence 
of $\Lie_\infty$ algebras
$$
0\lon \wICH_n((n+1)) \lon  \wICH((n+1)) \lon \wICH((n)) \lon 0.
$$

The morphism (\ref{9: g_n+1: gr ftt(n+1) to gr wICH(n+1)}) restricts to the following 
morphism of the associated graded complexes (in fact, of $\Lie_\infty$-algebras)
$$
\hat{g}_n(n+1):  \ft_n((n+1)) \lon \wICH_n((n+1)).
$$

Our next purpose is to prove the following 

\subsection{Theorem}\label{9: Theorem on H(ICH_n)} {\it The morphism $gr(g_n(n+1))$
induces an isomorphism of Lie algebras,
$$
H^0(gr \wICH_n((n+1))) \simeq {gr}\,\ft_n((n+1))\simeq  \mathfrak{freeLie} \langle T_{1n}, T_{2n}, \ldots, T_{n-1,n} \rangle
\ \ \oplus \ \ \K\langle T_{0n}\rangle   \ \ \forall n\geq 1.
$$
Moreover $H^{\bu<0}(gr\wICH_n((n+1)))=0$ for all $n\geq 1$.}

\mip

This Theorem implies the equalities  (\ref{9: H of gr wICH})  which in turn imply
Theorem {\ref{8: qis on ftt and ICH}}.

\sip

The upshot of this section is that {\it once we prove {Theorem} {\ref{9: Theorem on H(ICH_n)}}, the Main Theorem {\ref{7: Main Theorem}} follows immediately}.

\bip


{\large
\section{\bf Proof of Theorem {\ref{9: Theorem on H(ICH_n)}}}
}

\sip

\subsection{\sc Step 1: the case $n=1$} This ``smallest" possible case illustrates some tricks used in the study of the general case. The complex   
$gr\ICH_1(1+1)$ decomposes into a direct sum of the trivial 1-dimensional subcomplex 
$$
gr\wICH_1(1+1)= \K\langle   \Ba{c}\resizebox{11mm}{!}{ \xy
(-5,1)*+{_0}*\frm{o}="0";
 (5,1)*+{_1}*\frm{o}="1";
 \ar @{-} "1";"0" <0pt>
\endxy} \Ea       \rangle \ \ \oplus \ \ \ gr\wICH_1(1+1)'
$$
and its complement $gr\wICH_1(1+1)'$ which is spanned by hypergraphs having at least one black vertex. We have to show that $gr\ICH_1(1+1)'$ is acyclic in non-positive degrees.
\sip

The complex 
 $gr\ICH_1(1+1)'$ is spanned by unrooted at least trivalent trees 
 of the form
  $$
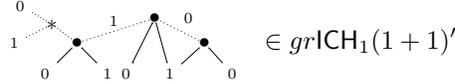
\Ba{c}\resizebox{33mm}{!}{ \xy
(-10,1)*+{_0}="0";
(0,1)*+{_1}="1";
(10,1)*+{_1}="2";
(3,1)*+{_0}="3";
(20,1)*+{_0}="n";
  (-5,6)*{\bu}="bu";
   (7.5,10)*{\bu}="r";
     (15.5,6)*{\bu}="rr";
(1,9.5)*{_1};
(13,10)*{_0};
(-9,9)*{*}="s";
(-14,12)*+{_0}="f1";
(-15,6)*+{_1}="f2";
 \ar @{-} "0";"bu" <0pt>
  \ar @{-} "1";"bu" <0pt>
   \ar @{.} "bu";"r" <0pt>
    \ar @{-} "rr";"2" <0pt>
     \ar @{-} "rr";"n" <0pt>
  \ar @{.} "rr";"r" <0pt>
 \ar @{-} "r";"2" <0pt>
 \ar @{-} "r";"3" <0pt>
 \ar @{.} "s";"bu" <0pt>
 \ar @{.} "s";"f1" <0pt>
 \ar @{.} "s";"f2" <0pt>
\endxy} \Ea \in gr\ICH_1(1+1)'
$$
 More formally, a generator of $gr\ICH_1(1+1)'$ is a tree $T$ with at least one black vertex such that
 \Bi
 \item[(i)] the black vertices of $T$ are all at least trivalent;
 \item[(ii)] internal edges (if any) of $T$ are dotted edges decorated by symbols $0$ or $1$ (see (\ref{7: type1 hyperedge}) for a pictorial explanation);
 \item[(iii)]  the labelled legs  of $T$ are either solid edges decorated by numbers from the set $\{0,1\}$ or hyperedges of type-2 whose type-2 flags  are decorated by symbols $0$ and $1$ (these hyperedges  can be understood as dotted hairs labelled by the pair of indices $(0,1)$ as in the picture above). Call the {\it internal}\, edges labelled by $1$ {\it bad}.
 \Ei
 
  Consider a filtration of the  complex $gr\wICH_1(1+1)'$ 
   by the total number of solid edges, and let $(E'_r, \p_r)$ be the associated (regular bounded) spectral sequence. The initial page $(E'_0, \p_0)$ has the differential
 $\p_0=\delta_\bu$ which creates internal edges labelled by zero.
The set 
of bad internal edges of generators $T$ of $E_0'$ is invariant under the differential $d_\bu$ (in fact, it is invariant under the full differential $d$).
 Hence the complex $(E_0', d_\bu)$ decomposes into a direct sum of
 subcomplexes
 $$
 (E_0', d_\bu)=\bigoplus_{N\geq 0} (C_N, d_\bu),
 $$
 where $C_N$ is the subcomplex spanned by graphs with precisely $N$ bad internal edges.
 
 \sip
 
 For each $N\geq 1$ we can consider a larger complex $\tilde{C}_N$ in which (i) the bad internal edges are distinguished (say totally ordered or, equivalently, labelled by a set of natural numbers $k\in [N]$), and (ii) each such $k$-numbered internal edge is  broken into to a pair of half-edges  which are distinguished  by some symbols, say labelled by $k'$ and $k''$, and assigned the cohomological degree $-1$ each). As 
 $$
 C_N=\left(\tilde{C}_N\right)_{\bS_N \rtimes(\bS_2)^N},
$$ 
we have, by  Maschke theorem, an isomorphism of cohomology groups
$$
H^\bu(C_N, d_\bu)= \left(H^\bu(\tilde{C}_N, d_\bu)\right)_{\bS_N \rtimes(\bS_2)^N}.
$$

  Breaking $N$ bad internal edges of generators of $C_N$ into half-edges decomposes the complex $\tilde{C}_N$ into $N+1$ connected components, and the differential $\delta_\bu$ preserves this decomposition. Hence
  $\tilde{C}_N$ can be identified with the unordered tensor product of $N+1$
  complexes, each of which is isomorphic to the complex $\tilde{C}_{connected}$
spanned by  unrooted trees $T''$ such that
   \Bi
   \item[(a)] {\it all} the internal edges of $T''$  (if any) are  labelled by zero, 
    \item[(c)] the legs of $T''$ can be of {\it three} types, the two types are as described in (iii) above, and the 3rd type are some labelled half-edges of the broken bad  edges. 
\Ei
 Thus the complex $C_0$ and each tensor factor $C_{connected}$ of  $\tilde{C}_{N\geq 1}$  is isomorphic essentially\footnote{Legs of graphs from $\Lie_\infty(m+1)$ are all labelled by {\it different}\, integers, while in our case some labels of legs can be repeated; such repeated labels can be understood as different ones but then symmetrized in one block (corresponding to a white vertex). Such a symmetrization commutes with the differential and hence does not affect our main conclusion about the structure of the cohomology groups.} isomorphic to the cyclic
 $\Lie_\infty$-operad  so that its cohomology is generated by trivalent unrooted trees modulo the $IHX$-relations. Hence we conclude that the next page  $E_1'=H^\bu(E_0', \delta_\bu)$ is generated by trivalent graphs of the form
   $$
\Ba{c}\resizebox{33mm}{!}{ \xy
(-10,1)*+{_0}="0";
(10,1)*+{_1}="2";
(3,1)*+{_0}="3";
(20,1)*+{_0}="n";
  (-5,6)*{\bu}="bu";
   (7.5,10)*{\bu}="r";
     (15.5,6)*{\bu}="rr";
(1,9.5)*{_1};
(13,10)*{_0};
(-9,9)*{*}="s";
(-14,12)*+{_0}="f1";
(-15,6)*+{_1}="f2";
 \ar @{-} "0";"bu" <0pt>
   \ar @{.} "bu";"r" <0pt>
    \ar @{-} "rr";"2" <0pt>
     \ar @{-} "rr";"n" <0pt>
  \ar @{.} "rr";"r" <0pt>
 \ar @{-} "r";"3" <0pt>
 \ar @{.} "s";"bu" <0pt>
 \ar @{.} "s";"f1" <0pt>
 \ar @{.} "s";"f2" <0pt>
\endxy} \Ea \in E_1'=H^\bu(E_0', \delta_\bu)
$$
modulo the $IHX$-relation associated with each internal edge labelled by zero.  The differential $\p_1$ in $E_1'$ is 
equal to $d^{(2)}$; it acts non-trivially only on  typer-2 flags of the trivalent graphs from $E_1'$ by changing them as follows 
$$
 d^{(2)}:\ \ \ \Ba{c}\resizebox{12mm}{!}{ \xy
(-7,0)*+{_0}="0";
(7,0)*+{_1}="1";
(0,12)*{\bu}="2";
  (0,5)*{*}="s";
 \ar @{.} "s";"2" <0pt> 
 \ar @{.} "s";"0" <0pt>
 \ar @{.} "s";"1" <0pt>
\endxy} \Ea 
\rar 
 \Ba{c}\resizebox{12mm}{!}{ \xy
(-7,0)*+{_0}="0";
(7,0)*+{_1}="1";
(0,12)*{\bu}="2";
  (0,5)*{\bu}="s";
(2,8)*{_0};
 \ar @{.} "s";"2" <0pt> 
 \ar @{-} "s";"0" <0pt>
 \ar @{-} "s";"1" <0pt>
\endxy} \Ea 
$$ 
Hence the complex $(E_1',\delta^{(2)})$ can be identified with the tensor product of a trivial complex with the unordered
tensor product of acyclic 2-dimensional complexes $\K[-1]\oplus \K$ whose generators are shown pictorially  just above. We conclude that the cohomology group
$H^\bu(E_1', \delta^{(2)}) \simeq gr H^\bu(gr\wICH_1(1+1)')$, if non-zero, is spanned by hypergraphs with no hyperedges of type-2, and no trivalent black vertices of the them 
  $$
   \Ba{c}\resizebox{12mm}{!}{ \xy
(-7,0)*+{_0}="0";
(7,0)*+{_1}="1";
(0,12)*{}="2";
  (0,5)*{\bu}="s";
(2,8)*{_0};
 \ar @{.} "s";"2" <0pt> 
 \ar @{-} "s";"0" <0pt>
 \ar @{-} "s";"1" <0pt>
\endxy} \Ea 
$$ 
Every such a hypergraph sits in cohomological degree +1 implying the vanishing of non-positively graded cohomology groups $H^{\bu\leq 0}(gr\wICH_1(1+1))'=0$. Hence  $H^{\bu<0}(gr\wICH_1(1+1))=0$ and 
$$
H^{0}(gr\wICH_1(1+1))
= \K\langle   \Ba{c}\resizebox{11mm}{!}{ \xy
(-5,1)*+{_0}*\frm{o}="0";
 (5,1)*+{_1}*\frm{o}="1";
 \ar @{-} "1";"0" <0pt>
\endxy} \Ea  \rangle = \K\langle T_{01}\rangle
$$
as claimed in Theorem {\ref{9: Theorem on H(ICH_n)}} in the case $n=1$.

\subsection{\sc Step 2: acyclcity of a ``type-0" subcomplex} We study a subcomplex  $C_{type0}$  (in fact, a direct summand) of $(gr\wICH_n, d)$  which is the generated by the hyperedges of type-0. More precisely,
Let $C_{type0}=\sum_{i=-2}^0 C_{type0}^i$ be a $\Z$-graded subspace of $\ICH_n$ defined as follows:

$$
C_{type0}^{-1}=\K \left\langle \ga_{ijn}:=
  \Ba{c}\resizebox{28mm}{!}{ \xy
(-23,0)*+{_0}*\frm{o}="0";
(-7.5,0)*{...};
(+7.5,0)*{...};
(-17.5,0)*{...};
(-12,0)*+{_i}*\frm{o}="0";
(0,0)*+{_j}*\frm{o}="1";
(12,0)*+{_n}*\frm{o}="2";
  (0,9)*{*}="s";
 \ar @{.} "s";"2" <0pt>
 \ar @{.} "s";"0" <0pt>
 \ar @{.} "s";"1" <0pt>
\endxy} \Ea, \ \ \ 0\leq i<j<n
   \right\rangle,  
$$
 $$
C_{type0}^{0}=\K \left\langle 
 d_i^{(2)} \ga_{ijn},  d_j^{(2)} \ga_{ijn},  d_n^{(2)} \ga_{ijn} \ \ \ 0\leq i<j<n
   \right\rangle,  
$$

 $$
C_{type0}^{1}=\delta \left( C_{type0}^{0}\right) =\delta^{(2)} \left( C_{type0}^{0}\right) 
$$

\subsubsection{\bf Lemma} {\it The complex $C_{type0}$ is acyclic.}
\begin{proof}
The action of the differential $d$ on $C_{type0}^{-1}$ is obviously an injection  $d:C_{type0}^{-1}\hook  C_{type0}^{0}$, so to prove the Lemma it is enough to show that $H^{0}(C_{type0})=0$. 

\sip

The complex $C_{type0}$ decomposes into a direct sum of complexes
$$
 C_{type0}=\bigoplus_{ 0\leq i<j<n} C_{ijn}
$$
and we have to consider two cases, $i=0$ and $i\neq 0$. In the  first case we have
$$
C_{ijn}^{0}=\K \left\langle 
\Ba{c}\resizebox{21mm}{!}{ \xy
(-10,1)*+{_0}*\frm{o}="0";
(0,1)*+{_j}*\frm{o}="j";
(10,1)*+{_n}*\frm{o}="n";
  (0,12)*{*}="s";
   (-2,6)*{\bu}="bu";
%
 \ar @{.} "0";"s" <0pt>
 \ar @{.} "bu";"s" <0pt>
 \ar @{.} "n";"s" <0pt>
 \ar @{-} "bu";"0" <0pt>
  \ar @{-} "bu";"j" <0pt>
\endxy} \Ea
, \ \ \
\Ba{c}\resizebox{21mm}{!}{ \xy
(-10,1)*+{_0}*\frm{o}="0";
(0,1)*+{_j}*\frm{o}="j";
(10,1)*+{_n}*\frm{o}="n";
  (0,12)*{*}="s";
   (3,6)*{\bu}="bu";
%
 \ar @{.} "0";"s" <0pt>
 \ar @{.} "bu";"s" <0pt>
 \ar @{.} "j";"s" <0pt>
 \ar @{-} "bu";"0" <0pt>
  \ar @{-} "bu";"n" <0pt>
\endxy} \Ea
  \right\rangle,
$$
The only cycle of degree $0$ in this complex is the sum
(with the appropriate ordering of edges and flags)
$$
\Ba{c}\resizebox{21mm}{!}{ \xy
(-10,1)*+{_0}*\frm{o}="0";
(0,1)*+{_j}*\frm{o}="j";
(10,1)*+{_n}*\frm{o}="n";
  (0,12)*{*}="s";
   (-2,6)*{\bu}="bu";
%
 \ar @{.} "0";"s" <0pt>
 \ar @{.} "bu";"s" <0pt>
 \ar @{.} "n";"s" <0pt>
 \ar @{-} "bu";"0" <0pt>
  \ar @{-} "bu";"j" <0pt>
\endxy} \Ea
+
\Ba{c}\resizebox{21mm}{!}{ \xy
(-10,1)*+{_0}*\frm{o}="0";
(0,1)*+{_j}*\frm{o}="j";
(10,1)*+{_n}*\frm{o}="n";
  (0,12)*{*}="s";
   (3,6)*{\bu}="bu";
%
 \ar @{.} "0";"s" <0pt>
 \ar @{.} "bu";"s" <0pt>
 \ar @{.} "j";"s" <0pt>
 \ar @{-} "bu";"0" <0pt>
  \ar @{-} "bu";"n" <0pt>
\endxy} \Ea
$$
which is equal precisely to 
$d \Ba{c}\resizebox{16mm}{!}{ \xy
(-10,1)*+{_0}*\frm{o}="0";
(0,1)*+{_j}*\frm{o}="j";
(10,1)*+{_n}*\frm{o}="n";
  (0,12)*{*}="s";
%
 \ar @{.} "0";"s" <0pt>
 \ar @{.} "n";"s" <0pt>
 \ar @{.} "j";"s" <0pt>
\endxy} \Ea$. Hence every complex $C_{0jn}$ is acyclic.

In the second case one has
$$
C_{ijn}^{0}=\K \left\langle 
  \Ba{c}\resizebox{28mm}{!}{ \xy
(-23,0)*+{_0}*\frm{o}="0";
(-7.5,0)*{...};
(+7.5,0)*{...};
(-17.5,0)*{...};
(-12,0)*+{_i}*\frm{o}="i";
(0,0)*+{_j}*\frm{o}="j";
(12,0)*+{_n}*\frm{o}="n";
  (0,9)*{*}="s";
(-12,6)*{\bu}="bu";
 \ar @{.} "s";"bu" <0pt>
 \ar @{.} "s";"j" <0pt>
 \ar @{.} "s";"n" <0pt>
 \ar @{-} "bu";"0" <0pt>
  \ar @{-} "bu";"i" <0pt>
\endxy} \Ea, \ \ 
 \Ba{c}\resizebox{28mm}{!}{ \xy
(-23,0)*+{_0}*\frm{o}="0";
(-7.5,0)*{...};
(+7.5,0)*{...};
(-17.5,0)*{...};
(-12,0)*+{_i}*\frm{o}="i";
(0,0)*+{_j}*\frm{o}="j";
(12,0)*+{_n}*\frm{o}="n";
  (0,9)*{*}="s";
(-12,6)*{\bu}="bu";
 \ar @{.} "s";"bu" <0pt>
 \ar @{.} "s";"i" <0pt>
 \ar @{.} "s";"n" <0pt>
 \ar @{-} "bu";"0" <0pt>
  \ar @{-} "bu";"j" <0pt>
\endxy} \Ea, \ \
 \Ba{c}\resizebox{28mm}{!}{ \xy
(-23,0)*+{_0}*\frm{o}="0";
(-7.5,0)*{...};
(+7.5,0)*{...};
(-17.5,0)*{...};
(-12,0)*+{_i}*\frm{o}="i";
(0,0)*+{_j}*\frm{o}="j";
(12,0)*+{_n}*\frm{o}="n";
  (0,9)*{*}="s";
(-12,6)*{\bu}="bu";
 \ar @{.} "s";"bu" <0pt>
 \ar @{.} "s";"j" <0pt>
 \ar @{.} "s";"i" <0pt>
 \ar @{-} "bu";"0" <0pt>
  \ar @{-} "bu";"n" <0pt>
\endxy} \Ea, 
   \right\rangle
$$
Again it is easy to see that the space of degree $0$ cycles is one-dimensional and is equal to 
$d  \Ba{c}\resizebox{28mm}{!}{ \xy
(-23,0)*+{_0}*\frm{o}="0";
(-7.5,0)*{...};
(+7.5,0)*{...};
(-17.5,0)*{...};
(-12,0)*+{_i}*\frm{o}="0";
(0,0)*+{_j}*\frm{o}="1";
(12,0)*+{_n}*\frm{o}="2";
  (0,9)*{*}="s";
 \ar @{.} "s";"2" <0pt>
 \ar @{.} "s";"0" <0pt>
 \ar @{.} "s";"1" <0pt>
\endxy} \Ea$. Hence every  complex $C_{ijn}$ with $0<i<j<n$ is also acyclic.
\end{proof}

\sip

The complex $gr\wICH_n$ decomposes into the {\it direct}\, sum of complexes
$$
gr\wICH_n= C_{type0} \oplus \overline{\ICH}_n.
$$
Indeed, the only possible obstruction to the direct summation above may come  from the hypergraphs from $d^{(2)}(C_{type0}^{0}))$ 
which have two black vertices and the internal dotted edge {\it labelled by zero}; in principle such hypergraphs may come as a result of applying $d_\bu$ to some generator from the complement $\overline{\ICH}_n$ (as a graded vector space), i.e.\ by splitting some black vertex of a hypergraph from  $\overline{\ICH}$. However a moment look at such possibly problematic hypergraphs, say this one
$$
 \Ba{c}\resizebox{28mm}{!}
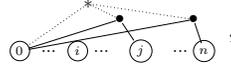
{ \xy
(-23,0)*+{_0}*\frm{o}="0";
(-7.5,0)*{...};
(+7.5,0)*{...};
(-17.5,0)*{...};
(-12,0)*+{_i}*\frm{o}="i";
(0,0)*+{_j}*\frm{o}="j";
(12,0)*+{_n}*\frm{o}="n";
  (-10,9)*{*}="s";
(10,6)*{\bu}="bu";
(-4,6)*{\bu}="bu2";
 \ar @{.} "s";"bu" <0pt>
 \ar @{.} "s";"bu2" <0pt>
  \ar @{-} "0";"bu2" <0pt>
  \ar @{-} "j";"bu2" <0pt>
 \ar @{.} "s";"0" <0pt>
 \ar @{-} "bu";"0" <0pt>
  \ar @{-} "bu";"n" <0pt>
\endxy} \Ea,
$$
convinces us that this is impossible as {\it both}\, black vertices have an edge connected to $\wo$.

\sip

Next we study the other direct summand  $\overline{\ICH}_n$ .

\subsection{\sc Step 3: Reducing $\overline{\ICH}_n((n+1))$, $n\geq 2$, to a smaller quasi-isomorphic subcomplex}

The  differential $d_\bu$ in $\overline{\ICH}_n((n+1))$ creates only dotted edges labelled by zero  (see (\ref{9: d_bu in grICH})) while $d^{(2)}$ creates dotted edges labelled by any integer $i\in \{1,2,\ldots, n\}$ when acting on type-2 flags attached to some white vertices $\wi$ and $\wj$ with $0<i<j\leq n$
$$
(d_i^{(2)} + d_j^{(2)})
 \Ba{c}\resizebox{24mm}{!}{ \xy
 (-15,0)*+{_0}*\frm{o}="0";
(-7,0)*+{_i}*\frm{o}="i";
(7,0)*+{_j}*\frm{o}="j";
(0,12)*{\bu}="bu";
  (0,5)*{*}="s";
 \ar @{.} "s";"bu" <0pt> 
 \ar @{.} "s";"i" <0pt>
 \ar @{.} "s";"j" <0pt>
\endxy} \Ea 
\lon 
\Ba{c}\resizebox{24mm}{!}{ \xy
 (-19,0)*+{_0}*\frm{o}="0";
(-7,0)*+{_i}*\frm{o}="i";
(7,0)*+{_j}*\frm{o}="j";
(0,12)*{\bu}="bu";
  (-9,5)*{\bu}="s";
(-3,8)*{_j};
 \ar @{.} "s";"bu" <0pt> 
 \ar @{-} "s";"i" <0pt>
 \ar @{-} "s";"0" <0pt>
\endxy} \Ea 
\ + \ 
 \Ba{c}\resizebox{24mm}{!}{ \xy
 (-19,0)*+{_0}*\frm{o}="0";
(-7,0)*+{_i}*\frm{o}="i";
(7,0)*+{_j}*\frm{o}="j";
(0,12)*{\bu}="bu";
  (-3,5)*{\bu}="s";
(0,8)*{_i};
 \ar @{.} "s";"bu" <0pt> 
 \ar @{-} "s";"0" <0pt>
 \ar @{-} "s";"j" <0pt>
\endxy} \Ea , \ \ \ \  \ i>0, j>0.
$$ 
Let us label such a type-2 hyperedge in the left hand side of the above formula by the corresponding pair of different {\it positive}\, numbers $(i,j)$ and call it {\it bad}\,\footnote{
Note that the differential $d^{(2)}$, when acting of type2 hyperedges which are not bad (i.e. the ones which have one flag attached to $\wo$),  creates only dotted edges labelled by {\it zero},
$
 d_i^{(2)}:
 \Ba{c}\resizebox{12mm}{!}{ \xy
(-7,0)*+{_0}*\frm{o}="i";
(7,0)*+{_i}*\frm{o}="j";
(0,12)*{\bu}="bu";
  (0,5)*{*}="s";
 \ar @{.} "s";"bu" <0pt> 
 \ar @{.} "s";"i" <0pt>
 \ar @{.} "s";"j" <0pt>
\endxy} \Ea 
\lon 
\Ba{c}\resizebox{11mm}{!}{ \xy
 (0,0)*+{_0}*\frm{o}="0";
(12,0)*+{_i}*\frm{o}="i";
(6,12)*{\bu}="bu";
  (6,5)*{\bu}="s";
(7.6,8)*{_0};
\ar @{.} "s";"bu" <0pt> 
 \ar @{-} "s";"i" <0pt>
 \ar @{-} "s";"0" <0pt>
\endxy} \Ea
$.
}.

\sip

We call a hypergraph $\Ga\in   \overline{\ICH}_n((n+1))$ {\it redundant}\, if it has at least one internal edge labelled by $j>0$
or at least one bad hyperedge.  Such hypergraphs span
a subcomplex $\ICH_n((n+1))^{rednt}\subset \overline{\ICH}_n((n+1))$
which is in fact a direct summand,
$$
\overline{\ICH}_n((n+1))= \overline{\ICH}_n((n+1))^{rednt}\ \ \oplus \ \ \overline{\ICH}_n((n+1))^{nice},
$$
where $\overline{\ICH}_n((n+1))^{nice}$ is  spanned by hypergraphs whose 
 internal edges (if any) are {\it all}\, labelled by $0$ and whose type-2 hyperedges (if any) have one flag connected to $\wo$. We shall see below that this complex can be understood as generated by ordinary graphs, not by hypergraphs; in this sense it is {\it nice}\, indeed.

\mip

\subsubsection{\bf Lemma} {\it The complex $\overline{\ICH}_n((n+1))^{rednt}$ has trivial cohomology in non-positive degrees, i.e.\  $H^{\bu\leq 0}\overline{\ICH}_n((n+1))^{rednt}=0$ for any $n\geq 2$.}

\begin{proof}
Let us call a trivalent black vertex of the form
$$
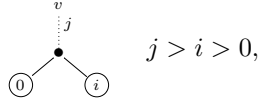
\Ba{c}\resizebox{14mm}{!}{ \xy
 (0,0)*+{_0}*\frm{o}="0";
(12,0)*+{_i}*\frm{o}="i";
(6,12)*{^v}="bu0";
  (6,5)*{\bu}="s";
(7.6,9.6)*{_{j}};
\ar @{.} "s";"bu0" <0pt> 
 \ar @{-} "s";"i" <0pt>
 \ar @{-} "s";"0" <0pt>
\endxy} \Ea \ \ \  j>i>0, 
$$
{\it a bad vertex  of index $(i,j)$}. It is connected to some other black vertex $v$ via a $j$-labelled dotted edge with $j>i$; if that black vertex $v$ is also bad (of any bi-index, $(i,j)$ including), then we get an internally 
connected hypergraph of cohomological degree $+1$ which is of no concern to us.
Hence we can proceed with an assumption that  every {\it bad}\, black vertex is connected by its unique
$j$-labelled dotted edge to a  black vertex $v$ which is {\it not}\, bad. For symmetry reasons the black vertex $v$ can not have attached more that $1$ such bad vertex with the same index $(i,j)$ as otherwise the hypergraph identically vanishes,
  $$
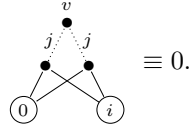
\Ba{c}\resizebox{16mm}{!}{ \xy
 (0,0)*+{_0}*\frm{o}="0";
(12,0)*+{_i}*\frm{o}="i";
(6,12)*{\bu}="bu0";
  (3,6)*{\bu}="s1";
(9,6)*{\bu}="s2";
(6,13.7)*{^v};
(8.9,9.6)*{_{j}};
(3.6,9.6)*{_{j}};
\ar @{.} "s1";"bu0" <0pt>
 \ar @{-} "s1";"i" <0pt>
 \ar @{-} "s1";"0" <0pt>
 \ar @{.} "s2";"bu0" <0pt>
 \ar @{-} "s2";"i" <0pt>
 \ar @{-} "s2";"0" <0pt>
\endxy} \Ea \equiv 0.
$$

Consider a filtration of $\overline{\ICH}_n((n+1))^{rednt}$ by
the number of non-bad black vertices, and let $(\cE^{rednt}_r, \p_r)$
be the associated spectral sequence. The initial page $(\cE_0^{rednt}, \p_0)$
has the differential $\p_0$ acting non-trivially only on bad hyperedges (if any) by creating a bad black vertex as follows
$$
 \Ba{c}\resizebox{24mm}{!}
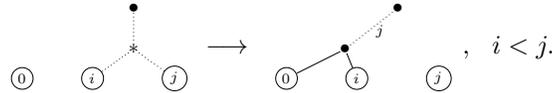
{ \xy
 (-19,0)*+{_0}*\frm{o}="0";
(-7,0)*+{_i}*\frm{o}="i";
(7,0)*+{_j}*\frm{o}="j";
(0,12)*{\bu}="bu";
  (0,5)*{*}="s";
 \ar @{.} "s";"bu" <0pt> 
 \ar @{.} "s";"i" <0pt>
 \ar @{.} "s";"j" <0pt>
\endxy} \Ea 
\lon 
\Ba{c}\resizebox{24mm}{!}{ \xy
 (-19,0)*+{_0}*\frm{o}="0";
(-7,0)*+{_i}*\frm{o}="i";
(7,0)*+{_j}*\frm{o}="j";
(0,12)*{\bu}="bu";
  (-9,5)*{\bu}="s";
(-3,8)*{_j};
 \ar @{.} "s";"bu" <0pt> 
 \ar @{-} "s";"i" <0pt>
 \ar @{-} "s";"0" <0pt>
\endxy} \Ea , \ \ i<j.
$$

The differential preserves the number of non-bad black vertices in $\cE^{rednt}_0$  so we can consider a larger complex $\widehat{\cE}^{rednt}_0$
spanned by hypergraphs in which non-bad black vertices are distinguished, say labelled by some integers $k\in \N$. Then every bad black vertex is uniquely determined by its index $(i,j,k)$, where $i<j$ is its original index and $k$ is the index of the unique non-special black vertex $v$ to which it is attached. Using the lexicographical order on the triples $(i,j,k)$ we obtain a well-defined notion of the {\it minimal bad vertex}; call it {\it very bad}. Consider next a filtration of $\widehat{\cE}^{rednt}_0$  by the number
of black vertices  which are not {\it very bad}. The associated graded decomposes into a direct sum
$$
gr\widehat{\cE}^{rednt}_0= {A}'_0\ \ \oplus \ \  C\ot A_0''
$$ 
where ${A}'_0$ is the trivial subcomplex spanned by hypergraphs with no bad hyperedges and no black vertices, the complex $C$ is a 2-dimensional acyclic complex of the form
$$
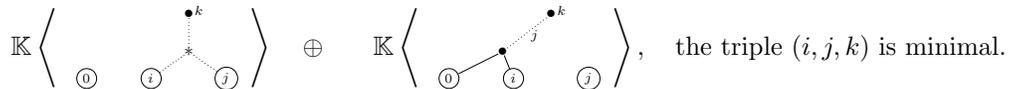
\K\left\langle  \Ba{c}\resizebox{22mm}{!}{ \xy
 (-19,0)*+{_0}*\frm{o}="0";
(-7,0)*+{_i}*\frm{o}="i";
(7,0)*+{_j}*\frm{o}="j";
(0,12)*{\bu}="bu";
(2,12)*{^k};
  (0,5)*{*}="s";
 \ar @{.} "s";"bu" <0pt> 
 \ar @{.} "s";"i" <0pt>
 \ar @{.} "s";"j" <0pt>
\endxy} \Ea \right\rangle 
\ \ \ \oplus\ \ \  \ \ 
\K\left\langle 
\Ba{c}\resizebox{22mm}{!}{ \xy
 (-19,0)*+{_0}*\frm{o}="0";
(-7,0)*+{_i}*\frm{o}="i";
(7,0)*+{_j}*\frm{o}="j";
(0,12)*{\bu}="bu";
(2,12)*{^k};
  (-9,5)*{\bu}="s";
(-3,8)*{_j};
 \ar @{.} "s";"bu" <0pt> 
 \ar @{-} "s";"i" <0pt>
 \ar @{-} "s";"0" <0pt>
\endxy} \Ea \right\rangle, 
\ \ \ \text{the triple}\ (i,j,k) \ \text{is minimal}.
$$
and $A_0''$ is a trivial complex. We conclude that the next page $({\cE}^{rednt}_1,\p_1)$ of the spectral sequence is equal to $A_0'$, i.e.\ it is generated by hypergraphs with no bad hyperedges and no bad black vertices, but with at least one
internal dotted edge labelled by $j>0$.

\sip

Let us consider next a filtration of $(\cE^{rednt}_1,\p_1)$ by the number of solid edges, and let $(E_r\cE_1^{rednt}, \p_{1r})$ be the associated spectral sequence.
The induced differential $\p_{10}$ in $E_0\cE^{rednt}_1$ is just $d_\bu$.
From now one we understand a connected to $v$ subgraph of the form 
$$
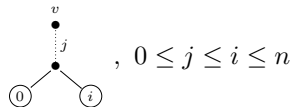
\Ba{c}\resizebox{13mm}{!}{ \xy
 (0,0)*+{_0}*\frm{o}="0";
(12,0)*+{_i}*\frm{o}="i";
(6,12)*{\bu}="bu";
(6,14)*{^v};
  (6,5)*{\bu}="s";
(7.6,8.6)*{_{j}};
\ar @{.} "s";"bu" <0pt> 
 \ar @{-} "s";"i" <0pt>
 \ar @{-} "s";"0" <0pt>
\endxy} \Ea, 
 \ 0\leq j \leq i\leq n
$$
as an attached to $v$ a {\it wavy leg} labelled by a pair of integers $(j\leq i)$ as in the following picture,
$$
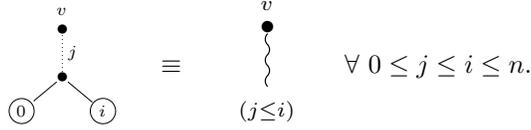
\Ba{c}\resizebox{15mm}{!}{ \xy
 (0,0)*+{_0}*\frm{o}="0";
(12,0)*+{_i}*\frm{o}="i";
(6,12)*{\bu}="bu";
(6,14)*{^v};
  (6,5)*{\bu}="s";
(7.6,8.6)*{_{j}};
\ar @{.} "s";"bu" <0pt> 
 \ar @{-} "s";"i" <0pt>
 \ar @{-} "s";"0" <0pt>
\endxy} \Ea 
\ \ \ \equiv \ \ \ 
\Ba{c}\resizebox{9mm}{!}{ \xy
(6,12)*{\bu}="bu";
(6,14)*{^v};
  (6,5)*{}="s";
(6,2)*{_{(j\leq i)}};
\ar @{~} "s";"bu" <0pt> 
\endxy} \Ea \ \ \ \  \forall \ 0\leq j \leq  i\leq n.
$$
 The reason why it works  is that  the induced differentials $\p_{1r}$ act on such subgraphs trivially, i.e.\ they are preserved. From now on  the $j$-labelled dotted internal edges of the type shown in the left hand side of the above picture are {\it not}\, understood as {\it internal}\, edges anymore, they are understood as wavy legs attached to black vertices.

\sip

We now play exactly the same trick as in {\sc STEP 1}: we brake all {\it internal}\, edges labelled 
by integers $j>0$ into pairs of new legs, identify the resulting
unordered tensor product of complexes of unrooted trees with, essentially, cyclic
$\Lie_\infty$ complexes, and conclude that the next page $(E_1\cE^{rednt}_1, \p_{11})$ is spanned by connected trivalent trees (modulo the $IHX$ relation at each internal edge labelled by zero)
which have legs of three types: the solid ones labeled by integers from $0$ to $n$, non-bad hyperedges of type-2 indexed by pairs $(0,j)$, $j>0$, 
and wavy edges labelled by pairs $(j<i)$ with $j\geq 0$ as in the picture
   $$
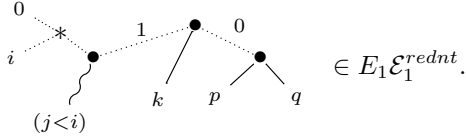
\Ba{c}\resizebox{43mm}{!}{ \xy
(-9,-2)*+{_{(j<i)}}="0";
(10,1)*+{_p}="2";
(3,1)*+{_k}="3";
(20,1)*+{_q}="n";
  (-5,6)*{\bu}="bu";
   (7.5,10)*{\bu}="r";
     (15.5,6)*{\bu}="rr";
(1,9.5)*{_1};
(13,10)*{_0};
(-9,9)*{*}="s";
(-14,12)*+{_0}="f1";
(-15,6)*+{_i}="f2";
 \ar @{~} "0";"bu" <0pt>
   \ar @{.} "bu";"r" <0pt>
    \ar @{-} "rr";"2" <0pt>
     \ar @{-} "rr";"n" <0pt>
  \ar @{.} "rr";"r" <0pt>
 \ar @{-} "r";"3" <0pt>
 \ar @{.} "s";"bu" <0pt>
 \ar @{.} "s";"f1" <0pt>
 \ar @{.} "s";"f2" <0pt>
\endxy} \Ea \in E_1\cE_1^{rednt}.
$$
The induced differential $\p_{11}$ acts non-trivially only on non-bad type-2 hyperedges as follows
$$
 \p_{11}:
 \Ba{c}\resizebox{14mm}{!}
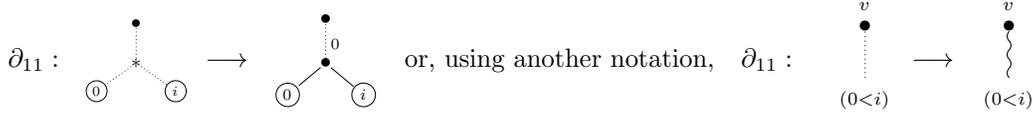
{ \xy
(-7,0)*+{_0}*\frm{o}="i";
(7,0)*+{_i}*\frm{o}="j";
(0,12)*{\bu}="bu";
  (0,5)*{*}="s";
 \ar @{.} "s";"bu" <0pt> 
 \ar @{.} "s";"i" <0pt>
 \ar @{.} "s";"j" <0pt>
\endxy} \Ea 
\lon 
\Ba{c}\resizebox{14mm}{!}{ \xy
 (0,0)*+{_0}*\frm{o}="0";
(12,0)*+{_i}*\frm{o}="i";
(6,12)*{\bu}="bu";
  (6,5)*{\bu}="s";
(7.6,8)*{_0};
\ar @{.} "s";"bu" <0pt> 
 \ar @{-} "s";"i" <0pt>
 \ar @{-} "s";"0" <0pt>
\endxy} \Ea \ \ \text{or, using another notation,}\ \ \ 
 \p_{11}:\ \
\Ba{c}\resizebox{8mm}{!}{ \xy
(6,12)*{\bu}="bu";
(6,14)*{^v};
  (6,5)*{}="s";
(6,2)*{_{(0<i)}};
\ar @{.} "s";"bu" <0pt> 
\endxy} \Ea
\lon
\Ba{c}\resizebox{8mm}{!}{ \xy
(6,12)*{\bu}="bu";
(6,14)*{^v};
  (6,5)*{}="s";
(6,2)*{_{(0<i)}};
\ar @{~} "s";"bu" <0pt> 
\endxy} \Ea
$$
This action changes only decorations of legs of trivalent vertices 
(``a dotted $(0,i)$-edge goes into a wavy $(0,i)$-edge"); hence 
the next and the final page  $E_2\cE_1^{rednt}= H^\bu(E_1\cE_1^{rough},\p_{11})$ of the spectral sequence is  generated by trivalent hypergraphs with no type-2 hyperedges at all; every such hypergraph has the cohomological degree equal to  +1. We conclude all the non-positively graded cohomology groups  $H^{\bu \leq 0}(\overline{\ICH}_n((n+1))^{rednt})$
vanish proving thereby the Lemma.

\end{proof}

\subsection{\sc Step 4: computing the cohomology of $\overline{\ICH}_n((n+1))^{nice}$   }
Recall that $\overline{\ICH}_n((n+1))^{nice}$ is spanned by hypergraphs $\Ga$ whose 
dotted internal edges (if any) are all labelled by $0$ so we omit such labels in the pictures from now on; moreover  type-2 hyperedges of $\Ga$ (if any) must have one flag connected to $\ \wo\ $. We identify every such a type-2 hyperedge with a dotted edge as explained  in the picture
$$
 \Ba{c}\resizebox{14mm}{!}
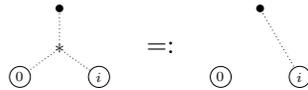
{ \xy
(-7,0)*+{_0}*\frm{o}="i";
(7,0)*+{_i}*\frm{o}="j";
(0,12)*{\bu}="bu";
  (0,5)*{*}="s";
 \ar @{.} "s";"bu" <0pt> 
 \ar @{.} "s";"i" <0pt>
 \ar @{.} "s";"j" <0pt>
\endxy} \Ea\ \  =: \
 \Ba{c}\resizebox{14mm}{!}{ \xy
(-7,0)*+{_0}*\frm{o}="0";
(7,0)*+{_i}*\frm{o}="i";
(0,12)*{\bu}="bu";
 \ar @{.} "i";"bu" <0pt> 
\endxy} \Ea  
$$
Thus the generators of $\overline{\ICH}_n((n+1))^{nice}$   are understood from now on as ordinary internally connected graphs of internal genus zero which have two types of edges, the solid edges and the dotted ones; moreover, every internal edge is dotted, and there are no dotted edges landing at $\ \wo\ $, e.g.
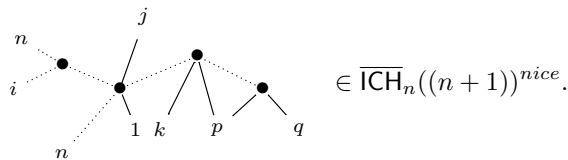
 $$
\Ba{c}\resizebox{43mm}{!}{ \xy
(-9,-2)*+{_{n}}="0";
(0,1)*+{_1}="1";
(10,1)*+{_p}="2";
(3,1)*+{_k}="3";
(20,1)*+{_q}="n";
  (-2,6)*{\bu}="bu";
   (7.5,10)*{\bu}="r";
     (15.5,6)*{\bu}="rr";
(1,14.5)*+{^j}="j";
%
(-9,9)*{\bu}="s";
(-14,12)*+{_n}="f1";
(-15,6)*+{_i}="f2";
 \ar @{.} "0";"bu" <0pt>
  \ar @{-} "1";"bu" <0pt>
  \ar @{-} "j";"bu" <0pt>
   \ar @{.} "bu";"r" <0pt>
    \ar @{-} "rr";"2" <0pt>
     \ar @{-} "rr";"n" <0pt>
  \ar @{.} "rr";"r" <0pt>
 \ar @{-} "r";"2" <0pt>
 \ar @{-} "r";"3" <0pt>
 \ar @{.} "s";"bu" <0pt>
 \ar @{.} "s";"f1" <0pt>
 \ar @{.} "s";"f2" <0pt>
\endxy} \Ea \in \overline{\ICH}_n((n+1))^{nice}.
$$
Here labels of legs  stand for the white vertices to which they are attached.
The differential in this complex is the sum  $d=d_\bu + d^{(2)}$
where $d_\bu$ splits black vertices as in (\ref{9: d_bu in grICH}) and
$d^{(2)}$ acts only on  labelled  {\it dotted legs}\, as  follows
\Beq\label{10: d(2) on dotted edge i}
 \Ba{c}\resizebox{14mm}{!}{ \xy
(-7,0)*+{_0}*\frm{o}="0";
(7,0)*+{_i}*\frm{o}="i";
(0,12)*{\bu}="bu";
 \ar @{.} "i";"bu" <0pt> 
\endxy} \Ea  
\lon 
\Ba{c}\resizebox{13mm}{!}{ \xy
 (0,0)*+{_0}*\frm{o}="0";
(12,0)*+{_i}*\frm{o}="i";
(6,12)*{\bu}="bu";
  (6,5)*{\bu}="s";
\ar @{.} "s";"bu" <0pt> 
 \ar @{-} "s";"i" <0pt>
 \ar @{-} "s";"0" <0pt>
\endxy} \Ea.
\Eeq

The Theorem is proven once we show
\Beq\label{9: iso H(ICH_n) nice and freeLie}
H^{\bu<0}(\overline{\ICH}_n((n+1))^{nice})=0, \ \ H^0(\overline{\ICH}_n((n+1))^{nice})\simeq  \mathfrak{freeLie} \langle T_{1n}, T_{2n}, \ldots, T_{n-1,n} \rangle
\  \oplus \  \K\langle T_{0n}\rangle \ \  \forall\ n\geq 2.
\Eeq

\sip

 The total number $N\geq 1$ of solid and dotted edges attached to the vertex $\ \wn\ $ stays invariant under the differential $d$, only their type may change from ``dotted" color to ``solid" one under the action of $d$. Hence the complex decomposes into a direct sum
$$
\overline{\ICH}_n((n+1))^{nice}=\bigoplus_{N\geq 1} \overline{\ICH}_n((n+1))^{nice}_N
$$ 
when $N$ is the total valency of the white vertex $\ \wn\ $.

\subsubsection{\bf The case $N=1$}
There is a short exact sequence
\Beq\label{10: s.e.s. for N=1}
0\lon \overline{\ICH}_n((n+1))^{nice}_{solid} \lon
\overline{\ICH}_n((n+1))^{nice}_1 \lon \overline{\ICH}_n((n+1))^{nice}_{dotted} \lon 0,
\Eeq
where the subcomplex $\overline{\ICH}_n(n+1)^{nice}_{solid}$ is spanned by graphs whose unique edge $e_n$ attached to $\wn$ is {\it solid}. The differential acts trivially on that edge. Hence we can understand the generators of  $\overline{\ICH}_n(n+1)^{nice}_{solid}$ as {\it rooted trees}\, with the root being the solid edge $e_n$ and with legs being of two types, solid and dotted ones. Solid legs are labelled by integers from the set $\{0,1,\dots, n-1\}$ and dotted legs are labelled by integers from
the set $\{1,\dots, n-1\}$. 

\bip

{\sc ``Solid" subcase}. Consider a filtration of  $\overline{\ICH}_n(n+1)^{nice}_{solid}$ by the number of solid edges and let $(E^{solid}_r,\p_r)$ be the associated (regular and bounded) spectral sequence. As $\p_0=\delta_\bu$ and we work with rooted at least trivalent trees, 
the complex $(E^{solid}_0,\p_0)$ is isomorphic to a summand in the complex $\Lie_\infty$.
Hence the next page $E^{solid}_1=H^\bu(E^{solid}_0,\p_0)$ is isomorphic (up to the degree shift) to the
the free Lie superalgebra
$$
E^{solid}_1\simeq \mathfrak{freeLie} \left\langle x_0,x_1,\ldots, x_{n-1}, y_1,\dots, y_{n-1}  \right\rangle, \ \ \ |x_\bu|=0, \ \ |y_\bu|=-1,
$$
generated by  $n$ formal variables $x_0,x_1,\ldots, x_{n-1}$  of cohomological degree zero (corresponding to solid edges  attached to white vertices), and  $n-1$ formal variables $y_1,\dots, y_{n-1}$ of degree
$-1$ (corresponding to dotted edges attached to white vertices). 
The induced differential $\p_1= d^{(2)}$ acts on the formal variables follows as (see (\ref{10: d(2) on dotted edge i}))
$$
\p_1 (x_k)=0,\  \p_1(y_i)=[x_0,x_i],\ \ \ \ \forall k\in\{0,1,\dots, n-1\},\ 
\forall i\in \{1,\dots,n-1\}.
$$
We claim that the spectral sequence $(E^{solid}_r,\p_r)$  degenerates at $r=2$, 
\Beqr
E^{solid}_2 \simeq gr H^\bu(gr \overline{\ICH}_n((n+1))^{nice}_{solid}) 
&\simeq & \frac{\mathfrak{freeLie}(x_0,x_1,\dots,x_n)}{([x_0,x_1],\dots,[x_0,x_n])} \nonumber \\
&\simeq & \mathfrak{freeLie}
\left\langle x_1,\ldots, x_{n-1}  \right\rangle\ \ \oplus \ \  
\K\langle x_0 \rangle \label{10: on E_2 solid}
\Eeqr
Indeed, by the Elimination Theorem (see Chapter 2, \S 9, Proposition 10 and its Corollary in \cite{B}) applied to the decomposition of the set of generators 
$X:=\{x_0,x_1,\ldots, x_{n-1}, y_1,\dots, y_{n-1}\}$ into the disjoint union $x_0\coprod (X\setminus x_0)$, we get a decomposition of $\mathfrak{freeLie}\left\langle x_0,x_1,\ldots, x_{n-1}, y_1,\dots, y_{n-1}  \right\rangle$ into the following direct sum of graded vector spaces,
$$
\mathfrak{freeLie}\left\langle x_0,x_1,\ldots, x_{n-1}, y_1,\dots, y_{n-1}  \right\rangle = 
\K\langle x_0\rangle \oplus \mathfrak{freeLie}\left\langle x_1,\ldots, x_{n-1} \right\rangle \oplus \mathfrak{freeLie}'
$$
where  $\mathfrak{freeLie}'$ is the subspace of the free Lie superalgebra
on the following set of generators,
$$
\left\langle ad_{x_0}^{p_1}\cdot x_1,\ldots,ad_{x_0}^{p_{n-1}} \cdot x_{n-1},ad_{x_0}^{q_1}\cdot y_1,\dots, ad_{x_0}^{p_{n-1}} \cdot y_{n-1}  \right\rangle_{p_1,\dots, p_{n-1}, q_1,\dots, q_{n-1}\geq 0}
$$
which contains at least one symbol with $x_0$ or at least one letter from the alphabet $\{ y_1,\dots, y_{n-1}\}$. Re-denote some of these generators as follows
$$
ad_{x_0}\cdot x_i=: z , \ \ \    ad_{x_0}^{p}\cdot x_i=: ad_{x_0}^{p-1}\cdot z\ \ \ \forall p\geq 1, i\in \{1,..., n-1\}.
$$
In this notation the differential $\p_1$ acts non-trivially only on the summand $\mathfrak{freeLie}'$ as a derivation acting on the $y$-generators as  $\p_1 y_\bu=z_\bu$ 
and trivially on all other generators. We can now define a contracting homotopy
$h$ on $\mathfrak{freeLie}'$ as a derivation acting non-trivially only on the new letter as  $h(z_\bu)=y_\bu$. For any basis element $f$  of $\mathfrak{freeLie}'$ containing $a$ letters $z_\bu$ and $b$ letters $y_\bu$ we have
$$
(\p_1\circ h + h\circ \p_1)f= (a+b)f.
$$
Since $a+b\geq 1$, we conclude that the subcomplex $(\mathfrak{freeLie}',\p_1)$ is acyclic. The claim  (\ref{10: on E_2 solid}) is proven.

\sip

As the root edge of every generator of  the complex $\overline{\ICH}_n((n+1))^{nice}_{solid}$  is solid, the isomorphism  (\ref{10: on E_2 solid}) implies that  the cohomology
of the summand  $gr \overline{\ICH}_n((n+1))^{nice}_{solid}$ is given by the equivalence classes of unrooted  trivalent graphs (modulo the IHX relation along each dotted internal edge)
whose {\it all} \, legs are solid edges. Every such a graph has the cohomological degree +1
so that we conclude that 
$$
H^{\bu \leq 0}(\overline{\ICH}_n((n+1))^{nice}_{solid})=0.
$$

\mip

{\sc ``Dotted" subcase}. We treat this subcase in the full analogy to the previous ``solid" subcase, i.e. we consider a filtration  of the quotient operad $\overline{\ICH}_n(n+1)^{nice}_{dotted}$ by the number of solid edges and study the associated spectral sequence $(E^{dotted}_r,\p_r)$ by identifying its elements with rooted trees,  the major difference being that the root edge is now a {\it dotted edge}\, attached to $\ \wn\ $, not a solid one.  This difference plays no role in the arguments leading us to the main conclusion,
$$
E^{dotted}_2 \simeq gr H^\bu(gr \overline{\ICH}_n((n+1))^{nice}_{dotted}) \simeq \mathfrak{freeLie}
\left\langle x_1,\ldots, x_{n-1}  \right\rangle\ \ \oplus \ \  
\K\langle x_0 \rangle
$$
the element $x_0$ being central in this Lie algebra. 
As the root  edge is now {\it dotted}, the vector space $E_2^{dotted}$ is
now concentrated in the cohomological degree zero (not in degree +1 as in the previous case) 
$$
 H^0(\overline{\ICH}_n((n+1))^{nice}_{dotted})\simeq  \mathfrak{freeLie} \langle T_{1n}, T_{2n}, \ldots, T_{n-1,n} \rangle
\ \ \oplus \ \ \K\langle T_{0n}\rangle \ \ \ \ \forall\ n\geq 2.
$$

Next we shall show  the connecting morphism associated with the short exact sequence (\ref{10: s.e.s. for N=1}),
$$
C: H^0(\overline{\ICH}_n((n+1))^{nice}_{dotted}) \rar H^1(\overline{\ICH}_n((n+1))^{nice}_{solid}),
$$
vanishes.
Indeed the connecting morphism is given by  $d_n^{(2)}$ acting on the unique dotted edge of the graph generators of $H^0(\overline{\ICH}_n((n+1))^{nice}_{dotted})$
as follows
$$
 \Ba{c}\resizebox{12mm}{!}{ \xy
(-7,0)*+{_0}*\frm{o}="0";
(7,0)*+{_n}*\frm{o}="i";
(0,12)*{\bu}="bu";
 \ar @{.} "i";"bu" <0pt> 
\endxy} \Ea  
\lon 
\Ba{c}\resizebox{11mm}{!}{ \xy
 (0,0)*+{_0}*\frm{o}="0";
(12,0)*+{_n}*\frm{o}="i";
(6,12)*{\bu}="bu";
  (6,5)*{\bu}="s";
\ar @{.} "s";"bu" <0pt> 
 \ar @{-} "s";"i" <0pt>
 \ar @{-} "s";"0" <0pt>
\endxy} \Ea
$$
This action, when re-written equivalently in terms of the above identification 
 of $E_2^{solid}$ with $\mathfrak{freeLie}
\left\langle x_1,\ldots, x_{n-1}  \right\rangle\ \oplus  \  
\K\langle x_0 \rangle$, is given by the map:
$$
\Ba{rccc}
C: & \mathfrak{freeLie}
\left\langle x_1,\ldots, x_{n-1}  \right\rangle\ \ \oplus \ \  
\K\langle x_0 \rangle & \lon & \mathfrak{freeLie}
\left\langle x_1,\ldots, x_{n-1}  \right\rangle\ \ \oplus \ \  
\K\langle x_0 \rangle\\
& f & \lon & [x_0,f]=0
\Ea
$$
which vanishes as $x_0$ as is central. Thus we conclude that
$$
H^{\bu \leq -1}(\overline{\ICH}_n((n+1))^{nice}_{1})=0, \ \ 
H^{0}(\overline{\ICH}_n((n+1))^{nice}_{1})=\mathfrak{freeLie}
\left\langle x_1,\ldots, x_{n-1}  \right\rangle\ \ \oplus \ \  
\K\langle x_0 \rangle.
$$

\bip

The proof of Theorem {\ref{9: Theorem on H(ICH_n)}} will be completed once we  show that 
$$
 H^{\bu\leq 0}(\overline{\ICH}_n((n+1))^{nice}_N)=0\ \ \ \ \forall  N\geq 2.
$$

\subsubsection{\bf The case $N\geq 2$}
 Consider a  filtration of $\overline{\ICH}_n((n+1))^{nice}_N$ 
$$
\overline{\ICH}_n((n+1))^{nice}_N\equiv \overline{\ICH}_n((n+1))^{nice}_{0,N}
\supset \overline{\ICH}_n((n+1))^{nice}_{1,N} \supset \dots \supset \overline{\ICH}_n((n+1))^{nice}_{N,N}
$$
where  $\overline{\ICH}_n((n+1))^{nice}_{p,N}$ is the subspace spanned by 
hypergraphs with at least $p$ solid edges attached to $\ \wn\ $.
Let $(E_r, \p_r)$ be the associated spectral sequence. The induced differential in $(E_0, \p_0)$ is
$$
\p_0= d_\bu + \sum_{i=1}^{n-1} d_i^{(2)}.
$$
Arguing as in the ``solid subcase" above  we conclude that the next page $(E_1,\p_1)$ is spanned by unrooted trivalent trees modulo the IHX-relation for each internal edge. The solid legs of such trees with at least one black vertex are labelled by integers from the set 
$\{0,1,2,\dots, n\}$, while the dotted legs (if any) are all labelled by $n$ only. 
  These trees can be equivalently understood 
as formal cyclic words of the form
$$
 (x_n, f_1), \ \ (y_n, f_2),
 $$
where  
$$
f_{1}, f_2 \in \mathfrak{freeLie} \left\langle x_0,x_1,\ldots, x_{n}, y_n\right\rangle/[x_0,x_{i<n}]
$$ 
and the formal scalar bracket $(\ ,\ )$ is assumed to be  graded  symmetric and satisfy the ``Killing form type" identities,
$$
 ([a_i,b_j],f)+ (-1)^{|a_i||b_j|} (b_j, [a_i,f])=0,
$$
for any $a_i,b_j\in  \{ x_0,x_1,\ldots, x_{n}, y_n \}$ and any Lie word $f$ as above. The differential $\p_1$ acts on the letters $y_n$ in these cyclic words as a derivation $d_n^{(2)}: y_n \rar [x_0, x_n]$.
It is, of course, assumed that total number of the letters  $x_n,y_n$ in the above cyclic words $(x_n, f_1)$ or $(y_n, f_1)$   is precisely  $N\geq 2$. Due to the ``Killing form type" identities the decomposition into the sum $(x_n, f_1) + (y_n, f_2)$ is not uniquely defined, say
$$
 (x_n, [y_n,f])=([x_n,y_n],f)= - (y_n, [x_n,f]),
$$
which corresponds to the fact that the trivalent trees are unrooted --- one can choose a solid leg labelled by $n$  as the root $x_n$ (and then immediately get an associated Lie word) or one can choose a dotted leg 
$y_n$ as a root (and get a different Lie word).
\sip

The differential $\p_1$ in $E_1$ preserves the number of solid edges attached to
white vertices $\ \wi\ $ with $i\in \{1,\dots, n-1\}$, i.e.\ it preserves the number of formal variables $x_i$  with $i\in \{1,\dots, n-1\}$ in the above representation of the trivalent trees  in terms of cyclic words. Hence we split the complex $(E_1,\p_1)$ into the direct sum
$$
(E_1,\p_1)= (E_1^{\geq 1},\p_1) \oplus (E_1^{0},\p_1)
$$
where $(E_1^{\geq 1},\p_1)$ is the direct summand spanned by cyclic words which have at least one letter $x_i$ with $i\in \{1,\dots, n-1\}$. Each such cyclic word can be represented
in the form $(x_i, f)$ so that we get a direct sum of complexes
$$
(E_1^{\geq 1},\p_1)=\bigoplus_{i=1}^{n-1}\left(\mathfrak{freeLie} \left\langle x_0,x_i,x_{i+1}\ldots,x_{n-1}, x_{n}, y_{n}  \right\rangle/([x_0,x_{k<n}]), \p_1=d_n^{(2)}\right). 
$$
Using the Elimination Theorem and the contracting homotopy $h$  as the previous case $N=1$ we conclude that each  cohomology classes of each summand in the above complex admits a cycle representative which is a linear combination of Lie words with no single letter $y_n$  present; in terms of graphs it means that the cohomology group $H^{\bu}(E_1^{\geq 1},\p_1)$ is generated by trivalent graphs (modulo the $IHX$-relations) whose {\it all} legs are solid; all such graphs have degree $+1$ so that we conclude that
$$
H^{\bu\leq 0}(E_1^{\geq 1},\p_1)=0.
$$

The second summand $(E_1^{0},\p_1)$ is best understood in terms of trivalent graphs on two white vertices, $\ \wo
 $ and $
  \wn\ $. This complex is spanned by trivalent graphs (modulo the IHX relations) of the form,
  $$
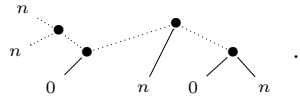
\Ba{c}\resizebox{38mm}{!}{ \xy
(-10,1)*+{_0}="0";
(10,1)*+{_0}="2";
(3,1)*+{_n}="3";
(20,1)*+{_n}="n";
  (-5,6)*{\bu}="bu";
   (7.5,10)*{\bu}="r";
     (15.5,6)*{\bu}="rr";
(-9,9)*{\bu}="s";
(-14,12)*+{_n}="f1";
(-15,6)*+{_n}="f2";
 \ar @{-} "0";"bu" <0pt>
   \ar @{.} "bu";"r" <0pt>
    \ar @{-} "rr";"2" <0pt>
     \ar @{-} "rr";"n" <0pt>
  \ar @{.} "rr";"r" <0pt>
 \ar @{-} "r";"3" <0pt>
 \ar @{.} "s";"bu" <0pt>
 \ar @{.} "s";"f1" <0pt>
 \ar @{.} "s";"f2" <0pt>
\endxy} \Ea.
$$
These graph can have dotted edges  attached only to $\wn$ white its solid edges can be attached
to $\wo$ and $\wn$.
The differential  $\p_1=d_n^{(2)}$ acts on dotted edges as follows
as follows
$$
 \Ba{c}\resizebox{5mm}{!}
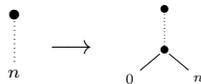
{ \xy
(0,0)*+{_n}="i";
(0,8)*{\bu}="bu";
 \ar @{.} "i";"bu" <0pt> 
\endxy} \Ea  
\lon 
\Ba{c}\resizebox{13mm}{!}{ \xy
 (0,0)*+{_0}="0";
(12,0)*+{_n}="i";
(6,12)*{\bu}="bu";
  (6,5)*{\bu}="s";
\ar @{.} "s";"bu" <0pt> 
 \ar @{-} "s";"i" <0pt>
 \ar @{-} "s";"0" <0pt>
\endxy} \Ea
$$
Hence the complex $(E_1^{0},\p_1)$ is acyclic.

\sip

The proof of Theorem {\ref{9: Theorem on H(ICH_n)}} is completed. $\Box$

\sip

\bip

\def\cprime{$'$}

\end{document}